\documentclass[12pt, reqno]{amsart}
\usepackage[ruled,vlined]{algorithm2e}
\usepackage{fancyhdr,graphicx,amsmath,amssymb,amsthm}
\usepackage[left=1.3in,right=1.3in,bottom=1.0in]{geometry}
\usepackage{tcolorbox, enumitem, verbatim}
\usepackage{yfonts}
\usepackage{xcolor}
\usepackage{hyperref}
\hypersetup{
    colorlinks=true,   
    linkcolor=blue,    
    citecolor=red,     
    urlcolor=magenta,  
    filecolor=cyan     
}
\usepackage{tikz}
\usepackage[utf8]{inputenc}
\usepackage{mathtools}
\usepackage{faktor}
\usepackage{quiver}
\usepackage{tikz-cd}
\usepackage{bbm}
\usepackage{faktor}
\usepackage{mathabx}
\usepackage{thmtools}
\usepackage{xcolor}
\usepackage{mathrsfs}
\usepackage[nameinlink,capitalise,noabbrev]{cleveref}
\usepackage[backend=biber,style=alphabetic,maxbibnames=4,doi=false,isbn=false,url=false]{biblatex}
\renewbibmacro*{in:}{}

\AtEveryBibitem{\clearlist{location}}
\usepackage{bm}
\usepackage{stmaryrd}
\usepackage[OT1]{fontenc}

\theoremstyle{plain}
\newtheorem{theorem}{Theorem}[section]
\newtheorem{coroll}[theorem]{Corollary}
\newtheorem{lemma}[theorem]{Lemma}
\newtheorem{prop}[theorem]{Proposition}

\newtheorem*{claim}{Claim}

\theoremstyle{definition}

\newtheorem{definition}[theorem]{Definition}

\newtheorem*{definition*}{Definition}

\theoremstyle{remark}
\newtheorem{remark}[theorem]{Remark}

\newenvironment{claimproof}[1][Proof of Claim]{%
  \begin{proof}[#1]%
}{%
  \end{proof}%
}

\newcommand{\CC}{\mathbb{C}}
\newcommand{\AAA}{\mathbb{A}}
\newcommand{\ZZ}{\mathbb{Z}}

\newcommand{\PP}{\mathbb{P}}
\newcommand{\QQ}{\mathbb{Q}}

\newcommand{\GG}{\mathbb{G}}

\newcommand{\HH}{\mathbb{H}}

\newcommand{\DR}{\mathrm{DR}}

\newcommand{\univ}{\mathrm{univ}}
\newcommand{\Hod}{\mathrm{Hod}}

\newcommand{\rank}{\operatorname{rank}}
\newcommand{\im}{\operatorname{im}}

\newcommand{\sh}[1]{\mathcal{#1}}
\newcommand{\bb}[1]{\mathbb{#1}}
\newcommand{\scr}[1]{\mathscr{#1}}
\newcommand{\Rm}[1]{\mathrm{#1}}
\newcommand{\Sf}[1]{\mathsf{#1}}
\newcommand{\Bf}[1]{\mathbf{#1}}

\newcommand{\Ide}[1]{\mathfrak{#1}}
\newcommand{\der}{\mathrm{d}}

\newcommand{\ol}[1]{\overline{#1}}

\newcommand{\Spec}{\operatorname{Spec}}

\newcommand{\Sum}[1]{\underset{#1}{\sum}}

\newcommand{\Prod}[1]{\underset{#1}{\prod}}
\newcommand{\tab}{\hspace*{\parindent}}
\newcommand{\red}{\mathrm{red}}

\DeclareMathOperator{\Hom}{Hom}%

\newcommand{\genlegendre}[4]{%
  \genfrac{(}{)}{}{#1}{#3}{#4}%
  \if\relax\detokenize{#2}\relax\else_{\!#2}\fi
}

\DeclarePairedDelimiter{\set}{\{}{\}}
\DeclarePairedDelimiter{\abs}{\lvert}{\rvert}

\DeclarePairedDelimiter{\brac}{(}{)}
\DeclarePairedDelimiter{\sqbrac}{[}{]}
\DeclarePairedDelimiter{\basis}{\langle}{\rangle}
\DeclareMathAlphabet{\pazocal}{OMS}{zplm}{m}{n}

\newcommand{\Comment}[1]{}

\def\XXint#1#2#3{{\setbox0=\hbox{$#1{#2#3}{\int}$}
     \vcenter{\hbox{$#2#3$}}\kern-.5\wd0}}

\allowdisplaybreaks

\makeatletter
\def\l@subsection{\@tocline{2}{0pt}{2.5pc}{5pc}{}}
\makeatother

\numberwithin{equation}{section}

\DeclareFieldFormat[article,inbook,incollection,inproceedings,patent,thesis,unpublished]{title}{#1\isdot}

\title[Rank two cohomology jump loci]{Rank two cohomology jump loci on projective varieties}
\author{Quoc Anh Tran}
\address{Department of Mathematics, University of Chicago, Chicago, Illinois, USA}
\email{quocanh@uchicago.edu}
\urladdr{https://quocanh-tran.github.io/}
\date{\today}

\begin{document}
\begin{abstract}
  Let $X$ be a smooth connected complex projective variety. We give a structure theorem for irreducible components of the reduced rank two cohomology jump loci on $X$. More precisely, for such a component $Z$, if the general point of $Z$ is projectively virtually abelian or projectively rigid, then either $Z$ is a symmetric square of an algebraic torus, or a quotient of a torus by an involution, or its normalization $Z^{\nu}$ is a torus. If the general point is neither, then either $Z^{\nu}$ is a torus, or an \'{e}tale quotient of $Q^{\nu} \times T$ where $Q$ is a fixed-determinant, fixed-inertia component of the character variety of some DM-curve. 
\end{abstract}
\maketitle
\vspace{-1cm}
\tableofcontents

\section{Introduction}
\label{sec:intro}
Let $X$ be a smooth projective variety, and 
\begin{align*}
  R_B(X, x, r) &\coloneq \Hom(\pi_1(X, x), \Rm{GL}_r)\\
  \sh{M}_B(X, r) &\coloneq [R_B(X, x, r)/\Rm{GL}_r]\\
  M_B(X, r) &\coloneq R_B(X, x, r)//\Rm{GL}_r
\end{align*}
be, respectively, the moduli of framed local systems, the quotient stack, and the GIT quotient/good moduli space for the stack. We have the cohomology jump loci 
\[\Sigma^i_k(X, x, r) \coloneq \set*{\bb{V} \in R_B(X, x, r)\ \Big\vert\ \dim H^i(X, \bb{V}) \geq k}\]
which can be equipped with a determinental subscheme structure. This loci is $\Rm{GL}_r$-invariant hence descends to a substack $\Sigma^i_k(X, r) \subset \sh{M}_B(X, r)$. Denote $\ol{\Sigma}^i_k(X, r)$ its reduced image in $M_B(X, r)$. \\
\tab The study of these loci originated from Green-Lazarsfeld's generic vanishing theorem \cite{GreenLazarsfeld1987} for topologically trivial holomorphic line bundles on compact K\"{a}hler $X$. Their results concern jump loci in $\Rm{Pic}^0(X)$. They later proved that any irreducible jump component must be a translated subtorus in \cite{GreenLazarsfeld1991}. The corresponding Betti version is about rank one jump loci $\Sigma^i_k(X, 1)$ in
\[R_B(X, x, 1) \simeq M_B(X, 1)\]
which is a finite disjoint union of algebraic tori. Simpson \cite{simpson1993subspaces} improved on Green-Lazarsfeld's result by showing that any irreducible component of $\Sigma^i_k(X, 1)$ is a torsion-translated subtorus.\\
\tab Since then, there have been several generalizations. Arapura \cite{Arapura1997} studied jump loci $\Sigma^i_{1}(X, 1)$ on quasi-compact K\"{a}hler manifolds and gave a description of irreducible components via maps to curves. Budur-Wang \cite{BudurWang2020} proved the translated subtorus theorem for $\Sigma^i_k(X, 1)$ for all $i, k$ for quasi-compact K\"{a}hler manifolds; they also recovered the torsion statement of Simpson to smooth quasi-projective varieties \cite{BudurWang2015QuasiProjective,budur2020absolutesets}. Esnault-Kerz developed an arithmetic approach to these results \cite{EsnaultKerz2020}.\\
\tab In higher rank, the local picture has been studied extensively. Dimca-Papadima \cite{DimcaPapadima2014} related the germs of nonabelian jump loci near the trivial representation to commutative differential graded algebras (via Sullivan's construction). Budur-Wang \cite{budur2015cohomology} formulated the deformation theory with cohomological constraints using differential graded Lie algebra pairs. Specifically in rank two, Papadima-Suciu \cite{PapadimaSuciu2019,PapadimaSuciu2020} obtained decomposition and pullback (from maps to curves) results for the jump loci germs near the trivial representation. \\
\tab However, the global geometry of higher rank jump loci are largely unknown. Our main result is a structure theorem for rank two jump loci components, up to normalization. 
\begin{theorem}
  Let $Z$ be an irreducible component of $\overline\Sigma^i_k(X,2)$.
\begin{enumerate}
\item In the reducible case, $Z\simeq\Rm{Sym}^2(N)$ for a translated torus $N$, or $Z^\nu$ is a product of two translated tori. In the irreducible projectively dihedral case, $Z^\nu$ is a torus or $Z$ is a quotient of a translated torus by an involution. In the projectively finite and projectively dense rigid cases, $Z^\nu$ is a torus.
\item In the projectively dense nonrigid case, either $Z^\nu$ is a torus or an \'{e}tale quotient
\[
Z^\nu\simeq(Q^\nu\times T)/H,
\]
where $T$ is an algebraic torus, $Q$ is a component of a fixed-determinant, fixed-inertia locus $M_B(\sh{Y}, 2)$ where $\sh{Y}$ a proper hyperbolic Deligne-Mumford curve. 
\end{enumerate}
Here $Z^\nu$ and $Q^\nu$ denote normalizations.
\end{theorem}

\subsection{Hodge structures on jump loci}
We first obtain several Hodge-theoretic results concerning Hodge subvarieties of $R_B(X, x, r)$, which should be of independent interest. 
\begin{theorem}[{\cref{thm:jump=Hodge}}]
  For smooth projective $X$, the jump locus $\Sigma_k^i(X, x, r) \subseteq R_B(X, x, r)$ is a Hodge subvariety (see \cref{def:Hodge-substack}). Furthermore, it is closed under semisimplification, and $\Sigma_k^i(X, x, r)_{\red}$ is a Hodge subvariety which is absolute $\QQ$-closed. \\
  \tab Every irreducible component of $\Sigma_k^i(X, x, r)_{\red}$ is a Hodge subvariety and is absolute $\ol{\QQ}-$closed. 
\end{theorem}
At $\bb{V} \in \Sigma^i_k(X, x, r)$ underlying a complex variation of mixed Hodge structures, this result gives a pro-MHS-algebra structure on $\widehat{\sh{O}}_{\Sigma^i_k(X, x, r), \bb{V}}$. Hence we have an alternative approach to the result of Lef\`{e}vre \cite{lefevre2021mixed}. We do not know if our Hodge structures agree. 

\begin{theorem}[{\cref{thm:GIT-Hodge}, \cref{thm:norm-Hodge}}]
  Let $Z \subseteq M_B(X, r)$ be the scheme-theoretic image of a Hodge subvariety in $R_B(X, x, r)$. Let $\bb{V} \in Z$ be a local system underlying a complex variation of Hodge structures, then 
  \[\widehat{\sh{O}}_{Z, \bb{V}}, \tab \widehat{\sh{O}}_{Z^{\nu}, z}\]
  both underlie pro-MHS-algebras, where $Z^{\nu} \to Z$ is the normalization map with $z \in Z^{\nu}$ lying over $\bb{V} \in Z$. Furthermore, the morphism 
  \[\widehat{\sh{O}}_{Z, \bb{V}} \to \widehat{\sh{O}}_{Z^{\nu}, z}\]
  is a morphism of pro-MHS-algebras. 
\end{theorem}
In particular, this gives a pro-MHS-algebra structure on $\widehat{\sh{O}}_{M_B(X, r), \bb{V}}$ for $\bb{V}$ underlying a complex variation of Hodge structures, thus recovering a result of Eyssidieux-Simpson \cite{EyssidieuxSimpson2011}. 
\begin{remark}
  We opted for the projective setting due to later classification result, and ease of exposition. However, most of section 3 (which contains these Hodge-theoretic results) should translate to the case where $X$ is smooth and quasi-projective, with appropriate setup. 
\end{remark}

\subsection{Outline of the classification of jump components}
Here we specialize to the case $r = 2$. Let $A = M_B(X, 1)$. As consequence of a point-wise classification theorem \cite{nguyen2008algebraic}, we can label the components based on their general point: 
\begin{prop}[{\cref{lem:jump-labels}}]
  Let $Z \subseteq \ol{\Sigma}^i_k(X, 2)$ be an irreducible component. Consider a general point $\rho \in Z$, with $\bar{\rho}: \pi_1(X, x) \to \Rm{PGL}_2(\CC)$ the projective quotient. Then $\rho$ is exactly one of the following cases:
  \begin{description}
    \item[(A)] $\rho$ is reducible, i.e., $\rho \simeq \chi_1 \oplus \chi_2$;
    \item[(B)] $\rho$ is irreducible and projectively dihedral, i.e. up to conjugation, $G_{\ol{\rho}} \subseteq D_{\infty}$;
    \item[(C)] $G_{\ol{\rho}} \simeq \ol{\rho}(\pi_1(X, x)) \simeq A_4, S_4$ or $A_5$;
    \item[(D)] $\ol{\rho}$ has Zariski dense image and is rigid;
    \item[(E)] $\ol{\rho}$ has Zariski dense image and is non-rigid.
  \end{description}
\end{prop}
Next, we have a structure theorem for the first 4 cases:
\begin{theorem}
  \leavevmode
  \begin{description}
    \item[(A)] (\cref{thm:reducible-jump}) either $Z \simeq \Rm{Sym}^2(N)$ or $Z^{\nu} \simeq N_1 \times N_2$ with $N, N_1, N_2 \subseteq A$ translated subtori of $A$;
    \item[(B)] (\cref{thm:dihedral-jump}) either $Z^{\nu} \simeq N$ or $Z \simeq N/\basis{\tau^*}$ with $N$ a translated subtorus of $M_B(X^{\epsilon}, 1)$ and $\tau^*$ an involution on $N$; 
    \item[(C)] (\cref{thm:finite-jump}) $Z^{\nu}$ is an algebraic torus;
    \item[(D)] (\cref{thm:rigid-jump}) $Z^{\nu}$ is an algebraic torus. 
  \end{description}
\end{theorem}
For the last case, which we will call nonabelian jump components, up to taking a finite cover, we can consider a point of $Z$ as a pair $(\bb{V}, \bb{L})$ where $\bb{V}$ has fixed determinant and $\bb{L}$ is a rank one local system. More precisely, there is a torsion lift $\epsilon \in A$ of some $\bar{\epsilon} \in A/A^2$ such that $Z \subseteq M_B(X, 2, \det \in \epsilon A^2)$. Consider the diagram: 
\[\begin{tikzcd}
	{M_B(X, 2, \det = \epsilon) \times A} && {M_B(X, 2, \det \in \epsilon A^2)} \\
	\\
	A && A
	\arrow["{(\bb{V}, \bb{L}) \mapsto \bb{V} \otimes \bb{L}}","{\Phi_{\epsilon}}"', from=1-1, to=1-3]
	\arrow["{\Rm{pr}_2}"', from=1-1, to=3-1]
	\arrow["\det", from=1-3, to=3-3]
	\arrow["{\bb{L} \mapsto \epsilon \otimes \bb{L}^{\otimes 2}}", from=3-1, to=3-3]
  \end{tikzcd}\]
  and let $C$ be an irreducible component of $\Phi_{\epsilon}^{-1}(Z)$. We can recover $Z$, in a sense, as follows (for more details, see \cref{sec:jump-product})
  \begin{prop}[{\cref{lem:normalized-det-sheet}}]
    $Z^{\nu} \simeq C^{\nu}/H_C$ is an \'{e}tale quotient of $C^{\nu}$.
  \end{prop}
  Denote the projections:
  \[\begin{tikzcd}
	{} & {C \subseteq \Phi_{\epsilon}^{-1}(Z) \subseteq M_B(X, 2, \det = \epsilon) \times A} && A \\
	& {M_B(X, 2, \det = \epsilon)}
	\arrow["{\Rm{pr}_1}"', from=1-2, to=1-4]
	\arrow["{\Rm{pr}_1}", from=1-2, to=2-2]
\end{tikzcd}\]
  with images:
  \[D = \ol{\Rm{pr}_1(C)}_{\red} \subseteq M_B(X, 2, \det = \epsilon), \tab T = \ol{\Rm{pr}_2(C)}_{\red} \subseteq A\]
  \leavevmode\vspace{-\baselineskip}
\begin{prop}[{\cref{lem:GIT-proj=torus}}]
  $T$ is a torsion-translated subtorus of $A$.
\end{prop}
As a consequence of \cite{corlette2008classification} we have a factorization of $D$ through a DM-curve:
\begin{prop}[{\cref{thm:Shafarevich-map}}]
  There is a proper hyperbolic DM curve $\sh{Y} = [Y/G]$ and a morphism $s: X \to \sh{Y}$ and an irreducible component $Q$ of a fixed-determinant, fixed-inertia locus $M_B(\sh{Y}, 2, \det = \delta, \bm{\tau})$ such that
  \begin{enumerate}
    \item $s_*: \pi_1(X, x) \to \pi_1(\sh{Y}, \bar{y})$ is a surjection.
    \item with the closed embedding $s^*: M_B(\sh{Y}, 2) \hookrightarrow M_B(X, 2)$, $D \subseteq s^*(Q)$. 
  \end{enumerate}
\end{prop}
Via the closed embedding, we can view $C$ as a subvariety of $Q \times T \subseteq M_B(\sh{Y}, 2) \times A$. More precisely, 
\[C \subseteq \set*{(\bb{W}, \bb{L}) \in Q \times T\ \big\vert\ \dim H^i(X, s^*\bb{W} \otimes \bb{L}) \geq k}\]
\tab On the other hand, 
\[\Rm{R}\Gamma(X, s^*\bb{W} \otimes \bb{L}) \simeq \Rm{R}\Gamma(\sh{Y}, \bb{W} \otimes \Rm{R}s_*\bb{L})\]
hence we can, instead, study a twisted jump loci on the DM-curve $\sh{Y}$. For a fixed $\bb{L}$, we understand when this cohomology jumps. 
\begin{prop}[{\cref{thm:fin-or-prod}}]
  Either $C = Q \times T, C = \set*{\bb{V}} \times T$, or $\Rm{pr}_2: C \to T$ is generically finite. 
\end{prop}
For the case where $\Rm{pr}_2: C \to T$ is generically finite, using an invariant Grassmannian construction, we construct a finite scheme over $T$ (\cref{thm:finite-closure}) which contains a Zariski open subset of $C$. As a consequence, 
\begin{prop}[{\cref{thm:not-gen-fin}}]
  If $\Rm{pr}_2: C \to T$ is generically finite, then it is finite and surjective. 
\end{prop}
Finally, suppose that $\Rm{pr}_2:C \to T$ is finite, and consider the normalized map $\Rm{pr}_2^{\nu}: C^{\nu} \to T$. If this is not \'{e}tale, we can produce a unitary $\bb{L} \in T$ with $c \in (\Rm{pr}_2^{\nu})^{-1}(\bb{L})$ such that $C^{\nu}$ is smooth at $c$, $\Rm{pr}_2^{\nu}$ is ramified at $c$, and 
\[\widehat{\sh{O}}_{T, \bb{L}} \to \widehat{\sh{O}}_{C^{\nu}, c}\]
is an injective morphism of pro-MHS-algebras. This results in a contradiction (\cref{thm:etaleness}). Putting everything together we get: 
\begin{theorem}
  Either $C = Q \times T$, $C = \set*{\bb{V}} \times T$, or $\Rm{pr}_2^{\nu}: C^{\nu} \to T$ is \'{e}tale. 
\end{theorem}
\begin{coroll}
  If the general point of $Z$ is of type $\Bf{(E)}$, then either $Z$ is an algebraic torus, $Z^{\nu}$ is an algebraic torus, or $Z^{\nu}$ is an \'{e}tale quotient of $Q^{\nu} \times T$.  
\end{coroll}

\subsection{Acknowledgement} I would like to thank my advisor, Benjamin Bakker, for suggesting the problem of putting Hodge structures on jump loci which became the basis of this project. I am grateful for numerous conversations and suggestions, as well as his careful reading of previous versions. I would also like to thank Daniel Litt for his interest in the result. \\
\\
\textbf{AI discloure:} No part of this paper was AI generated. OpenAI's ChatGPT Plus was used for literature searches, and later the Pro version of ChatGPT Astra was used to audit the paper. Along with minor issues, it pointed out previous errors with \cref{thm:GIT-Hodge} and \cref{lem:end-Hodge}. The fixes are the author's. 

\section{Preliminaries}
\label{sec:prelim}
We work with schemes over $\CC$. For a scheme $S$, let $S_{\red}$ be the reduced scheme, and $S^{\nu}$ be the normalization of $S$. 
\subsection{Jump ideals}
\label{sec:jump-ideal}
Let $A$ be a noetherian ring, and $K^\bullet$ a bounded above complex of $A-$modules such that $H^i(K^\bullet)$ is a finitely generated $A$-module for all $i$. By \cite[Lemma III.12.3]{hartshorne1977algebraic}, there exists a quasi-isomorphism $\phi: L^\bullet \to K^\bullet$ where $(L^\bullet, \der^\bullet)$ is a bounded above complex of finitely generated free $A$-modules. 

\begin{definition}
  \label{def:jump-ideal}
  The cohomological jump ideals of $K^\bullet$ are defined to be 
  \[J^i_k(K^\bullet) \coloneq I_{\rank_A(L^i) - k + 1}(\der^{i - 1} \oplus \der^i) \subset A\]
  where $I$ denotes the determinantal ideal generated by $(\rank_A(L^i) - k + 1)^2$-minors. This definition does not depend on the choice of $L^\bullet$ \cite[Definition-Proposition 2.2]{budur2015cohomology}. We also have that, for a base change $A \to B$,
  \begin{align}
    \label{eq:base-change-jump-ideal}
    J^i_k(K^\bullet)B &= J^i_k(K^\bullet \otimes_{A}^{\Rm{L}} B)\\
    V(J^i_k(K^\bullet)) &= \set*{s \in \Spec A\ \big\vert\ \dim_{\kappa(s)}H^i(K^\bullet \otimes_A^{\Rm{L}} \kappa(s)) \geq k} 
  \end{align}
\end{definition}
Let $S = \Spec A \to R_B(X, x, r) = \Hom(\pi_1(X, x), \Rm{GL}_r)$ be a classifying map. This is equivalent to a representation 
\[\rho_S: \pi_1(X, x) \to \Rm{GL}_r(A)\]
which we might think of as a $A[\pi_1(X, x)]$-module $V_S$. Let the corresponding local system of free $A$-modules on $X$ be denoted $\bb{V}_S$. Here we have an explicit free model for $\Rm{R}\Gamma(X, \bb{V}_S)$. Fix a finite CW complex structure for $X$, we can form the cellular complex $C_\bullet(\widetilde{X}, \ZZ)$ on the universal cover $\widetilde{X} \to X$; this is a complex of free $\ZZ[\pi_1(X, x)]$-modules. Then 
\[C^\bullet(X, \bb{V}_S) \coloneq \Hom_{\ZZ[\pi_1(X, x)]}\brac*{C_\bullet(\widetilde{X}, \ZZ), V_S}\]
is a complex of free $A$-modules. Cellular cohomology with local system coefficients identifies this with $\Rm{R}\Gamma(X, \bb{V}_S)$. Furthermore, for $f: \Spec B \to S$, we have that 
\begin{equation}
  C^\bullet(X, \bb{V}_S \otimes_{A_X} B_X) = C^\bullet(X, \bb{V}_S) \otimes_A B \label{eq:cellular-base-change}
\end{equation}
with $A_X, B_X$ being the constant local systems on $X$. Recall that $R_B(X, x, r)$ is affine, and denote the universal local system, which is a local system of free $\Gamma(R_B(X, x, r), \sh{O}_{R_B(X, x, r)})$-modules on $X$, to be $\bb{V}^{\univ}$ (which is just $\bb{V}_S$ for $S = R_B(X, x, r)$). We can then give a scheme structure to 
\[\Sigma^i_k(X, x, r) = \set*{\bb{V} \in R_B(X, x, r)\ \big\vert\ \dim H^i(X, \bb{V}) \geq k}\]
as $\Sigma^i_k(X, x, r) = V\brac*{J^i_k\brac*{\Rm{R}\Gamma(X, \bb{V}^{\univ})}}$. 
\subsection{Variations of complex mixed Hodge and twistor structures}
\label{sec:variations}
\begin{definition}
  \label{def:pHS}
  A complex pure Hodge structure ($\CC$-HS) of weight $n$ is a finite-dimensional $\CC$-vector space $V$ with $n-$opposed decreasing filtrations $F^\bullet, F'^{\bullet}$, i.e.,
  \[\Rm{gr}_F^p \Rm{gr}_{F'}^q V = 0 \tab p + q \neq n, \tab V = \bigoplus_{p + q = n} V^{p, q}, \tab V^{p, q} = F^p \cap F'^{q}\]
  \tab Morphisms are filtered morphisms. A polarization of $V$ is a hermitian form $h$ such that $V = \bigoplus_{p + q = n}V^{p, q}$ is $h$-orthogonal, and $\bigoplus_{p + q = n} (-1)^ph\vert_{V^{p, q}}$ is positive definite. 
\end{definition}
\begin{definition}
  \label{def:MHS}
  A complex mixed Hodge structure ($\CC$-MHS) is $(V, W_\bullet, F^\bullet, F'^{\bullet})$, with $W_\bullet$ an increasing (weight) filtration, such that for any $n$, the filtrations induced by $F^\bullet, F'^{\bullet}$ on $\Rm{gr}^W_n V$ define a pure Hodge structure of weight $n$. \\
  \tab Morphisms are those compatible with all three filtrations. A graded polarization is a polarization $h_n$ on $\Rm{gr}^W_n V$ for each $n$.   
\end{definition}

\begin{definition}
  \label{def:MTS}
  A complex pure twistor structure ($\CC$-TS) of weight $n$ is a vector bundle $\sh{V} \simeq \sh{O}_{\PP^1}(n)^{\oplus r}$. A complex mixed twistor structure ($\CC$-MTS) is a pair $(\sh{V}, W_\bullet)$ where $\sh{V}$ is a vector bundle on $\PP^1$, $W_\bullet \sh{V}$ is a locally split increasing (weight) filtration by subbundles satisfying 
  \[\Rm{gr}^W_n \sh{V} \simeq \sh{O}_{\PP^1}(n)^{\oplus r_n}\]
  \tab Morphisms of complex mixed twistor structures are $\sh{O}_{\PP^1}-$linear morphisms which preserve the weight filtrations. 
\end{definition}
The complex mixed Hodge and twistor structures form categories $\Sf{MHS}_{\CC}$ and $\Sf{MTS}_{\CC}$. These are $\CC$-linear abelian categories. We have fiber functors $\omega$, with image $\CC$-vector spaces, where $\omega(V, W_\bullet, F^\bullet, F'^{\bullet}) = V$, and $\omega(\sh{V}, W_\bullet) = \sh{V}\vert_{\lambda = 1}$. 

\begin{theorem}[{\cite{Deligne1971}, \cite[Proposition 1.10]{Sabbah2018}}]
  \label{thm:tannakian}
  $(\Sf{MHS}_{\CC}, \omega)$ and $(\Sf{MTS}_{\CC}, \omega)$ are neutral Tannakian categories. 
\end{theorem}

\begin{definition}
  \label{def:VMHS}
  A graded polarized complex variation of mixed Hodge structures ($\CC$-VMHS) on $X$ is $(\bb{V}, W_\bullet, F^\bullet, F'^{\bullet})$ where 
  \begin{itemize}
    \item $\bb{V}$ is a complex local system on $X$ and $W_\bullet$ is a flat increasing (weight) filtration;
    \item $F^\bullet, F'^{\bullet}$ are locally split decreasing filtrations
      \[F^p \subset \sh{O}_X \otimes_{\CC_X} \bb{V}, \tab F'^q \subset \overline{\sh{O}}_X \otimes_{\CC_X} \bb{V}\]
    by holomorphic and anti-holomorphic subbundles respectively;
    \item for the natural flat connections $\nabla, \overline{\nabla}$ coming from the exterior derivative, these must satisfy Griffiths transversality 
    \[\nabla F^p \subset F^{p - 1} \otimes \Omega^1_X, \tab \overline{\nabla} F'^{q} \subset F'^{q - 1} \otimes \Omega^1_X\]
    \item there exists a flat hermitian form $h_n$ on each $\Rm{gr}^W_n \bb{V}$ such that for each $x \in X$, the restriction of the tuple to $x$ gives a graded polarized complex mixed hodge structure. 
  \end{itemize}
\end{definition}
For the twistor analogue, we first review the theory of harmonic bundle developed by Corlette and Simpson in \cite{corlette1988flat,simpson1992higgs}. 

\begin{definition}
  Let $\sh{E}$ be a complex $\sh{C}^\infty$-vector bundle on a projective complex manifold $X$ equipped with a flat $C^\infty$-connection $D$. For $h$ a hermitian metric on $\sh{E}$, we have a decomposition
  \[D = D_h + \Psi\]
  where $D_h$ is unitary with respect to $h$, and $\Psi$ is self-adjoint for $h$. Take a further type decomposition into $(1, 0)$ and $(0, 1)$ components
  \[D_h = \partial_h + \bar{\partial}_h, \tab \Psi = \theta + \theta^*\]
  \tab We say that $h$ is pluriharmonic if $(\theta + \bar{\partial}_h)^2 = 0$. A harmonic bundle $(\sh{E}, D, h)$, equivalently $(\sh{E}, \bar{\partial}_h, \theta, h)$, is then a smooth complex vector bundle $\sh{E}$ with a flat smooth connection $D$ and a pluriharmonic metric $h$. 
\end{definition}

Let $\scr{A}_X \coloneq \sh{C}^\infty_X \boxtimes \sh{O}_{\PP^1}$ be the sheaf of $C^\infty$-functions on $X_{\PP^1} \coloneq X \times \PP^1$ which are holomorphic in the $\PP^1$-direction, and choose generating sections $x_0, x_\infty \in H^0(\PP^1, \sh{O}_{\PP^1}(1))$ vanishing at $0$ and $\infty$ respectively. Then we get a derivative 
\[x_0 \partial_{X_{\PP^1}/\PP^1} + x_\infty \bar{\partial}_{X_{\PP^1}/\PP^1}: \scr{A}_X \to \scr{A}_X \otimes \Omega^1_{X_{\PP^1}/\PP^1}(1)\]
with $\Omega^1_{X_{\PP^1}/\PP^1}(1) \coloneq \Omega_X^1 \boxtimes \sh{O}_{\PP^1}(1)$, where $\Omega_X^1$ is the sheaf of $C^\infty$-one-forms on $X$. Let $\brac*{\sh{E}, \bar{\partial}_h, \theta, h}$ be a harmonic bundle on $X$, and define $\sh{E}_{\PP^1} \coloneq \sh{E} \boxtimes \sh{O}_{\PP^1}$. Then we have a $x_0 \partial_{X_{\PP^1}/\PP^1} + x_\infty \bar{\partial}_{X_{\PP^1}/\PP^1}$-connection 
\[\sh{D} \coloneq x_0(\partial_h + \theta^*) + x_\infty(\theta + \bar{\partial}_h): \sh{E}_{\PP^1} \to \sh{E}_{\PP^1} \otimes \Omega^1_{X_{\PP^1}/\PP^1}(1), \tab \sh{D}^2 = 0\]
where we abuse notation so that $\partial_h + \theta^*, \theta + \bar{\partial}_h$ really mean their natural extensions to $\sh{E}_{\PP^1}$. 
\begin{definition}
  \label{def:VTS}
  The resulting $\brac*{\sh{E}_{\PP^1}, \sh{D}}$ is the variation of pure twistor structures ($\CC-$VTS) of weight 0 associated to the harmonic bundle $\brac*{\sh{E}, \bar{\partial}_h, \theta, h}$. We get the weight $k$ $\CC$-VTS associated to $\brac*{\sh{E}, \bar{\partial}_h, \theta, h}$ by twisting $\sh{E}_{\PP^1, k} \coloneq \sh{E} \boxtimes \sh{O}_{\PP^1}(k)$. 
\end{definition}
\begin{remark}
  Every semisimple complex local system on $X$ underlies a harmonic bundle, hence a variation of pure twistor structures. 
\end{remark}

\begin{definition}
  \label{def:VMTS}
  A graded-polarized variation of mixed twistor structures ($\CC$-VMTS) on $X$ is a triple $\brac*{\scr{V}, W_\bullet, \sh{D}}$ such that 
  \begin{itemize}
    \item $\scr{V}$ is an $\mathscr{A}_X$-module on $X_{\PP^1}$;
    \item $\sh{D}: \scr{V} \to \scr{V} \otimes \Omega^1_{X_{\PP^1}/\PP^1}(1)$ is a $x_0 \partial_{X_{\PP^1}/\PP^1} + x_\infty \bar{\partial}_{X_{\PP^1}/\PP^1}$-connection with $\sh{D}^2 = 0$;
    \item $W_\bullet \scr{V}$ is a $\sh{D}-$flat increasing weight filtration of $\sh{E}_{\PP^1}$;
    \item each $\brac*{\Rm{gr}^W_k \scr{V}, \Rm{gr}_k^W \sh{D}}$ has a pluriharmonic hermitian metric $h_k$ making it a $\CC-$VTS of weight $k$ associated to a harmonic bundle. 
  \end{itemize}
  Morphisms are $\scr{A}_X$-linear morphisms which preserve $W_\bullet$ and commutes with $\sh{D}$. 
\end{definition}
\begin{remark}
  \label{rmk:VMTS=lambda}
  On an affine chart $\AAA^1_{\lambda} = \PP^1 - \set*{\infty}$, let $D_{\lambda} = x^{-1}_{\infty} \sh{D}\vert_{X \times \set*{\lambda}}$. The data of a $\CC$-VMTS is equivalently a family of integrable $\lambda$-connections 
  \[\nabla_\lambda(fs) = \lambda \der f \otimes s + f \nabla_{\lambda}(s), \tab \nabla_{\lambda}^2 = 0\]
  where $\nabla_{\lambda} = D^{0, 1}_{\lambda}$ on the holomorphic vector bundle $(\scr{V}\vert_{\lambda}, D^{0, 1}_{\lambda})$. The fibers $\lambda = 1$ and $\lambda = 0$ give a flat bundle, hence an underlying local system by Riemann-Hilbert correspondence, and a Higgs bundle respectively. 
\end{remark}

\subsection{Variations with artinian coefficients}
\label{sec:A-variations}
Let $\Sf{MS}_{\CC}$ be either $\Sf{MHS}_{\CC}$ or $\Sf{MTS}_{\CC}$, and $\mathbf{T}(0)$ be the tensor unit in $\Sf{MS}_{\CC}$. 
\begin{definition}
  \label{def:MS-algebra}
  A $\mathbf{T}(0)$-MS-algebra is a commutative unital algebra object $\Bf{A} \in \Sf{MS}_{\CC}$, i.e., $\Bf{A} \in \Sf{MS}_{\CC}$ with morphisms in $\Sf{MS}_{\CC}$: 
  \[\mathbf{T}(0) \xrightarrow{u} \Bf{A}, \tab \Bf{A} \otimes \Bf{A} \xrightarrow{\mu} \Bf{A}\]
  satisfying the usual unit, associativity, and commutativity diagrams. Similarly, we can define $\Bf{A}$-MS-algebra (and module) from a $\mathbf{T}(0)$-MS-algebra $\Bf{A}$. \\
  \tab We say that $\Bf{A}$ is artinian local if $\omega(\Bf{A})$ is an artinian local $\CC-$algebra, and the residue map $\Bf{A} \twoheadrightarrow \mathbf{T}(0)$ is a morphism in $\Sf{MS}_{\CC}$. Since $\Sf{MS}_{\CC}$ is abelian, $\Ide{m}_{\Bf{A}} \coloneq \ker(\Bf{A} \twoheadrightarrow \mathbf{T}(0))$ is a subobject of $\Bf{A}$ in $\Sf{MS}_{\CC}$, and we have a splitting $\Bf{A} = \mathbf{T}(0) \oplus \Ide{m}_{\Bf{A}}$. 
\end{definition}

\begin{definition}
  \label{def:pro-MS}
  \leavevmode
  \begin{enumerate}
    \item A pro-MS is a pro-object in $\Sf{MS}_{\CC}$. 
    \item A pro-$\mathbf{T}(0)$-MS-algebra is a $\mathbf{T}(0)$-algebra object in the category of pro-MS. Equivalently, it is an inverse system of $\mathbf{T}(0)$-MS-algebras. 
    \item A pro-$\mathbf{T}(0)$-MS-algebra $\bm{\Lambda}$ is complete, noetherian, artinian, or local when $\omega(\bm \Lambda)$ is as a $\CC$-algebra. 
    \item An ideal subobject $\Bf{I}$ of a (pro-)$\mathbf{T}(0)$-MS-algebra is a (pro-)$\Bf{T}(0)$-MS-ideal if it is a sub pro-$\CC$-MS. 
    \item If $\bm \Lambda$ is a pro-$\mathbf{T}(0)$-MS-algebra, we define a pro-$\bm \Lambda$-MS-module to be a pro-$\bm \Lambda$-module object in $\Sf{MS}_{\CC}$. 
  \end{enumerate}
\end{definition}

Let $a_X: X \to \Rm{pt}$, and denote $\sh{A} = a_X^*\Bf{A}$ for $\Bf{A} \in \Sf{MS}_{\CC}$, i.e., $\sh{A}$ is the constant variation with fiber $\omega(\Bf{A})$. 
\begin{definition}
  \label{def:A-VMS}
  Let $\Bf{A} \in \Sf{MS}_{\CC}$ be an artinian local $\mathbf{T}(0)$-MS-algebra. An $\Bf{A}$-variation of mixed structures ($\Bf{A}$-VMS) is a complex variation of mixed Hodge/twistor structures $\sh{V}$ together with an action morphism 
  \[\sh{A} \otimes \sh{V} \to \sh{V}\]
  which is a morphism of variations of mixed structures and which satisfies unit and associativity diagrams. We also require the underlying local system of $\sh{V}$ to be a local system of free $\omega(\Bf A)$-modules. \\
  \tab If $\bm \Lambda$ is a pro-$\mathbf{T}(0)$-MS-algebra with 
  \[\bm \Lambda = \varprojlim_{n} \bm{\Lambda}_n, \tab \bm{\Lambda}_n = \bm{\Lambda}/\Ide{m}_{\bm \Lambda}^{n + 1}\]
  where each $\bm{\Lambda}_n$ is an artinian local $\mathbf{T}(0)$-MS-algebra with $\bm{\Lambda}_{n + 1} \to \bm{\Lambda}_n$ being a morphism in $\Sf{MS}_{\CC}$, then a pro-$\bm \Lambda$-VMS is of the form 
  \[\sh{V} = \varprojlim_{n} \sh{V}_n, \tab \sh{V}_n = \sh{V}/\Ide{m}_{\bm \Lambda}^{n + 1}\sh{V}\]
  where $\sh{V}$ has an underlying local system of free $\omega(\bm \Lambda)$-modules. Each $\sh{V}_n$ is a $\CC$-VMS for which the maps $\sh{V}_{n + 1} \to \sh{V}_n$ are morphisms of $\CC$-VMS, each $\bm \Lambda/\Ide{m}_{\bm \Lambda}^{n + 1} \subset \Rm{End}(\sh{V}_n)$ is a sub $\CC$-MS. 
\end{definition}

\begin{prop}
  \label{lem:A-direct-img}
  Let $X$ be a smooth projective complex variety. Let $\Bf{A} \in \Sf{MS}_{\CC}$ be an artinian local $\mathbf{T}(0)$-algebra, and $\sh{V}$ an $\Bf{A}$-VMS with underlying local system $\bb{V}$. Then there is a functorial
  \[K_{\Bf{A}}(\sh{V}) = \Rm{R}a_{X, *}^{\Rm{MS}} \sh{V} \in D^b(\Bf{A}\Rm{\text{-}MS\text{-}mod})\]
  such that $\omega(K_{\Bf{A}}(\sh{V})) \simeq \Rm{R}\Gamma(X, \bb{V})$ in $D^b(\omega(\Bf A)\text{-}\Rm{mod})$. Furthermore, for every morphism $\Bf{A} \to \Bf{B}$ in $\Sf{MS}_{\CC}$ with the corresponding $\sh{A} \to \sh{B}$ of constant variations on $X$, we have base change 
  \[K_{\Bf{A}}(\sh{V}) \otimes_{\Bf{A}}^{\Rm{L}} \Bf{B} \simeq K_{\Bf{B}}(\sh{V} \otimes_{\sh{A}} \sh{B})\]
\end{prop}
\begin{proof}
  Forget about the coefficient $\Bf{A}$, we can regard $\sh{V}$ as a smooth object in the derived category of algebraic complex mixed modules. Then the fact that we have $K(\sh{V}) = \Rm{R}a_{X, *}^{\Rm{MS}} \sh{V} \in D^b(\Sf{MS}_{\CC})$ with $\omega(K(\sh{V})) \simeq \Rm{R}\Gamma(X, \bb{V})$ follows from the six functor formalism. In the Hodge setting, this follows from work of Saito \cite[Theorem 0.1]{saito1990mixed}; in the twistor setting, this follows from works of Sabbah and Mochizuki \cite{sabbah2005polarizable, mochizuki2015mixed}. We now show compatibility with the $\Bf{A}$-structure. Apply the projection formula, which in this specific case is a consequence of six functor formalism and dualizability of $\Bf{A}$, to the action morphism $(a_X)^*\Bf{A} \otimes \sh{V} \to \sh{V}$ we get 
  \[\Bf{A} \otimes \Rm{R}a_{X, *}^{\Rm{MS}} \sh{V} \to \Rm{R}a_{X, *}^{\Rm{MS}} \sh{V}\]
  and similarly for the unit and associativity diagrams. Regarding the base change $\Bf{A} \to \Bf{B}$, notice that we have a base change morphism by adjunction
  \[\Bf{B} \otimes_{\Bf{A}}^{\Rm{L}} \Rm{R}a_{X, *}^{\Rm{MS}} \sh{V} \to \Rm{R}a_{X, *}^{\Rm{MS}}(\sh{V} \otimes_{\sh{A}} \sh{B})\] 
  which, after applying $\omega$, becomes the base change morphism $\Rm{R}\Gamma(X, \bb{V}) \otimes_{\omega(\Bf A)}^{\Rm{L}} \omega(\Bf B) \to \Rm{R}\Gamma(X, \bb{V} \otimes_{A_X} B_X)$. This follows from the compatibility of $\omega$ with the six functors. The base change $\Rm{R}\Gamma(X, \bb{V}) \otimes_{\omega(\Bf A)}^{\Rm{L}} \omega(\Bf B) \to \Rm{R}\Gamma(X, \bb{V} \otimes_{A_X} B_X)$ is an isomorphism due to identifications with free cellular complexes, see \eqref{eq:cellular-base-change}. Finally, $\omega$ is faithful and exact, hence its derived version is conservative, and $\Bf B \otimes_{\Bf A}^{\Rm{L}} \Rm{R}a_{X, *}^{\Rm{MS}} \sh{V} \to \Rm{R}a_{X, *}^{\Rm{MS}}(\sh{V} \otimes_{\sh{A}} \sh{B})$ must be an isomorphism to begin with.
\end{proof}

\subsection{Absolute subset and rank one theory}
\label{sec:rank-one}
We have a quasi-projective moduli space $R_{\DR}(X, x, r)$ of framed flat connections on $X$ \cite[Theorem 6.13]{simpson1994moduli2}. The Riemann-Hilbert correspondence then induces a complex analytic isomorphism \cite[Theorem 7.1]{simpson1994moduli2}:
\begin{equation}
  \label{eq:RH-moduli}
  \Sf{RH}: R_{\DR}(X, x, r) \xrightarrow{\simeq} R_B(X, x, r)
\end{equation}
\tab Given $(E, \nabla) \in R_{\DR}(X, x, r)$ and $\sigma \in \Rm{Aut}(\CC/\QQ)$, we have a base change $X^{\sigma} \to X$ and pulling back gives $(E^{\sigma}, \nabla^{\sigma}) \in R_{\DR}(X^{\sigma}, x, r)$. In fact, this gives a bijection $\sigma: R_{\DR}(X, x, r) \xrightarrow{\simeq} R_{\DR}(X^{\sigma}, r)$. 
\begin{definition}
  \label{def:absolute-set}
  A $\Rm{GL}_r$-invariant subset $\Sigma \subseteq R_B(X, x, r)$ is absolute $\ol{\QQ}$-constructible if for any $\sigma \in \Rm{Aut}(\CC/\QQ)$, 
  \[\Sigma^{\sigma} \coloneq \Sf{RH} \circ \sigma \circ \Sf{RH}^{-1}(\Sigma) \subseteq R_B(X^\sigma, x, r)\]
  is $\ol{\QQ}$-constructible, i.e., a finite union of locally closed subvarieties defined over $\ol{\QQ}$.  
\end{definition}
We have similar notions for the stacks $\sh{M}_{B/\DR}(X, r)$ as well as the good moduli spaces $M_{B/\DR}(X, r)$. Now, consider the rank one character variety 
\begin{equation}
  \label{eq:rank1-decomp}
  M_B(X, 1) = \Hom(\pi_1(X, x), \CC^*) \simeq (\CC^*)^{\rank H_1(X, \ZZ)} \times \Hom(H_1(X, \ZZ)_{\Rm{tors}}, \CC^*)
\end{equation}
which is a finite disjoint union of algebraic tori.
\begin{theorem}[{\cite{budur2020absolutesets,simpson1993subspaces}}]
  \label{thm:rank1-absolute}
  Every irreducible absolute $\ol{\QQ}$-closed subset $\Sigma \subseteq M_B(X, 1)$ is a torsion-translated algebraic subtorus
  \[\Sigma = \tau \cdot N\]
  where $N \subseteq M_B(X, 1)$ is a subtorus, and $\Rm{Lie}(N) \subseteq H^1(X, \CC)$ is a complex Hodge substructure. Absolute $\ol{Q}$-constructible subsets are obtained from a finite sequence of unions, intersections and complements.
\end{theorem}
\begin{proof}
  The first and last statement follows from \cite[Theorem 6.1]{simpson1993subspaces} and \cite[Theorem 9.1.2]{budur2020absolutesets}. Regarding the Lie algebra, by \cite[Theorem 6.1]{simpson1993subspaces}, $N$ is a triple torus, hence by \cite[Theorem 2.1]{simpson1993subspaces} there exists abelian subvariety $P \subseteq \Rm{Alb}(X)$ such that $N$ is the image of $M_B(\Rm{Alb}(X)/P, 1) \hookrightarrow M_B(X, 1)$. The map on tangent space is
  \[H^1(\Rm{Alb}(X)/P, \CC) \hookrightarrow H^1(X, \CC)\]
  is induced by $\pi^*: H^1(\Rm{Alb}(X)/P, \CC) \to H^1(\Rm{Alb}(X), \CC) \simeq H^1(X, \CC)$ where $\pi: \Rm{Alb}(X) \to \Rm{Alb}(X)/P$ is the projection. Pullback is a Hodge morphism hence we are done.  
\end{proof}
\begin{coroll}[{\cite[Theorem 11.1.1]{budur2020absolutesets}}]
  \label{lem:rank1-jump}
  Let $\sh{F} \in D^b_c(X, \CC)$ be of geometric origin, then the twisted jump loci 
  \[\Sigma^i_k(X, \sh{F}) = \set*{\bb{V} \in R_B(X, 1)\ \big\vert\ \dim \HH^i(X, \sh{F} \otimes \bb{V}) \geq k}\]
  is absolute $\ol{\QQ}$-closed and thus a finite union of torsion-translated algebraic subtori. If $\sh{F}$ is not of geometric origin, then $\Sigma^i_k(X, \sh{F})$ is only a finite union of translated algebraic subtori. 
\end{coroll}
\subsection{Determinant sheet}
\label{sec:jump-product}
In this section, let $A = M_B(X, 1)$ and the decomposition \eqref{eq:rank1-decomp} be $A = A^\circ \times \hat{T}$, where $T = H_1(X, \ZZ)_{\Rm{tors}}$. Since $T$ is finite, its dual $\hat{T}$ is also finite. We have the squaring map 
\[[2]: A \to A, \tab \bb{L} \mapsto \bb{L}^{\otimes 2}\]
and let $A^2 = \im [2], A[2] = \ker [2]$. Then $A/A^2 \simeq \hat{T}/\hat{T}^2$ is a finite group. For each $\bar{\epsilon} \in A/A^2$ we can pick a torsion lift $\epsilon \in A$. Let $\epsilon A^2$ be the coset $\epsilon A^2 = \set*{\epsilon \otimes \bb{L}^{\otimes 2}\ \vert\ \bb{L} \in A}$, then 
\[A = \bigsqcup_{\bar{\epsilon} \in A/A^2} \epsilon A^2\]
and thus 
\[M_B(X, 2) = \bigsqcup_{\bar{\epsilon} \in A/A^2} M_B(X, 2, \det \in \epsilon A^2)\]
\tab For each $\bar{\epsilon} \in A/A^2$ with torsion lift $\epsilon \in A$, we have the following proposition 
\begin{prop}
  \label{lem:determinant-sheet}
  There is a Cartesian diagram 
  \[\begin{tikzcd}
	{M_B(X, 2, \det = \epsilon) \times A} && {M_B(X, 2, \det \in \epsilon A^2)} \\
	\\
	A && A
	\arrow["{(\bb{V}, \bb{L}) \mapsto \bb{V} \otimes \bb{L}}","{\Phi_{\epsilon}}"', from=1-1, to=1-3]
	\arrow["{\Rm{pr}_2}"', from=1-1, to=3-1]
	\arrow["\det", from=1-3, to=3-3]
	\arrow["{\bb{L} \mapsto \epsilon \otimes \bb{L}^{\otimes 2}}", from=3-1, to=3-3]
  \end{tikzcd}\tag{2.5}\label{eq:det-prod-diagram}\]
  and $\Phi_{\epsilon}$ is a finite \'{e}tale $A[2]$-torsor over $M_B(X, 2, \det \in \epsilon A^2)$, with the action 
  \[\rho \cdot (\bb{V}, \bb{L}) = (\bb{V} \otimes \rho, \bb{L} \otimes \rho^{-1})\]
  The same is true for the quotient stacks and framed moduli schemes. 
\end{prop}
\begin{proof}
  Notice that 
  \[M_B(X, 2, \det \in \epsilon A^2) \times_A A = \set*{(\bb{W}, \bb{L})\ \big\vert\ \det \bb{W} = \epsilon \otimes \bb{L}^{\otimes 2}}\]
  which is isomorphic to $M_B(X, 2, \det = \epsilon) \times A$ via the maps 
  \begin{align*}
    M_B(X, 2, \det \in \epsilon A^2) \times_A A &\to M_B(X, 2, \det = \epsilon) \times A\\ 
    (\bb{W}, \bb{L}) &\mapsto (\bb{W} \otimes \bb{L}^{-1}, \bb{L})\\
    (\bb{V} \otimes \bb{L}, \bb{L}) &\mapsfrom (\bb{V}, \bb{L})
  \end{align*}
  \tab Next, the map downstair map is a $A[2]$-torsor over $\epsilon A^2$, with the action $\rho \cdot \bb{L} = \bb{L} \otimes \rho^{-1}$. Hence $\Phi_{\epsilon}$ is also a $A[2]$-torsor. The induced action on $M_B(X, 2, \det \in \epsilon A^2) \times_A A$ is $\rho \cdot (\bb{W}, \bb{L}) = (\bb{W}, \bb{L} \otimes \rho^{-1})$ which translate to the stated action on $M_B(X, 2, \det = \epsilon) \times A$. 
\end{proof}
Now, let $Z \subset M_B(X, 2)$ be an irreducible subvariety, then it must be contained in some $M_B(X, 2, \det \in \epsilon A^2)$. We can then look at $\phi_{\epsilon}^{-1}(Z)$ and pick an irreducible component $C$. Set $H_C = \set*{\rho \in A[2]\ \vert\ \rho \cdot C = C}$. 
\begin{prop}
  \label{lem:normalized-det-sheet}
  The induced map between normalizations $C^{\nu} \to Z^{\nu}$ is a finite \'{e}tale $H_C$-torsor, and $Z^{\nu} \simeq C^{\nu}/H_C$. 
\end{prop}
\begin{proof}
  We have a base change diagram 
  \[\begin{tikzcd}
	W && {\Phi_{\epsilon}^{-1}(Z)} \\
	\\
	{Z^{\nu}} && Z
	\arrow["f", from=1-1, to=1-3]
	\arrow["g"', from=1-1, to=3-1]
	\arrow["{\phi_{\chi}}", from=1-3, to=3-3]
	\arrow["\nu", from=3-1, to=3-3]
\end{tikzcd}\]
  where $g$ is a finite \'{e}tale $A[2]-$torsor. Let $\eta_Z$ be the generic point of $Z$, then the generic point $\eta_{C}$ of $C$ lies over $\eta_Z$. Notice that we have an isomorphism of function fields $\kappa(\eta_{Z^{\nu}}) \simeq \kappa(\eta_{Z})$, hence the an isomorphism of generic fibers $W_{\eta_{Z^{\nu}}} \simeq \Phi^{-1}_{\epsilon}(Z)_{\eta_Z}$. Thus there is a $W_{\eta_{Z^{\nu}}} \ni \eta \mapsto \eta_C$ under this identification. Let $\ol{C}$ be the closure of $\eta$ in $W$. Since $W \to Z^{\nu}$ is finite \'{e}tale, $\ol{C}$ is an irreducible component of $W$. It follows that $\ol{C}$ is open and closed in $W$, and $\ol{C} \to Z^{\nu}$ is a finite \'{e}tale $H_C-$torsor. \\
  \tab On the other hand, $W \to \Phi_{\epsilon}^{-1}(Z)$ is finite, since it is the base change of $\nu$ which is finite. Thus $\ol{C} \to C$ is finite. It is also birational by construction, hence $\ol{C} = C^\nu$. 
\end{proof}

\section{Hodge theory of $R_B$ and $M_B$} 
\label{sec:hodge-jump}
In this section, we prove that the cohomological jump loci enjoy several nice Hodge theoretic properties. We first briefly recall some results and definitions, following the exposition of \cite{bakker2024linear}.
\subsection{Betti universal deformation algebra and Deligne-Hitchin space}
\label{sec:Deligne-Hitchin}
\begin{theorem}[{\cite[Theorem 4.2]{bakker2024linear}}]
  \label{thm:Betti-MTS}
  Let $(\bb{V}, \phi)$ be a framed rank $r$ local system underlying a $\CC$-VMS $\sh{V}$, then the completed local ring $\Lambda_{(\bb{V}, \phi)} \coloneq \widehat{\sh{O}}_{R_B(X, x, r), (\bb{V}, \phi)}$ has a pro-$\mathbf{T}(0)$-MS-algebra structure $\bm{\Lambda}_{(\bb{V}, \phi)}$ with 
  \[\omega(\bm{\Lambda}_{(\bb{V}, \phi)}) = \Lambda_{(\bb{V}, \phi)}\] 
  and the universal formal local system $\widehat{\bb{V}}$ underlies a free pro-$\bm{\Lambda}_{(\bb{V}, \phi)}$-VMS $\widehat{\sh{V}}$. Let $\Bf{M}$ be $\CC^r$ equipped with the mixed structure induced by $\phi$. 
  The universal framing $\widehat{\phi}: \widehat{\bb{V}}_x \to \Lambda_{(\bb{V}, \phi)} \widehat{\otimes}_\CC \CC^r$ is the image under $\omega$ of a morphism of pro-$\bm{\Lambda}_{(\bb{V}, \phi)}$-MS-modules $\widehat{\phi}: \widehat{\sh{V}}_x \to \bm{\Lambda}_{(\bb{V}, \phi)} \widehat{\otimes}_{\mathbf{T}(0)} \Bf{M}$.\\
  \tab These structures are uniquely determined by the following universal property: For any artinian local $\mathbf{T}(0)-$MS-algebra $\Bf{A}$ and any $\Bf{A}$-VMS $\sh{U}$ with a framing $\psi: \sh{U}_x \to \Bf{A} \otimes_{\Bf{T}(0)} \Bf{M}$ which restricts to $(\sh{V}, \phi)\ \Rm{mod}\ \Ide{m}_{\Bf{A}}$, there is an unique morphism $\bm{\Lambda}_{(\bb{V}, \phi)} \to \Bf{A}$ of local pro-$\Bf{T}(0)$-MS-algebras such that 
  \[\brac*{\Bf{A}, (\sh{U}, \psi)} \simeq \Bf{A} \widehat{\otimes}_{\bm{\Lambda}_{(\bb{V}, \phi)}} \brac*{\bm{\Lambda}_{(\bb{V}, \phi)}, (\widehat{\sh{V}}, \widehat{\phi})}\]
  or, in other words, $\bm{\Lambda}_{(\bb{V}, \phi)}$ represents the functor 
  \[\Bf{A} \mapsto \{\text{framed free $\Bf{A}$-VMS deformations of $(\sh{V}, \phi)$}\}/\simeq\] 
\end{theorem}
For a $\CC-$scheme $T$, denote $X_{T} = X \times T$ with projection $\pi_T: X_T \to T$.
\begin{prop}[{\cite[Proposition 4.1]{simpson1996hodgefilt}}]
  The moduli functor 
  \[\mathsf{R}_{\Rm{Hod}}(r)(T) = \left\{
    \begin{array}{c|l}
      (\lambda, E_{T}, \nabla_{\lambda,T}, \phi) & E_T \text{ locally free of rank $r$ on }X_T,\\ 
      & \lambda \in \Gamma(T, \sh{O}_T),\\
      & \pi_T^{-1}\sh{O}_T-\text{linear }\nabla_{\lambda, T}: E_T \to E_T \otimes \Omega^1_{X_T/T},\\
      & \nabla_{\lambda, T}^2 = 0, \\
      & \nabla_{\lambda, T}(fs) = \lambda \der_{X_T/T} f \otimes s + f \nabla_{\lambda, T} s\\
      &\phi: E_T\vert_{x \times T} \simeq \sh{O}_T^r\\
      &\text{over }\lambda = 0\ \text{the resulting Higgs bundles are}\\
      &\tab \text{semistable with vanishing $\QQ$-Chern classes}
    \end{array}
  \right\}\Big/\simeq \]
  is representable by a relative scheme $R_{\Hod}(X, x, r) \xrightarrow{\lambda} \AAA^1_{\lambda}$ parametrizing $\lambda$-connections with a framing. Furthermore, there is a $\GG_m$-action on $R_{\Hod}(X, x, r)$, covering the scaling action on $\AAA^1_\lambda$, by scaling the $\lambda-$connection. 
\end{prop}
\begin{definition}
  \label{def:DH-space}
  Consider the complex conjugate variety $X^c$, we can similarly define $R_{\Hod}(X^c, x, r)$ and there is a holomorphic isomorphism 
\begin{equation}
  \label{eq:Hod-conj}
  R_{\Hod}(X, x, r)\vert_{\GG_m} \simeq R_{\Hod}(X^c, x, r)\vert_{\GG_m}
\end{equation}
and there is a gluing to get the analytic \textit{Deligne-Hitchin} space $R_{\Rm{DH}}(X, x, r) \to \PP^1$. We have the derivations associated to the $\GG_m$-actions on $R_{\Hod}(X, x, r)$ and $R_{\Hod}(X^c, x, r)$ which glue to a derivation $\Theta$ on $R_{\Rm{DH}}(X, x, r)$. 
\end{definition}
\begin{prop}[{\cite[Corollary 8.8]{bakker2024linear}}]
  \label{lem:classifying-map}
  Let $\Bf{A}$ be an artinian local $\Bf{T}(0)$-MS-algebra. A framed $\Bf{A}$-VMTS $\scr{V}$ induces a $\PP^1$-morphism $\Spec_{\PP^1} \Bf{A} \to R_{\Rm{DH}}(X, x, r)$. 
\end{prop}
\begin{definition}
  If $\scr{V}$ is a $\CC$-VTS on $X$, then the resulting $\PP^1 \to R_{\Rm{DH}}(X, x, r)$ is called a preferred section. 
\end{definition}
More concretely, if $(\bb{V}, \phi) \in R_B(X, x, r)$ is a semisimple framed local system, then it underlies a canonical weight 0 $\CC$-VTS, giving a preferred section 
\begin{equation}
  \label{eq:preferred-sec}
  s_{(\bb{V}, \phi)}: \PP^1 \to R_{\Rm{DH}}(X, x, r)
\end{equation}
where $s_{(\bb{V}, \phi)}(1)$ is the flat bundle corresponding to $(\bb{V}, \phi)$ via the Riemann-Hilbert correspondence. 
\begin{prop}[{\cite[Lemma 8.11]{bakker2024linear}, \cite[Theorem 8.1]{simpson1997mixedtwistor}}]
  \label{lem:formal-section}
  We have a classifying map for the universal pro-$\bm{\Lambda}_{(\bb{V}, \phi)}$-VMS $\widehat{\scr{V}}$,
  \begin{equation}
    \label{eq:formal-sec}
    \widehat{s}_{(\bb{V}, \phi)}: \Rm{Spf}_{\PP^1} \bm{\Lambda}_{(\bb{V}, \phi)} \xrightarrow{\simeq} \widehat{R}_{\Rm{DH}}(X, x, r)_{s_{(\bb{V}, \phi)}}
  \end{equation}
  which is an isomorphism onto the completion of $R_{\Rm{DH}}(X, x, r)$ along the preferred section. 
\end{prop}

\begin{prop}[{\cite[Corollary 4.2]{simpson1992higgs}}]
  \label{lem:C*-action}
  There is a continuous $\CC^*$ action on $M_B(X, r)$ with fixed points underlying canonical weight zero complex variations of Hodge structures.
\end{prop}

\subsection{Jump loci are Hodge subvarieties}
\label{sec:hodge=jump}
Let $\sh{M}_B(X, r)$ be the quotient stack of local systems of rank $r$, and $\sh{Z} \subseteq \sh{M}_B(X, r)$ locally closed substack. Let $\bb{V} \in \sh{Z}(\CC)$ be a semisimple local system, and pick a framing $\phi: \bb{V}_x \xrightarrow{\simeq} \CC^r$. Denote $R\sh{Z} \subseteq R_B(X, x, r)$ the base change to the framed space. 
\begin{definition}[{\cite[Section 8.5]{bakker2024linear}}]
  \label{def:Hodge-substack}
  $\sh{Z}$ is \textit{Hodge} at $\bb{V}$ if 
  \begin{description}
    \item[(H1)] equipping $\bb{V}$ with the canonical weight 0 $\CC-$VTS, the quotient $\Lambda_{(\bb{V}, \phi)} \twoheadrightarrow \widehat{\sh{O}}_{R\sh{Z}, (\bb{V}, \phi)}$ underlies a quotient map of pro-$\Bf{T}(0)$-MTS-algebras. If $\bb{V}$ is a $\CC^*$-fixed point, then we can equip $\bb{V}$ with a canonical weight 0, and the quotient $\Lambda_{(\bb{V}, \phi)} \twoheadrightarrow \widehat{\sh{O}}_{R\sh{Z}, (\bb{V}, \phi)}$ underlies a quotient map of pro-$\Bf{T}(0)$-MHS-algebra. Denote $\bm{\Gamma}_{\sh{Z}, (\bb{V}, \phi)}$ the corresponding pro-$\Bf{T}(0)$-MS-algebra on $\widehat{\sh{O}}_{R\sh{Z}, (\bb{V}, \phi)}$. 
    \item[(H2)] The formal twistor subspace $\Rm{Spf}_{\PP^1} \bm{\Gamma}_{\sh{Z}, (\bb{V}, \phi)} \subseteq \Rm{Spf}_{\PP^1} \bm{\Lambda}_{(\bb{V}, \phi)}$ is tangent to the derivation $\Theta$ in \cref{def:DH-space}. 
    \item[(H3)] We have an induced map by composition 
    \[\widehat{s}_{R\sh{Z}, (\bb{V}, \phi)}: \Rm{Spf}_{\PP^1} \bm{\Gamma}_{\sh{Z}, (\bb{V}, \phi)} \to \widehat{R}_{\Rm{DH}}(X, x, r)_{s_{(\bb{V}, \phi)}}\]
    Denote $Z_{\Hod}$ (resp. $Z_{\ol{\Hod}}$) the Zariski closure of the image of $\widehat{s}_{R\sh{Z}, (\bb{V}, \phi)}\vert_{\AAA^1_\lambda}$ (resp. $\widehat{s}_{R\sh{Z}, (\bb{V}, \phi)}\vert_{\PP^1 \setminus \set*{0}}$) in $R_{\Hod}(X, x, r)$ (resp. $R_{\Hod}(X^c, x, r)$). Then $\widehat{s}_{R\sh{Z}, (\bb{V}, \phi)}\vert_{\GG_m}$ induces an isomorphism onto the completion of $Z_{\Hod}\vert_{\GG_m}$ along the image of $s_{(\bb{V}, \phi)}$, and the same is true for $\widehat{s}_{R\sh{Z}, (\bb{V}, \phi)}\vert_{\GG_m}$ and $Z_{\ol \Hod}\vert_{\GG_m}$. 
  \end{description}
  $\sh{Z}$ is a \textit{Hodge substack} if it is closed under semisimplification and Hodge at every semisimple $\bb{V} \in \sh{Z}(\CC)$. Notice that these conditions rely only on $R\sh{Z} \subseteq R_B(X, x, r)$, so we have a similar notion of \textit{Hodge subvariety} for a $\Rm{GL}_r$-invariant subvariety $Z \subseteq R_B(X, x, r)$. We also have obvious generalization to Hodge subvarieties of products of character varieties. 
\end{definition}

\begin{lemma}
  \label{lem:enriched-ideal}
  Let $\Bf{A} \in \Sf{MS}_{\CC}$ be an artinian local $\Bf{T}(0)$-MS-algebra, with $A = \omega(\Bf{A})$, and $\Bf{K} \in D^b(\Bf{A}\text{-}\Rm{MS}\text{-}\Rm{mod})$ with $K = \omega(\Bf{K})$ quasi-isomorphic to a bounded above complex of free $A$-modules. Then $J^i_k(K) \subseteq A$ is the image under $\omega$ of an ideal subobject $J^i_k(\Bf{K}) \subseteq \Bf{A}$. This construction commutes with $\Bf{T}(0)$-MS-algebra coefficient base change. 
\end{lemma}
\begin{proof}
  Let $G = \Rm{Aut}^{\otimes}(\omega)$ be the Tannakian dual of $\Sf{MS}_{\CC}$. Then we have an equivalence of category $\Sf{MS}_{\CC} \simeq \Rm{Rep}_{\CC}^{\Rm{fd}}(G)$. Under this equivalence, $\Bf{A}$ is identified with $A = \omega(\Bf{A})$ with an algebraic $G-$action, i.e., a comodule structure $\delta_A: A \to A \otimes \sh{O}(G)$; this is an action by algebra automorphism. Similarly, $\Bf{K}$ gives a $G$-equivariant complex $K = \omega(\Bf{K})$ of $A-$modules. \\
  \tab Let $B$ be a $\CC$-algebra, and $g \in G(B)$, all the tensor product below will be over $\CC$ unless specified otherwise. We have a $B$-algebra automorphism $g: A \otimes B \to A \otimes B$. Since $B$ is flat over $\CC$, the action of $g$ on $K \otimes^{\Rm{L}} B = K \otimes B$ satisfies
  \[g \cdot (bm) = g(b)(g \cdot m) \tab b \in A \otimes B, m \in K^i \otimes B\]
  coming from applying $\omega$ to the action morphism $\Bf{A} \otimes \Bf{K} \to \Bf{K}$. It follows that 
  \[g^*(K \otimes B) \coloneq (A \otimes B) \otimes_{g} (K \otimes B) \simeq K \otimes B\]
  since $g: A \otimes B \to A \otimes B$ is an isomorphism. Determinantal ideals commute with arbitrary base change, and by \eqref{eq:base-change-jump-ideal} we get 
  \[g(J^i_k(K\otimes B)) = J^i_k(g^*(K \otimes B)) \simeq J^i_k(K \otimes B)\]
  thus $J^i_k(K \otimes B)$ is stable under the action of $g \in G(B)$. It follows that $J^i_k(K)$ is a $G$-subrepresentation. The equivalence then gives a unique subobject $J^i_k(\Bf{K}) \subseteq \Bf{A}$ in $\Sf{MS}_{\CC}$ with $\omega(J^i_k(\Bf{K})) = J^i_k(K)$. To check that $J^i_k(\Bf{K})$ is an ideal subobject, it suffices to check that the composition
  \[\Bf{A} \otimes J^i_k(\Bf{K}) \to \Bf{A} \otimes \Bf{A} \xrightarrow{\mu} \Bf{A}\]
  factors through $J^i_k(\Bf{K}) \subseteq \Bf{A}$. Equivalently, 
  \[\Bf{A} \otimes J^i_k(\Bf{K}) \to \Bf{A} \to \Bf{A}/J^i_k(\Bf{K})\]
  is zero. This follows from $\omega$ being faithful and the image of that map under $\omega$ is zero. Compatibility with coefficient base change follows from \eqref{eq:base-change-jump-ideal}. 
\end{proof}
The jump loci $\Sigma^i_k(X, r) \subseteq R_B(X, x, r)$ is invariant under the conjugation action by $\Rm{GL}_r$, so we can ask whether it is a Hodge subvariety. 
\begin{theorem}
  \label{thm:jump=Hodge}
  For smooth projective $X$, the jump locus $\Sigma_k^i(X, r) \subseteq R_B(X, x, r)$ is a Hodge subvariety. Furthermore, it is closed under semisimplification, and $\Sigma_k^i(X, r)_{\red}$ is a Hodge subvariety which is absolute $\QQ$-closed. \\
  \tab Every irreducible component of $\Sigma_k^i(X, r)_{\red}$ is a Hodge subvariety and is absolute $\ol{\QQ}-$closed. 
\end{theorem}
\begin{proof}
  We first verify that $\Sigma^i_k(X, r)$ satisfies the three conditions in \cref{def:Hodge-substack}. 

  \medskip \noindent
  $\Bf{(H1)}$. Consider $(\bb{V}, \phi) \in R_B(X, x, r)$ semisimple, and let 
  \[\Lambda = \widehat{\sh{O}}_{R_B(X, x, r), (\bb{V}, \phi)}, \tab \Lambda_n = \Lambda/\Ide{m}_{\Lambda}^{n + 1}\]
  \tab Equip $\bb{V}$ with the canonical weight 0 $\CC$-VTS. If $\bb{V}$ is a $\CC^*$-fixed point, then it underlies a canonical weight 0 $\CC$-VHS. By \cref{thm:Betti-MTS}, $\Lambda$ underlies a pro-$\Bf{T}(0)$-MS-algebra $\bm \Lambda$, and the universal formal local system $(\widehat{\bb{V}}, \widehat{\phi})$ underlies a pro-$\bm \Lambda$-VMS $(\widehat{\sh{V}}, \widehat{\phi})$. Each $\Lambda_n$ underlies a $\Bf{T}(0)$-MS-algebra $\Lambda_n = \bm{\Lambda}/\Ide{m}_{\bm \Lambda}^{n + 1}$. Write $\sh{V}_n = \widehat{\sh{V}} \widehat{\otimes}_{\bm \Lambda} \bm{\Lambda}_n$, and $\bb{V}_n$ the underlying local system of $\sh{V}_n$. By \cref{lem:A-direct-img} we have 
  \[\Bf{K}_n = \Rm{R}a_{X, *}^{\Rm{MS}} \sh{V}_n \in D^b(\bm{\Lambda}_n\Rm{\text{-}MS\text{-}mod}), \tab \omega(\Bf{K}_n) \simeq \Rm{R}\Gamma(X, \bb{V}_n)\]
  \tab Let $R = \Gamma(R_B(X, x, r), \sh{O}_{R_B(X, x, r)})$, then $R_B(X, x, r) = \Spec R$. By the discussion of universal cellular complex in \cref{sec:jump-ideal}, we get a complex $C^\bullet(X, \bb{V}^{\univ}) \simeq \Rm{R}\Gamma(X, \bb{V}^{\univ})$ of free $R$-modules which we can think of as a complex of free sheaves on $R_B(X, x, r)$. Its fiber over a point $(\bb{W}, \psi) \in R_B(X, x, r)$ gives $\Rm{R}\Gamma(X, \bb{W})$. Notice that in $\Lambda$, the jump ideal becomes 
  \[J^i_k(\Rm{R}\Gamma(X, \bb{V}^{\univ}))\cdot \Lambda = \varprojlim J^i_k(\Rm{R}\Gamma(X, \bb{V}^{\univ})) \cdot \Lambda_n\]  
  \tab Recall that we have base change 
  \begin{align*}
    J^i_k(\Rm{R}\Gamma(X, \bb{V}^{\univ})) \cdot \Lambda_n &\simeq J^i_k(\Rm{R}\Gamma(X, \bb{V}^{\univ}) \otimes^{\Rm{L}}_{R} \Lambda_n)\\
    &\simeq J^i_k(\Rm{R}\Gamma(X, \bb{V}^{\univ}\otimes^{\Rm{L}}_{R} \Lambda_n))\\
    &\simeq J^i_k(\Rm{R}\Gamma(X, \bb{V}_n))
  \end{align*}
  thus $J^i_k(\Rm{R}\Gamma(X, \bb{V}^{\univ}))\cdot \Lambda \simeq J^i_k(\Rm{R}\Gamma(X, \widehat{\bb{V}}))$. Now, $\omega(\Bf{K}_n)$ is quasi-isomorphic to a complex of free $\Lambda_n$-modules, which is $C^\bullet(X, \bb{V}_n)$, hence by \cref{lem:enriched-ideal} we get that $J^i_k(\Rm{R}\Gamma(X, \widehat{\bb{V}}))$ underlies a pro-$\Bf{T}(0)$-MS-ideal $J^i_k(\Bf{K}) \subseteq \bm \Lambda$. It follows that we get $\Bf{(H1)}$ since $\Sf{MS}_{\CC}$ is abelian. \\
  \\
  $\Bf{(H3)}$. We first take a detour. Let $\Bf{A} \in \Sf{MTS}_{\CC}$ with $A = \omega(\Bf{A})$, and $\sh{W}$ an $\Bf{A}$-VMTS with underlying local system $\bb{W}$ of free $A$-modules. By \cref{rmk:VMTS=lambda}, this is the same as a family of holomorphic integrable $\lambda$-connections. Let $X_{\lambda} = X \times \AAA^1_{\lambda}$ with projection $\pi: X_{\lambda} \to \AAA^1_{\lambda}$, then the family of $\lambda$-connections, over the patch $\AAA^1_\lambda$, is the same as a holomorphic bundle $\sh{W}_{\lambda}$ on $X_{\lambda}$ with an integrable $\nabla_{\lambda}: \sh{W}_{\lambda} \to \sh{W}_{\lambda} \otimes \Omega^1_{X_{\lambda}/\AAA^1_{\lambda}}$. The $\Bf{A}$-action becomes an $A[\lambda]-$action. We can then form the $\lambda$-De Rham complex 
  \[\DR_{\lambda}(\sh{W}_{\lambda}) = \sqbrac*{\sh{W}_{\lambda} \xrightarrow{\nabla_{\lambda}}\sh{W}_{\lambda} \otimes \Omega^1_{X_{\lambda}/\AAA^1_{\lambda}} \to \dots \to \sh{W}_{\lambda} \otimes \Omega^{\dim X}_{X_{\lambda}/\AAA^1_{\lambda}}}\]
  and $(\Rm{R}a_{X, *}^{\Rm{MS}}\sh{W}) \vert_{\AAA^1_{\lambda}} \simeq \Rm{R}\pi_* \DR_{\lambda}(\sh{W}_{\lambda}) \in D^b(A[\lambda]\Rm{\text{-}mod})$, by \cite[proof of Theorem 4.1]{simpson1997mixedtwistor} and the proof of \cref{lem:A-direct-img}. Denote 
  \[K_{\Hod}(\sh{W}) \coloneq \Rm{R}\pi_* \DR_{\lambda}(\sh{W}_{\lambda})\]
  \tab Next, let $q: R_{\Hod}(X, x, r) \to \AAA^1_{\lambda}$ be the Hodge moduli space. This carries a universal $\lambda$-connection $(\sh{V}^{\univ}_{\lambda}, \nabla_{\lambda})$ on $X \times R_{\Hod}$. If the classifying Deligne-Hitchin section for $\sh{W}$ is $s_{\sh{W}}: \Spec_{\PP^1}\Bf{A} \to R_{\Rm{DH}}(X, x, r)$, denote $s_{\sh{W}, \lambda} \coloneq s_{\sh{W}}\vert_{\AAA^1_{\lambda}}: \Spec A[\lambda] \to R_{\Hod}(X, x, r)$, and
  \[(\sh{W}_{\lambda}, \nabla_{\lambda}) = s^*_{\sh{W}, \lambda}(\sh{V}^{\univ}_{\lambda}, \nabla_{\lambda})\]
  \tab The complex $K_{\Hod}$ satisfies base change (same proof as in \cite[\href{https://stacks.math.columbia.edu/tag/0FM0}{Tag 0FM0}]{stacks-project} works), hence 
  \[K_{\Hod}(\sh{W}) \simeq s_{\sh{W}, \lambda}^* K_{\Hod}(\sh{V}^{\univ})\]
  where $K_{\Hod}(\sh{V}^{\univ}) \coloneq \Rm{R}\pi_*\DR_{\lambda}(\sh{V}^{\univ}_{\lambda}, \nabla_{\lambda})$. We can think of this as a complex of $A[\lambda]$-modules, hence a notion of $\lambda$-jump ideal makes sense. Denote $\Sigma^i_{k, \Hod}(X, r) = V(J^i_k(K_{\Hod}(\sh{V}^{\univ})))$. The restriction to $\lambda = 1$ gives the De Rham complex for the flat bundle corresponding to $\bb{W}$ via the Riemann-Hilbert correspondence, hence by the comparison theorem this gives $\Sigma^i_k(\bb{W})$. \\
  \\
  \tab Back to our situation, let $s_n: \Rm{Spec}_{\PP^1}\bm{\Lambda}_n \to R_{\Rm{DH}}(X, x, r)$ be the classifying map for the $\bm \Lambda_n$-VMTS $\sh{V}_n$. The reduction mod $\Ide{m}_{\bm \Lambda}$ gives $s_{(\bb{V}, \phi)}: \PP^1 \to R_{\Rm{DH}}(X, x, r)$, the preferred section associated to the canonical weight 0 $\CC$-VTS on semisimple $\bb{V}$. 
  Let $\bm{\Gamma} = \bm{\Lambda}/J^i_k(\Bf{K})$ be the pro-$\Bf{T}(0)$-MTS-algebra on $\widehat{\sh{O}}_{\Sigma^i_k, (\bb{V}, \phi)}$, and put 
  \[s_{n, \Sigma}: \Spec_{\PP^1} \bm{\Gamma}_n \to \Rm{Spec}_{\PP^1}\bm{\Lambda}_n \xrightarrow{s_n} R_{\Rm{DH}}(X, x, r)\]
  
  Notice that 
  \[s_n^* J^i_k(K_{\Hod}(\sh{V}^{\univ})) = J^i_k(s_n^*K_{\Hod}(\sh{V}^{\univ})) \simeq J^i_k(\Bf{K}_n\vert_{\AAA^1_{\lambda}}) \simeq J^i_k(\Bf{K}_n)\vert_{\AAA^1_{\lambda}}\]
  thus the image of $s_{n, \Sigma}\vert_{\GG_m}$ agrees with the $(n+1)$-th order neighborhood of the image of $s_{(\bb{V}, \phi)}\vert_{\GG_m}$ in $\Sigma^i_{k, \Hod}(X, r)$. It follows that the image of $\widehat{s}_{\Sigma}\vert_{\GG_m}$ agrees with the completion of $\Sigma^i_{k, \Hod}(X, r)$ along the image of $s_{(\bb{V}, \phi)}$, as desired. \\
  \\
  $\Bf{(H2)}$. The scaling $\GG_m$-action turns $t \cdot (\sh{V}^{\univ}_{\lambda}, \nabla_{\lambda}) \mapsto (\sh{V}^{\univ}_{\lambda}, t\nabla_{\lambda})$. Multiplication by $t^p$ in degree $p$ gives a quasi-isomorphism 
  \[\DR_{\lambda}(\sh{V}^{\univ}_{\lambda}, \nabla_{\lambda}) \simeq \DR_{\lambda}(\sh{V}^{\univ}_{\lambda}, t\nabla_{\lambda})\]
  thus $J^i_k(K_{\Hod}(\sh{V}^{\univ}))$ is $\GG_m$-invariant, hence so is $\Sigma^i_{k, \Hod}(X, r)$. It follows that 
  \[\Theta(J^i_k(K_{\Hod})) \subseteq J^i_k(K_{\Hod})\]
  and completing along the image of $s_{(\bb{V}, \phi)}$ gives $\Bf{(H2)}$ by the identification in the previous section. \\
  \\
  \textit{Closed under semisimplification and absoluteness}. Closure under semisimplification follows from semicontinuity of local system cohomology. Asboluteness follows from the same proof as in \cite[Theorem 11.1.1]{budur2020absolutesets}. Finally, irreducible components and reduction remains Hodge, due to \cite[Proposition 8.21(3)]{bakker2024linear}, and absolute. 
\end{proof}

\begin{coroll}
  \label{lem:proj=torus}
  Using the notation of \cref{sec:jump-product}, let $Z$ be an irreducible component of $\Sigma^i_k(X, 2) \subseteq R_B(X, x, 2)$. Let $C$ be an irreducible component of $\phi_{\chi}^{-1}(Z) \subset R_B(X, 2, \det = \chi) \times A$, and denote 
  \[T \coloneq \ol{\Rm{pr}_2(C)}_{\red}\]
  the Zariski closure of the reduced image of the projection onto $A = M_B(X, 1)$. Then $T$ is a torsion-translated subtorus of $A$. \\
  \tab Furthermore, if $\bb{L} \in T$ is a unitary point, then $\widehat{\sh{O}}_{T, t}$ underlies a pro-$\Bf{T}(0)$-MHS-algebra such that for every $(\bb{V}, \bb{L}) \in C$ underlying $\CC$-VHS, $\widehat{\sh{O}}_{T, \bb{L}} \to \widehat{\sh{O}}_{C, (\bb{V}, \bb{L})}$ is a morphism of pro-$\Bf{T}(0)$-MHS-algebras. 
\end{coroll}
\begin{proof}
  By \cite[Proposition 9.10]{bakker2024linear}, all the morphisms in 
  \[
  \begin{tikzcd}
	{R_B(X, 2, \det = \chi) \times A} && {R_B(X, 2, \det \in \chi A^2)} \\
	\\
	A && A
	\arrow["{(\bb{V}, \bb{L}) \mapsto \bb{V} \otimes \bb{L}}", "{\Phi_{\chi}}"', from=1-1, to=1-3]
	\arrow["{\Rm{pr}_2}"', from=1-1, to=3-1]
	\arrow["\det", from=1-3, to=3-3]
	\arrow["{\bb{L} \mapsto \chi \otimes \bb{L}^{\otimes 2}}", from=3-1, to=3-3]
  \end{tikzcd}
  \]
  are absolute $\ol{\QQ}$-morphism \cite[Definition 9.9]{bakker2024linear}. By the previous theorem, $Z$ is absolute $\ol{\QQ}$-closed. Thus $\phi_{\chi}^{-1}(Z)$ is absolute $\ol{\QQ}$-closed, and $\Rm{pr}_2(C)$ is absolute $\ol{\QQ}$-constructible. By \cite[Proposition 9.14]{bakker2024linear}, $T$ is absolute $\ol{\QQ}$-closed. Thus by \cref{thm:rank1-absolute}, $T = \tau \cdot N$ where $N$ is a subtorus of $A$, and $\Rm{Lie(N)}$ is a complex sub Hodge structure of $H^1(X, \CC)$.\\
  \tab Now, suppose $\bb{L} \in T$ is unitary, then $\bb{L}$ underlies a canonical weight 0 $\CC$-VHS, and by \cref{thm:Betti-MTS}, $\widehat{\sh{O}}_{A, \bb{L}}$ underlies a pro-$\Bf{T}(0)$-MHS-algebra. Furthermore, since $A$ and $T$ are smooth, $\widehat{\sh{O}}_{T, \bb{L}}$ and $\widehat{\sh{O}}_{A, \bb{L}}$ are determined by tangent spaces, i.e., 
  \[\widehat{\sh{O}}_{A, \bb{L}} \simeq \widehat{\Rm{Sym}}\ H^1(X, \CC)^{\vee}, \tab \widehat{\sh{O}}_{T, \bb{L}} \simeq \widehat{\Rm{Sym}}\ \Rm{Lie}(N)^{\vee}\]
  and $\widehat{\sh{O}}_{A, \bb{L}} \to \widehat{\sh{O}}_{T, \bb{L}}$ is a quotient of pro-$\Bf{T}(0)$-MHS-algebras. Furthermore, $H^1(X, \CC)$ is pure of weight 1, hence so is $\Rm{Lie}(N)$ and we get $\Ide{m}_{T, \bb{L}} = W_{-1}\widehat{\sh{O}}_{T, \bb{L}}$ with $\Ide{m}_{T, \bb{L}}/\Ide{m}_{T, \bb{L}}^2$ is pure of weight $-1$. Suppose that $(\bb{V}, \bb{L}) \in C$ is a $\CC$-VHS point, equivalently $\bb{V}$ underlies a weight 0 $\CC$-VHS, then 
  \[\widehat{\sh{O}}_{R_B(X, 2, \det = \chi), \bb{V}}\ \widehat{\otimes}\ \widehat{\sh{O}}_{A, \bb{L}} \twoheadrightarrow \widehat{\sh{O}}_{C, (\bb{V}, \bb{L})}\]
  is a quotient of pro-$\Bf{T}(0)$-MHS-algebras. The composition
  \[\widehat{\sh{O}}_{A, \bb{L}} \to \widehat{\sh{O}}_{R_B(X, 2, \det = \chi), \bb{V}}\ \widehat{\otimes}\ \widehat{\sh{O}}_{A, \bb{L}} \twoheadrightarrow \widehat{\sh{O}}_{C, (\bb{V}, \bb{L})}\]
  is thus a morphism of pro-$\Bf{T}(0)$-MHS-algebras. This factors through $\widehat{\sh{O}}_{A, \bb{L}} \to \widehat{\sh{O}}_{T, \bb{L}}$ which is a quotient of pro-$\Bf{T}(0)$-MHS-algebras, hence the factorization 
  \[\widehat{\sh{O}}_{T, \bb{L}} \to \widehat{\sh{O}}_{C, (\bb{V}, \bb{L})}\]
  is a morphism of pro-$\Bf{T}(0)$-MHS-algebras, as desired. 
\end{proof}

\subsection{Hodge structures on good moduli image}
\label{sec:GIT-Hodge} For a Hodge subvariety $W \subseteq R_B(X, x, r)$, at a point $(\bb{V}, \phi) \in W$ with $\bb{V}$ underlying a canonical weight 0 complex variation of Hodge structures, we have a pro-$\Bf{T}(0)$-MHS-algebra structure on $\widehat{\sh{O}}_{W, (\bb{V}, \phi)}$. We prove a similar statement for the image of $W$ under the GIT quotient $q: R_B(X, x, r) \to M_B(X, r) = R_B(X, x, r)//\Rm{GL}_2$. 
\begin{prop}
  \label{thm:GIT-Hodge}
  Let $W \subseteq R_B(X, x, r)$ be a Hodge subvariety, and $Z = q(W) \subset M_B(X, r)$ be the scheme-theoretic image. Let $\bb{V} \in Z$ be a $\CC^*$-fixed point, hence underlying a weight 0 complex variation of Hodge structures. Then $\widehat{\sh{O}}_{Z, \bb{V}}$ underlies a pro-$\Bf{T}(0)$-MHS-algebra, with 
  \[\Ide{m}_{Z, \bb{V}} = W_{-1}\widehat{\sh{O}}_{Z, \bb{V}}\]
\end{prop}
\begin{proof}
  The stabilizer of $(\bb{V}, \phi) \in R_B(X, x, r)$ under the $\Rm{GL}_r$-action is $\Rm{Aut}(\bb{V})$. Denote $H_{\bb{V}} = \Rm{Aut}(\bb{V}) \subset \Rm{GL}_r$. We have an isotypic decomposition 
  \[\bb{V} \simeq \bigoplus \bb{V}_{\alpha} \otimes M_\alpha\]
  where $\set*{\bb{V}_{\alpha}}$ are pairwise non-isomorphic irreducible local systems, and $\set*{M_{\alpha}}$ are complex vector spaces. By equipping each $\bb{V}_{\alpha}$ with a weight 0 $\CC$-VHS and $M_{\alpha}$ with a pure complex Hodge structure of type $(0, 0)$, we get an equality of $\CC$-VHS. Then by Schur's lemma, $H_{\bb{V}} = \prod \Rm{GL}(M_{\alpha})$, hence acts by weight 0 automorphisms of the $\CC$-VHS on $\bb{V}$. 

  \medskip
  On the other hand, by picking a framing we get $(\bb{V}, \phi) \in W$ lying over $\bb{V}$. Recall that included in the definition of Hodge subvariety is the fact that $W$ is $\Rm{GL}_r$-invariant. Since $\bb{V}$ is semisimple, the $\Rm{GL}_r$-orbit of $(\bb{V}, \phi)$ is closed. Then by the Luna \'{e}tale slice theorem (see \cite{luna1973slices} or \cite[Remark 4.6]{alper2020luna}) there exists a $H_{\bb{V}}$-invariant affine subscheme $Y \subseteq R_B(X, x, r)$ containing $(\bb{V}, \phi)$ such that $\Rm{GL}_r \times^{H_{\bb{V}}} Y \to R_B(X, x, r)$ is \'{e}tale; its image $U \subseteq R_B(X, x, r)$ is a $\Rm{GL}_r$-invariant saturated (i.e., $U = q^{-1} \circ q(U)$) affine open. Here 
  \[\Rm{GL}_r \times^{H_{\bb{V}}} Y = (\Rm{GL}_r \times Y)/H_{\bb{V}}\]
  where the $H_{\bb{V}}$-action is $h \cdot (g, y) = (g \cdot h, h^{-1} \cdot y)$. Furthermore, the diagram 
  \[
  \begin{tikzcd}
	{\Rm{GL}_r \times^{H_{\bb{V}}} Y} && U \\
	\\
	{Y//H_{\bb{V}}} && {U//\Rm{GL}_r}
	\arrow[from=1-1, to=1-3]
	\arrow[from=1-1, to=3-1]
	\arrow[from=1-3, to=3-3]
	\arrow[from=3-1, to=3-3]
\end{tikzcd}\]
  is cartesian with two horizontal maps \'{e}tale. This induces the following cartesian diagram of quotient stacks by restricting to $\Rm{GL}_r$-invariant $W$
  \[\begin{tikzcd}
	{[(W \cap Y)/H_{\bb{V}}]} && {[W/\Rm{GL}_r]} \\
	\\
	{(W \cap Y)//H_{\bb{V}}} && Z
	\arrow[from=1-1, to=1-3]
	\arrow[from=1-1, to=3-1]
	\arrow[from=1-3, to=3-3]
	\arrow[from=3-1, to=3-3]
  \end{tikzcd}
  \]
  and \'{e}taleness of the bottom horizontal map then gives an isomorphism $\widehat{\sh{O}}_{Z, \bb{V}} \simeq \widehat{\sh{O}}_{W \cap Y, (\bb{V}, \phi)}^{H_{\bb{V}}}$. Now, notice that $\widehat{\sh{O}}_{Y, (\bb{V},\phi)}$ is a miniversal algebra for $\sh{M}_B$ \cite[Remark 4.17]{alper2020luna}, and
  \[[W/\Rm{GL}_r] \times_{\sh{M}_B} Y \simeq [W/\Rm{GL}_r] \times_{\sh{M}_B} R_B \times_{R_B} Y \simeq W \cap Y\]
  \tab It follows that the map $\widehat{\sh{O}}_{Y, (\bb{V},\phi)} \to \widehat{\sh{O}}_{W \cap Y, (\bb{V}, \phi)}$ is the base change of the miniversal family $\Rm{Spf}\ \widehat{\sh{O}}_{Y, (\bb{V},\phi)} \to \sh{M}_B(X, r)$ along $[W/\Rm{GL}_r] \hookrightarrow \sh{M}_B$. Then, by \cite[Lemma 4.59]{bakker2024linear} and the fact that $W$ is Hodge at $(\bb{V}, \phi)$, $\widehat{\sh{O}}_{Y, (\bb{V},\phi)} \to \widehat{\sh{O}}_{W \cap Y, (\bb{V}, \phi)}$ must be a quotient of pro-$\Bf{T}(0)$-MHS-algebras; the pro-$\Bf{T}(0)$-MHS-algebra on the miniversal $\widehat{\sh{O}}_{Y, (\bb{V},\phi)}$ is given by \cite[Theorem 4.3]{bakker2024linear}. 
  \medskip
  Let $e \in H_{\bb{V}}$ be the identity. We have a diagram
  \begin{equation}
  \begin{tikzcd}
	{\Rm{Spf}\ \widehat{\sh{O}}_{H_{\bb{V}}, e} \times \Rm{Spf}\ \widehat{\sh{O}}_{W \cap Y, (\bb{V}, \phi)}} && \\
	\\
	{H_{\bb{V}} \times\Rm{Spf}\ \widehat{\sh{O}}_{W \cap Y, (\bb{V}, \phi)}} && {\Rm{Spf}\ \widehat{\sh{O}}_{W \cap Y, (\bb{V}, \phi)} \times \Rm{Spf}\ \widehat{\sh{O}}_{W \cap Y, (\bb{V}, \phi)}} \\
	\\
	{\sh{M}_B(X, r)} && {\sh{M}_B(X, r) \times \sh{M}_B(X, r)}
	\arrow["{\text{completion}}"{description}, dashed, from=1-1, to=3-1]
	\arrow["{(s, t)}", from=1-1, to=3-3]
	\arrow[from=3-1, to=3-3]
	\arrow[from=3-1, to=5-1]
	\arrow[from=3-3, to=5-3]
	\arrow["\Delta", from=5-1, to=5-3]
\end{tikzcd}
  \end{equation}
  and by \cite[Theorem 4.3(3)]{bakker2024linear} there is a pro-$\Bf{T}(0)$-MHS-algebra structure on $\Rm{Spf}\ \widehat{\sh{O}}_{H_{\bb{V}}, e} \times \Rm{Spf}\ \widehat{\sh{O}}_{W \cap Y, (\bb{V}, \phi)}$ making $(s^*, t^*)$ (the ring maps) a morphism of pro-$\Bf{T}(0)$-MHS-algebras. Notice that $(s, t)$ gives the formal groupoid action of $\widehat{H}_{\bb{V}, e} = \Rm{Spf}\ \widehat{\sh{O}}_{H_{\bb{V}}, e}$ on $\Rm{Spf}\ \widehat{\sh{O}}_{W \cap Y, (\bb{V}, \phi)}$, hence (or alternatively using \cite[\textsection 4]{Alper2012}) we have 
  \begin{align*}
    \widehat{\sh{O}}_{W \cap Y, (\bb{V}, \phi)}^{H_{\bb{V}}} &\simeq \widehat{\sh{O}}_{W \cap Y, (\bb{V}, \phi)}^{\widehat{H}_{\bb{V}, e}} \\
    &\simeq \ker (s^* - t^*: \widehat{\sh{O}}_{W \cap Y, (\bb{V}, \phi)} \to \widehat{\sh{O}}_{W \cap Y, (\bb{V}, \phi)} \otimes \widehat{\sh{O}}_{H_{\bb{V}}, e})
  \end{align*}
  \tab Since $s^*, t^*$ are morphisms of pro-$\Bf{T}(0)$-algebras with $\Sf{MHS}_{\CC}$ abelian, we get an induced pro-$\Bf{T}(0)$-MHS-algebra structure on $\widehat{\sh{O}}_{Z, \bb{V}} \simeq \widehat{\sh{O}}_{W \cap Y, (\bb{V}, \phi)}^{H_{\bb{V}}}$. By \cite[Theorem 4.3]{bakker2024linear} again, 
  \[\Ide{m}_{Y, (\bb{V}, \phi)} = W_{-1} \widehat{\sh{O}}_{Y, (\bb{V}, \phi)}\]
  hence the same is true for the quotient $\widehat{\sh{O}}_{W \cap Y, (\bb{V}, \phi)}$. Finally, due to strictness of $s^*, t^*$ as morphisms of pro-$\Bf{T}(0)$-MHS-algebras, we get 
  \begin{align*}
    W_{-1}\widehat{\sh{O}}_{Z, \bb{V}} \simeq W_{-1}\widehat{\sh{O}}_{W \cap Y, (\bb{V}, \phi)}^{H_{\bb{V}}} &\simeq \widehat{\sh{O}}_{W \cap Y, (\bb{V}, \phi)}^{H_{\bb{V}}} \cap \Ide{m}_{\widehat{\sh{O}}_{W \cap Y, (\bb{V}, \phi)}} \\
    &= \Ide{m}_{\widehat{\sh{O}}_{W \cap Y, (\bb{V}, \phi)}^{H_{\bb{V}}}} = \Ide{m}_{\widehat{\sh{O}}_{Z, \bb{V}}}
  \end{align*}
  as desired. 
\end{proof}
The same proof gives the following twistor version:
\begin{coroll}
  \label{lem:GIT-twistor}
  Let $W \subseteq R_B(X, x, r)$ be a Hodge subvariety, and $Z = q(W) \subseteq M_B(X, r)$. Let $\bb{V} \in Z$, then $\widehat{\sh{O}}_{Z, \bb{V}}$ underlies a pro-$\Bf{T}(0)$-MTS-algebra with 
  \[\Ide{m}_{Z, \bb{V}} = W_{-1}\widehat{\sh{O}}_{Z, \bb{V}}\]
\end{coroll}
With this pro-$\Bf{T}(0)$-MHS-algebra structure on the good moduli space $M_B(X, r)$, we have a version of \cref{lem:proj=torus} for $M_B(X, 2)$ instead of the framed space $R_B(X, x, 2)$, as follows. 
\begin{coroll}
  \label{lem:GIT-proj=torus}
  Using the notation of \cref{sec:jump-product}, let $Z$ be an irreducible component of the reduced image $\ol{\Sigma}^i_k(X, 2) \subseteq M_B(X, 2)$. Let $C$ be an irreducible component of $\phi_{\chi}^{-1}(Z) \subset M_B(X, 2, \det = \chi) \times A$, and denote 
  \[T \coloneq \ol{\Rm{pr}_2(C)}_{\red}\]
  the Zariski closure of the reduced image of the projection onto $A = M_B(X, 1)$. Then $T$ is a torsion-translated subtorus of $A$. \\
  \tab Furthermore, if $\bb{L} \in T$ is a unitary point, then $\widehat{\sh{O}}_{T, t}$ underlies a pro-$\Bf{T}(0)$-MHS-algebra such that for every $(\bb{V}, \bb{L}) \in C$ underlying $\CC$-VHS, $\widehat{\sh{O}}_{T, \bb{L}} \to \widehat{\sh{O}}_{C, (\bb{V}, \bb{L})}$ is a morphism of pro-$\Bf{T}(0)$-MHS-algebras. 
\end{coroll}
\begin{proof}
  Firstly, the image of $\Rm{pr}_2$ is the same whether we work with $M_B(X, 2)$ or $R_B(X, x, 2)$, hence $T$ is still a torsion-translated subtorus of $A$. Secondly, every such irreducible component $Z$ is the good moduli image of some irreducible component $W$ of $\Sigma^i_k(X, 2) \subseteq R_B(X, x, 2)$. Next, since the miniversal pro-$\Bf{T}(0)$-MHS-algebra structure of \cite[Theorem 4.3]{bakker2024linear} is compatible with tensor product, and that the proof of \cref{thm:GIT-Hodge} is clearly functorial with respect to the inclusion of Hodge subvarieties $W \subseteq R_B(X, x, r)$, we get that 
  \[\widehat{\sh{O}}_{M_B(X, 2, \det = \chi), \bb{V}}\ \widehat{\otimes}\ \widehat{\sh{O}}_{A, \bb{L}} \to \widehat{\sh{O}}_{C, (\bb{V}, \bb{L})}\]
  is a quotient of pro-$\Bf{T}(0)$-MHS-algebras for when $\bb{V}, \bb{L}$ underlie weight 0 complex variation of Hodge structures. Then the same proof as in \cref{lem:proj=torus} works. 
\end{proof}

\subsection{Hodge structures on normalization}
\label{sec:norm-Hodge}
The goal of this section is to show that if we have a normalization $C^{\nu} \to C$, and $(\bb{V}, \bb{L})$ is a $\CC$-VHS point of $C$, then for every $c \in C^{\nu}$ lying over $(\bb{V}, \bb{L})$ there is a pro-$\Bf{T}(0)$-MHS-algebra structure on $\widehat{\sh{O}}_{C^{\nu}, c}$ as well. We will start with some algebra facts. 

\begin{lemma}
  \label{lem:finite-weight}
  Let $\phi: (A, \Ide{m}_A)\to (B, \Ide{m}_B)$ be a finite homomorphism of pro-$\Bf{T}(0)$-MHS-algebras, with $(A, \Ide{m}_A)$ and $(B, \Ide{m}_B)$ complete local domains. Assume further that 
  \[\Ide{m}_A = W_{-1}A, \tab W_0A = A\]
  then $W_0B = B$ and $\Ide{m}_B = W_{-1}B$ as well. 
\end{lemma}
\begin{proof}
  We have the functorial Deligne splitting \cite[Lemma-Definition 3.4]{PetersSteenbrink2008} for the weight filtrations on $A$ on $B$, giving $A \simeq \Prod{w \in \ZZ_{\leq 0}} A_w, B \simeq \Prod{w \in \ZZ} B_w$, and a finite map 
  \[\bar{\phi}: \Prod{w \in \ZZ_{\leq 0}} A_w \to \Prod{w \in \ZZ} B_w\]
  compatible with the weights. Let $0 \neq b \in B_w$ be a homogeneous element, then $b$ is integral, say of degree $d$, over $\Prod{w \in \ZZ_{\leq 0}} A_w$. Take the $dw$-homogeneous part of the monic polynomial for $b$, we get 
  \[b^d + a_1b^{d - 1} + \dots + a_d = 0\]
  with $a_j \in A_{jw}$. If $w > 0$ then $A_{jw} = 0$ hence $a_j = 0$ and $b^d = 0$ which gives a contradiction since $A$ is a domain. If $w = 0$ then all $a_j \in A_0 = \CC$ hence $b \in \CC$. It follows that $W_{-1}B = \Ide{m}_B$ and $W_0B = B$. 
\end{proof}

\begin{lemma}
  \label{lem:reduced-sing-ideal}
  Let $B$ be a complete noetherian local pro-$\Bf{T}(0)$-MHS-algebra with $W_{-1}B = \Ide{m}_B$. Then the ideal of the reduced singular locus $I(\Rm{Sing}(\Spec B)_{\red})$ is a pro-$\Bf{T}(0)$-MHS-ideal of $B$. 
\end{lemma}
\begin{proof}
  This follows from the proof of \cite[Lemma 8.23]{bakker2024linear} which implies that $I(\Rm{Sing}(\Spec B))$ is a pro-$\Bf{T}(0)$-MHS-ideal of $B$, and \cite[Lemma 4.36]{bakker2024linear} which says that taking the reduction preserves the Hodge structure. 
\end{proof}

\begin{lemma}
  \label{lem:end-Hodge}
  Let $(B, \Ide{m}_B)$ be a pro-$\Bf{T}(0)$-MHS-algebra which is a complete noetherian local domain, and $I \subset B$ a pro-$\Bf{T}(0)$-MHS-ideal. Then $\Rm{End}_B(I)$ underlies a local pro-$\Bf{T}(0)$-MHS-algebra structure such that the natural map defined by multiplication 
  \[B \to \Rm{End}_B(I)\]
  is an finite injective morphism of pro-$\Bf{T}(0)$-MHS-algebras. 
\end{lemma}
\begin{proof}
  Write 
  \[B = \varprojlim B/\Ide{m}_B^{n + 1} = \varprojlim B_n, \tab I = \varprojlim I/\Ide{m}_B^{n + 1}I = \varprojlim I_n\]
  then each $\Rm{End}_{B_n}(I_n)$ is just the kernel of the morphism of $\CC$-MHS
  \[\Rm{End}_{\CC}(I_n) \to \Hom_{\CC}(B_n \otimes_{\CC} I_n, I_n)\]
  hence carries the structure of an associative algebra object on the category of $\Bf{T}(0)$-MHS-algebra. Next, we have tensor-Hom adjunction
  \[\Hom_{\Sf{MHS}_{\CC}}(B_n \otimes I_n, I_n) \simeq \Hom_{\Sf{MHS}_{\CC}}(B_n, \Rm{End}_{\CC}(I_n))\] 
  under which the multiplication map $B_n \to \Rm{End}_{\CC}(I_n)$ corresponds to the $B_n$-module structure map $B_n \otimes I_n \to I_n$ which is a morphism of $\CC$-MHS since $I_n$ is a $\Bf{T}(0)-$MHS-ideal of $B_n$. Hence $B_n \to \Rm{End}_{\CC}(I_n)$ is a morphism of $\CC$-MHS which factors through $B_n \to \Rm{End}_{B_n}(I_n)$ in $\Sf{MHS}_{\CC}$.

  \medskip
  The map $\Rm{End}_B(I) \to \Rm{End}_{B_n}(I_n)$ is not surjective, and with kernel 
  \[K_n = \set*{f \in \Rm{End}_B(I)\ \big\vert\ f(I) \subseteq \Ide{m}_B^{n + 1} I}\]
  \tab Let $\basis*{x_1, \dots, x_l}$ be a generating set of $I$ as a $B$-module. Then the map 
  \[\Rm{End}_B(I) \to I^l, \tab f \mapsto \brac*{f(x_1, x_2, \dots, x_l)}\]
  is an injection, and gives $K_n = \Rm{End}_B(I) \cap \Ide{m}_B^{n + 1}I^l$. By \cite[Theorem 10.11]{atiyah1969introduction}, the two filtrations $\Ide{m}_B^{\bullet + 1} \Rm{End}_B(I)$ and $K_{\bullet}$ are cofinal. Thus we can take the inverse limit of the first part of the proof to get the desired statements. The kernel is $\Rm{Ann}_{B}(I) = 0$ since $B$ is a domain. 
\end{proof}

\begin{prop}
  \label{thm:norm-Hodge}
  Let $\Lambda$ be a reduced complete noetherian local $\CC$-algebra which underlies a pro-$\Bf{T}(0)$-MHS-algebra such that $\Ide{m}_{\Lambda} = W_{-1}\Lambda$. For every minimal prime $\Ide{p} \subset \Lambda$, denote $A = \Lambda/\Ide{p}$ and $A^{\nu}$ its normalization. Then $A^{\nu}$ underlies a pro-$\Bf{T}(0)$-MHS-algebra structure with 
  \[W_0A^{\nu} = A^{\nu}, \tab W_{-1}A^{\nu} = \Ide{m}_{\nu}\]
  such that the normalization map $\nu: A \to A^{\nu}$ and its composition $\Lambda \to A^{\nu}$ are morphisms of pro-$\Bf{T}(0)$-MHS-algebras. 
\end{prop}
\begin{proof}
  We will make use Oka's criterion for normality \cite[Chapter 6 \textsection 5.2]{grauert1984coherent}, which says that for a reduced ring $B$, $B$ is normal if and only if the multiplication map 
  \[B \to \Rm{End}_B(I(\Rm{Sing}(\Spec B)_{\red}))\]
  is an isomorphism. Start with $B^0 = A = \Lambda/\Ide{p}$ which is a domain, and define recursively 
  \[B^{j + 1} = \Rm{End}_{B^j}(I(\Rm{Sing}(\Spec B^j)_{\red}))\]
  then each $B^{j + 1} \subseteq A^{\nu}$ and is finite . We will prove this by induction. Denote $I^j = I(\Rm{Sing}(\Spec B^j)_{\red})$. Suppose $B^j \subseteq A^{\nu}$ is a domain and $\Rm{Frac}(B_j) = \Rm{Frac}(A)$, which is automatically satisfied for $j = 0$. Pick $0 \neq b \in I^j$, then for 
  \[f \in \Rm{End}_{B^j}(I^j)\] 
  $f(b)/b \in \Rm{Frac}(A)$. We get that $f(b') = \frac{f(b)}{b} b'$ for all $b' \in I^j$, which gives an embedding $\Rm{End}_{B^j}(I^j) \hookrightarrow \Rm{Frac}(A)$ by $f \mapsto f(b)/b$. Thus we get 
  \[B^j \subseteq \Rm{End}_{B^j}(I^j) = B^{j + 1} \subseteq \Rm{Frac}(A)\]
  where the first inclusion is by the multiplication map and the fact that $B^j$ is a domain. Finally $B^{j + 1}$ is finite over $B^j$ hence $B^{j + 1} \subseteq A^{\nu}$ as well. We have a chain 
  \[A \subseteq B^1 \subseteq B^2 \subseteq \dots \subseteq A^{\nu}\]
  which must stabilize due to $A^{\nu}$ being finite over $A$. At the stable stage we get $B^j \simeq \Rm{End}_{B^j}(I^j)$ hence $B^j = (B^j)^{\nu} = A^{\nu}$. 

  \medskip
  By \cite[Lemma 4.36]{bakker2024linear}, $\Lambda \to A = \Lambda/\Ide{p}$ is a morphism of pro-$\Bf{T}(0)$-MHS-algebras. By \cref{lem:reduced-sing-ideal} and \cref{lem:end-Hodge}, the inclusions 
  \[A \subseteq B^1 \subseteq B^2 \subseteq \dots \subseteq A^{\nu}\]
  are all morphisms of pro-$\Bf{T}(0)$-MHS-algebras. It follows that $A \to A^{\nu}$ is a morphisms of pro-$\Bf{T}(0)$-MHS-algebras, which is a finite map hence by \cref{lem:finite-weight} we get the weight statements. 
\end{proof}

\begin{coroll}
  \label{lem:normal-Hodge}
  Using the notation of \cref{sec:jump-product}, let $Z$ be an irreducible component of the reduced image $\ol{\Sigma}^i_k(X, 2) \subseteq M_B(X, 2)$. Let $C$ be an irreducible component of $\phi_{\chi}^{-1}(Z) \subset M_B(X, 2, \det = \chi) \times A$, and $(\bb{V}, \bb{L}) \in C$ a $\CC$-VHS point. Consider the normalization $\nu: C^{\nu} \to C$ with $c$ lying over $(\bb{V}, \bb{L})$, then $\widehat{\sh{O}}_{C^{\nu}, c}$ underlies a pro-$\Bf{T}(0)$-MHS-algebra with 
  \[W_0\widehat{\sh{O}}_{C^{\nu}, c} = \widehat{\sh{O}}_{C^{\nu}, c}, \tab W_{-1}\widehat{\sh{O}}_{C^{\nu}, c} = \Ide{m}_{C^{\nu}, c}\]
  such that the map $\widehat{\sh{O}}_{C, (\bb{V}, \bb{L})} \to \widehat{\sh{O}}_{C^{\nu}, c}$ is a morphism of pro-$\Bf{T}(0)$-MHS-algebras. 
\end{coroll}
\begin{proof}
  Any such $\widehat{\sh{O}}_{C^{\nu}, c}$ would be isomorphic to $\brac*{\widehat{\sh{O}}_{C, (\bb{V}, \bb{L})}/\Ide{p}}^{\nu}$ for some minimal prime $\Ide{p}$. \cref{thm:norm-Hodge} then gives the desired statements. 
\end{proof}

\section{Projectively virtually abelian and projectively rigid jump components}
\label{sec:abelian-comp}
Let $\Sigma^i_k(X, 2) \subseteq R_B(X, x, 2)$ be a jump locus, and $\ol{\Sigma}^i_k(X, 2) \subseteq M_B(X, 2)$ be the reduced good moduli image. Denote $Z \subseteq \ol{\Sigma}^i_k(X, 2)$ an irreducible component. We first justify labeling $Z$ based on its general point. 

\subsection{Generic projective type} Consider the following algebraic subgroups of $\Rm{PGL}_2$:
\begin{align}
  B &= \Rm{Stab}(\infty) = \set*{\sqbrac*{\begin{pmatrix}
    a & b\\
    0 & d
  \end{pmatrix}
  }\ \bigg\vert\ ad \neq 0}\\
  T_{0, \infty} &= \Rm{Stab}(\infty) \cap \Rm{Stab}(0) = \set*{\sqbrac*{\begin{pmatrix}
    a & 0\\
    0 & d
  \end{pmatrix}
  }\ \bigg\vert\ ad \neq 0}\\
  D_{\infty} &= T_{0, \infty} \cup \set*{\sqbrac*{\begin{pmatrix}
    0 & b\\
    c & 0
  \end{pmatrix}
  }\ \bigg\vert\ bc \neq 0}
\end{align}
where square bracket denotes the projective class. Notice that $T_{0, \infty}$ is an algebraic torus, $D_{\infty}$ is the normalizer of $T_{0, \infty}$, fixes the set of points $\set*{0, \infty} \subset \PP^1$ and is $D_{\infty} \simeq T_{0, \infty} \rtimes \ZZ/2\ZZ$.  

Let $\rho: \pi_1(X, x) \to \Rm{GL}_2(\CC)$ be a semisimple representation, and $\ol{\rho}: \pi_1(X, x) \to \Rm{PGL}_2(\CC)$ its projective quotient. Denote the Zariski closure
\[G_{\ol{\rho}} \coloneq \ol{\ol{\rho}(\pi_1(X, x))}^{\Rm{Zar}} \subseteq \Rm{PGL}_2\]
the algebraic monodromy group of $\ol{\rho}$. 
\begin{lemma}
  \label{lem:rep-type}
  For each semisimple $\rho: \pi_1(X, x) \to \Rm{GL}_2(\CC)$, exactly one of the following is true:
  \begin{description}
    \item[(A)] $\rho$ is reducible, i.e., $\rho \simeq \chi_1 \oplus \chi_2$;
    \item[(B)] $\rho$ is irreducible, and, up to conjugation, $G_{\ol{\rho}} \subseteq D_{\infty}$;
    \item[(C)] $G_{\ol{\rho}} \simeq \ol{\rho}(\pi_1(X, x)) \simeq A_4, S_4$ or $A_5$;
    \item[(D)] $\ol{\rho}$ has Zariski dense image and is rigid;
    \item[(E)] $\ol{\rho}$ has Zariski dense image and is non-rigid.
  \end{description}
\end{lemma}
\begin{proof}
  $\rho$ is reducible if and only if it fixes a line, or equivalently, $G_{\ol{\rho}}$ fixes a point in $\PP^1$. Up to conjugation, this is equivalent to $G_{\ol{\rho}}$ being contained in the Borel $B$. The fact that this classification is exhaustive follows from \cite[Theorem 1]{nguyen2008algebraic}. 

  For mutual exclusivity, it is clear that \textbf{(A)}, \textbf{(B)}, \textbf{(D)}, and \textbf{(E)} are mutually exclusive; same for \textbf{(C)}, \textbf{(D)} and \textbf{(E)}. For \textbf{(B)}, and \textbf{(C)}, notice that the only finite subgroups of $D_{\infty} = T_{0, \infty} \rtimes \ZZ/2\ZZ$ is either dihedral $D_n$ or a cyclic group. Finally for \textbf{(A)} and \textbf{(E)}, a finite subgroup fixing a point must fix another (it acts as an affine transformation on $\AAA^1$), hence, up to conjugation, is contained in the torus $T_{0, \infty}$. It follows that it must be cyclic. 
\end{proof}
We will show that a general point of an irreducible component must be of exactly one of the previous 5 types. In other words, we cannot have two types whose points are both dense in the component. It suffices to show that each type is a locally closed condition. 
\begin{definition}
  \label{def:jump-labels}
  Denote 
  \begin{align*}
    \Ide{D} &\subset M_B(X, 2) \text{ be } \bigcup_{0 \neq \epsilon \in \Hom(\pi_1, \ZZ/2\ZZ)} \Ide{D}_{\epsilon} \text{ as in \cref{lem:dihedral-locus},}\\
    \Ide{R} &\subset M_B(X, 2) \text{ the locus of reducible representations,}\\
    \Ide{F} &\subset M_B(X, 2) \text{ the locus of representations of type \textbf{(C)},}\\
    \Ide{G} &\subset M_B(X, 2) \text{ the locus of projectively dense representations.}
  \end{align*}
  then $M_B(X, 2) = \Ide{R} \cup \Ide{D} \cup \Ide{F} \cup \Ide{G}$. 
\end{definition}
\begin{prop}
  \label{lem:jump-labels}
  Let $Z$ be an irreducible component of $\ol{\Sigma}^i_k(X, 2)$, the reduced good moduli image of $\Sigma^i_k(X, 2) \subseteq R_B(X, x, 2)$. Then $Z$ has a nonempty Zariski open subset on which exactly one of the types in \cref{lem:rep-type} occurs. If the type is $\mathbf{(A)}, \mathbf{(B)}$, or $\mathbf{(C)}$, then $Z$ is contained in, respectively, $\Ide{R}$, $\Ide{D}$ or $\Ide{F}$. 
\end{prop}
\begin{proof}
  The locus of type $\mathbf{(A)}$ is $\Ide{R}$, which is closed by \cref{lem:reducible-locus}. The locus of type $\mathbf{(C)}$ is closed by \eqref{eq:finite-locus}. The locus of type $\mathbf{(B)}$ is open in $\Ide{D}$ by \cref{lem:dihedral-locus}. Type $\mathbf{(D)}$ locus is closed in $\Ide{G}$, since it is a finite union of images of $\theta_{\rho}$ as defined in \cref{lem:twist-map}; this is because there are only finitely many rigid points in $M_B(X, \Rm{PGL}_2)$. Finally the locus of type $\mathbf{(E)}$ is open in $\Ide{G}$ hence open. 
\end{proof}

\subsection{Reducible jump component}
\label{lem:reducible-jump}
\begin{lemma}
  \label{lem:reducible-locus}
  Let $A = M_B(X, 1)$, then the map 
  \begin{equation}
    \iota: (A \times A)/S_2 \to M_B(X, 2) \label{eq:reducible-embed}, \tab (\chi_1, \chi_2) \mapsto [\chi_1 \oplus \chi_2]
  \end{equation}
  is a closed embedding, with image the reducible locus $\Ide{R}$. 
\end{lemma}
\begin{proof}
  These are both affine schemes, and the corresponding map of rings is 
  \[\iota^{\#}: H^0(R_B(X, x, 2), \sh{O}_{R_B})^{\Rm{GL}_2} \to (H^0(A, \sh{O}_A) \otimes H^0(A, \sh{O}_A))^{S_2}\]
  mapping 
  \[f \mapsto \sqbrac*{(\chi_1, \chi_2) \mapsto f\brac*{\gamma \mapsto \begin{pmatrix}
    \chi_1(\gamma) & 0 \\
    0 & \chi_2(\gamma)
  \end{pmatrix}}}
  \]
  where $f$ is a function on $\Hom(\pi_1(X, x), \Rm{GL}_2)$. Notice that $H^0(A, \sh{O}_A) \otimes H^0(A, \sh{O}_A)$ is generated by the functions 
  \[\gamma^1: (\chi_1, \chi_2) \mapsto \chi_1(\gamma), \tab \gamma^2: (\chi_1, \chi_2) \mapsto \chi_2(\gamma) \tab \tab \gamma \in \pi_1(X, x)\]
  i.e., every regular function on $A \times A$ has the form $\sum \gamma^1 \delta^2$ with $\gamma, \delta \in \pi_1(X, x)$. Let $1 \neq \sigma \in S_2$, then $\sigma^*(\gamma^1 \delta^2) = \delta^1\gamma^2$. The $S_2$-invariant functions are then generated by 
  \[\gamma^1 \delta^2 + \delta^1 \gamma^2 \tab\tab \delta, \gamma \in \pi_1(X, x)\]
  \tab Next, notice that the trace functionals
  \[\Rm{tr}_{\gamma}: R_B(X, x, 2) \to \CC, \tab \rho \mapsto \Rm{tr}(\rho(\gamma))\]
  are $\Rm{GL}_2$-invariant and satisfies $\iota^{\#}(\Rm{tr}_{\gamma}) = \gamma^1 + \gamma^2$, hence 
  \[\iota^{\#}(\Rm{tr}_{\gamma}\Rm{tr}_{\delta} - \Rm{tr}_{\gamma\delta}) = \gamma^1\delta^2 + \delta^1\gamma^2\]
  thus $\iota^{\#}$ is a quotient map. As a result, $\iota$ is a closed embedding. 
\end{proof}
Let $q: A \times A \to (A \times A)/S_2$ be the quotient. Since $H^i(X, \chi_1 \oplus \chi_2) = H^i(X, \chi_1) \oplus H^i(X, \chi_2)$, we get that 
\begin{equation}
  \ol{\Sigma}^i_k(X, 2) \cap \Ide{R} = \iota\brac*{\bigcup_{a + b = k} \Sigma^i_a(X, 1) \times \Sigma^i_b(X, 1)} \label{eq:reduc-intersection}
\end{equation}
which is a finite union of images under $q$ of products of torsion-translated subtori, by \cref{lem:rank1-jump}. 
\begin{prop}
  \label{thm:reducible-jump}
  Let $Z$ be an irreducible component of $\ol{\Sigma}^i_k(X, 2)$, the reduced good moduli image of $\Sigma^i_k(X, 2) \subseteq R_B(X, x, 2)$. If the general point of $Z$ is reducible, then there is a product $N = N_1 \times N_2 \subseteq M_B(X, 1) \times M_B(X, 1)$ of torsion-translated subtori such that 
  \[\begin{cases}
    Z^{\nu}\simeq N_1 \times N_2  &\text{ if }N_1 \neq N_2\\
    Z\simeq \Rm{Sym}^2(N_1) &\text{ if }N_1 = N_2
  \end{cases}\]
\end{prop}
\begin{proof}
  Since $\Ide{R}$ is closed, $Z \subset \Ide{R}$. $Z$ is an irreducible component, hence it picks out one component $N = \Sigma^i_a(X, 1) \times \Sigma^i_b(X, 1)$ in \eqref{eq:reduc-intersection}. It follows that $Z = \iota(q(N))$ hence $Z \simeq q(N)$. \\
  \tab If $N_1 \neq N_2$ then $N_1 \cap \sigma N_2$ is closed subset of $N_1$, for $1 \neq \sigma \in S_2$. Therefore the map $N_1 \times N_2 \to q(N_1 \times N_2)$ is finite, surjective with generic degree 1. It follows that it must be the normalization. If $N_1 = N_2$ then $q(N_1 \times N_1) \simeq \Rm{Sym}^2(N_1)$. 
\end{proof}

\subsection{Projectively dihedral jump component}
\label{sec:dihedral-jump}
If, up to conjugation, $G_{\ol{\rho}} \subseteq D_{\infty}$ then $\rho$ fixes a set of two lines $\set*{l_1, l_2}$ with $l_1 \oplus l_2 \simeq \CC^2$. Each $\rho(\gamma)(l_i) \in \set*{l_1, l_2}$ hence we get a permutation character
\[\epsilon: \pi_1(X, x) \to \ZZ/2\ZZ\]
\tab If $\epsilon \neq 0$ (equivalently, $\rho$ is of type $\mathbf{(B)}$), $\rho\vert_{\ker\epsilon}$ then fixes both line hence $\rho\vert_{\ker\epsilon} \simeq \chi_1 \oplus \chi_2$. An element $\delta \in \pi_1(X, x)\setminus \ker\epsilon$ has nontrivial $\epsilon(\delta)$. Conjugation by $\delta$ on $\ker \epsilon \trianglelefteq \pi_1(X, x)$ then permutes $\chi_1$ and $\chi_2$. \\
\\
\tab Now pick any $0 \neq \epsilon \in \Hom(\pi_1(X, x), \ZZ/2\ZZ)$, and denote
\[A^{\epsilon} \coloneq \Hom(\ker\epsilon, \GG_m)\]
then conjugation by any $\delta \in \pi_1(X, x) \setminus \ker\epsilon$ on $\ker \epsilon \triangleleft \pi_1(X, x)$, denoted $\tau: \ker\epsilon \to \ker\epsilon$, induces an involution $\tau^*$ on $A^{\epsilon}$. This is because $\delta^2 \in \ker \epsilon$; the choice of $\delta$ does not matter. 
\begin{lemma}
  \label{lem:dihedral-locus}
  Induction gives a morphism 
  \[\iota_{\epsilon}: A^{\epsilon}/\basis{\tau^*} \to M_B(X, 2), \tab [\chi] \mapsto [\Rm{Ind}^{\pi_1(X, x)}_{\ker \epsilon} \chi]\]
  which is a closed embedding with image $\Ide{D}_{\epsilon}$. The irreducible points in the image correspond to representations of type $\mathbf{(B)}$ with the permutation character $\epsilon$, and form a Zariski open subset of $\sh{D}_{\epsilon}$. 
\end{lemma}
\begin{proof}
  The corresponding map on rings is 
  \[\iota_{\epsilon}^{\#}: H^0(R_B(X, x, 2), \sh{O}_{R_B})^{\Rm{GL}_2} \to H^0(A^{\epsilon}, \sh{O}_{A^{\epsilon}})^{\tau^*}\]
  mapping $f \mapsto \sqbrac*{\chi \mapsto f\brac*{\Rm{Ind}^{\pi_1(X, x)}_{\ker \epsilon} \chi}}$. Explicitly, with the same $\delta$ as before, and $\rho = \Rm{Ind}^{\pi_1(X, x)}_{\ker \epsilon} \chi$,  
  \[\begin{array}{ll}
    \rho(\gamma) = \begin{pmatrix}
    \chi(\gamma) & 0\\
    0 & \chi(\delta^{-1}\gamma \delta)
    \end{pmatrix}, & \gamma \in \ker(\epsilon)\\
    \rho(\delta') = \begin{pmatrix}
      0 & \chi(\delta h\delta)\\
      \chi(h) & 0
    \end{pmatrix}, & \delta' \in \pi_1(X, x)\setminus \ker \epsilon, \delta' = \delta \cdot h
  \end{array}\]
  \tab Now, regular functions on $A^{\epsilon}$ are generated by those of the form $\gamma: \chi \mapsto \chi(\gamma)$ with $\gamma \in \ker(\epsilon)$. Thus $\tau^*-$invariant functions are generated by 
  \[\gamma^{\tau}: \chi \mapsto \chi(\gamma) + \chi(\delta^{-1}\gamma \delta), \tab \gamma \in \ker(\epsilon)\]
  which are clearly $\iota_{\epsilon}^{\#}(\Rm{tr}_{\gamma})$ hence $\iota_{\epsilon}^{\#}$ is a quotient map. 

  \medskip
  Notice that $\iota_{\epsilon}(\chi)\vert_{\ker \epsilon} = \chi \oplus \tau^*\chi$. If $\chi \neq \tau^*\chi$ then the only possible invariant line must be generated by either basis vectors, i.e., $\CC \binom{1}{0}$ or $\CC \binom{0}{1}$. However, any element in $\pi_1(X, x) \setminus \ker\epsilon$ swaps these two lines, thus there is no invariant line and $\iota_{\epsilon}$ is of type $\mathbf{(B)}$. By the discussion at the beginning of the section, this gives all type $\mathbf{(B)}$ respresentations with permutation character $\epsilon$.\\
  \tab Now if $\chi = \tau^*\chi$ then for $\gamma \in \ker\epsilon$ and $\delta' = \delta \cdot h \in \pi_1(X, x) \setminus \ker \epsilon$,
  \[\iota_{\epsilon}(\chi)(\gamma) = \chi(\gamma) \cdot \Rm{Id},\tab \iota_{\epsilon}(\delta') = \chi(h) \begin{pmatrix}
    0 & \chi(\delta^2) \\
    1 & 0
  \end{pmatrix}\]
  then the eigendecomposition for $\iota_{\epsilon}(\delta')$, which is independent of $h$ hence independent of $\delta' \in \pi_1(X, x) \setminus \ker\epsilon$, gives a decomposition of $\iota_{\epsilon}(\chi)$. Finally, the set of $\chi$ fixed by $\tau^*(\chi)$ is closed, hence we get the Zariski open statement. 
\end{proof}
For $0 \neq \epsilon \in \Hom(\pi_1(X, x), \ZZ/2\ZZ)$, let $\pi^{\epsilon}: X^{\epsilon} \to X$ be the connected \'{e}tale double cover coming from $\ker \epsilon\ \triangleleft \pi_1(X, x)$. For $\chi \in A^{\epsilon}$ with the corresponding local system $\bb{L}$ on $X^{\epsilon}$, $\Rm{Ind}^{\pi_1(X, x)}_{\ker \epsilon} \chi$ corresponds to the local system $\pi^{\epsilon}_* \bb{L}$ on $X$. By a straightforward application of Leray's spectral sequence, 
\begin{equation}
  H^i(X, \Rm{Ind}^{\pi_1(X, x)}_{\ker \epsilon} \chi) \simeq H^i(X^{\epsilon}, \chi) \label{eq:leray-induction}
\end{equation}
and we get the following consequence for a type $\mathbf{(B)}$ jump component. 
\begin{prop}
  \label{thm:dihedral-jump}
  Let $Z$ be an irreducible component of $\ol{\Sigma}^i_k(X, 2)$, the reduced good moduli image of $\Sigma^i_k(X, 2) \subseteq R_B(X, x, 2)$. If the general point of $Z$ is of type $\mathbf{(B)}$, then there exists a torsion-translated subtorus $N \subseteq M_B(X^{\epsilon}, 1)$ such that either $Z^{\nu} \simeq N$ or, if $\tau^*N = N$, $Z \simeq N/\basis{\tau^*}$. 
\end{prop}
\begin{proof}
  Since $Z$ is irreducible, and its general point is of type $\mathbf{(B)}$, it must lie in some $\Ide{D}_{\epsilon}$ in \cref{lem:dihedral-locus}. By \eqref{eq:leray-induction}, it must be the image of an irreducible component $N$ of $\Sigma^i_k(X^{\epsilon}, 1)$, hence $N$ is a torsion-translated subtorus by \cref{lem:rank1-jump}. If $N$ is $\tau^*$-invariant, then since $\iota_{\epsilon}$ is a closed embedding, $Z \simeq N/\basis{\tau^*}$. If $N$ is not $\tau^*$-invariant, then $N \cap \tau^* N$ is closed in $N$, thus the composition $N \to A^{\epsilon}/\basis{\tau^*} \xrightarrow{\iota_{\epsilon}} Z$ is a finite map of generic degree 1, hence $N \simeq Z^{\nu}$.  
\end{proof}
\subsection{Projectively finite and rigid jump components}
\label{sec:finite-rigid}
Let $\bb{V}_{\rho} \in M_B(X, 2)$ be an irreducible semisimple representation with corresponding representation $\rho: \pi_1(X, x) \to \Rm{GL}_2(\CC)$. Let $A = M_B(X, 1)$.
\begin{lemma}
  \label{lem:twist-map}
  The morphism 
  \[\theta_{\rho}: A \to M_B(X, 2), \tab \chi \mapsto [\rho \otimes \chi]\]
  is finite, with fibers being the orbits of the finite group 
  \[K(\rho) = \set*{\chi \in A\ \big\vert\ \rho \otimes \chi \simeq \rho} \subseteq A[2]\]
\end{lemma}
\begin{proof}
  Consider the diagram 
  \[\begin{tikzcd}
	A &&&& \\
	&& {M_B(X, 2) \times_{A} A} && {M_B(X, 2)} \\
	\\
	&& A && A
	\arrow["{\chi \mapsto (\rho \otimes \chi, \chi)}"', from=1-1, to=2-3]
	\arrow["{\chi \mapsto \rho \otimes \chi}", curve={height=-18pt}, from=1-1, to=2-5]
	\arrow["\Phi"', from=2-3, to=2-5]
	\arrow[from=2-3, to=4-3]
	\arrow["\det",from=2-5, to=4-5]
	\arrow["{\chi \mapsto \det(\rho) \otimes \chi^{\otimes 2}}", from=4-3, to=4-5]
\end{tikzcd}
  \]
  then $A \to M_B(X, 2) \times_{A} A$ is a closed embedding, and $\Phi$ is finite since the horizontal $A \to A$ is finite. Thus $\theta_{\rho}$, being the composition, is finite. For the fiber identification, notice that $\rho \otimes \chi \simeq \rho \chi'$ if and only if $\chi' \otimes \chi^{-1} \in K(\rho)$. 
\end{proof}
Let $\Gamma \in \set*{A_4, S_4, A_5}$. Since $\pi_1(X, x)$ is finitely generated and $\Gamma$ is finite, there are only finitely many homomorphism $\varphi: \pi_1(X, x) \twoheadrightarrow \Gamma \subset \Rm{PGL}_2(\CC)$. Since every subgroup of $\Rm{PGL}_2$ isomorphic to $\Gamma$ is conjugate to the standard embedding, and we are working with $M_B$, there is no ambiguity with $\varphi$. Suppose $\varphi$ has a $\Rm{GL}_2$-lift $\bb{V}_{\Gamma, \varphi}$. Since any two lifts differ by a rank one twist, we get that 
\begin{equation}
  \Ide{F} = \bigcup_{\Gamma, \varphi} \theta_{\bb{V}_{\Gamma, \varphi}}(A) \subset M_B(X, 2) \label{eq:finite-locus}
\end{equation}
with $\Ide{F}$ the locus of type $\mathbf{(C)}$ representations. 
\begin{prop}
  \label{thm:finite-jump}
  Let $Z$ be an irreducible component of $\ol{\Sigma}^i_k(X, 2)$, the reduced good moduli image of $\Sigma^i_k(X, 2)$. If the general point of $Z$ is of type $\mathbf{(C)}$, then there is a torsion-translated $N \subset A$ such that $Z$ has a normalization
  \[Z^{\nu} = N/H\]
  where $H = \set*{\chi \in K(\rho)\ \vert\ \chi N = N}$. Consequently, $Z^{\nu}$ is a torsion-translated subtorus of $A/H$.
\end{prop}
\begin{proof}
  Since $Z$ is irreducible, it must lie in some $\theta_{\bb{V}_{\Gamma, \varphi}}(A)$, with $\bb{V}_{\Gamma, \varphi}$ a $\Rm{GL}_2$-lift of some $\varphi: \pi_1(X, x) \to \Rm{PGL}_2(\CC)$ with image $\Gamma$. We first claim that it is possible to pick $\bb{V}_{\Gamma, \varphi}$ with finite monodromy. Equivalently, there is some $\chi \in A$ such that $\chi \otimes \bb{V}_{\Gamma, \varphi}$ has finite monodromy. Recall that every class in $A/A^2$ has a torsion lift in $A$, and pick one such lift $\epsilon$ for $[\det{\bb{V}_{\Gamma, \varphi}}] \in A/A^2$. We get that for some $\bb{L} \in A$, $\epsilon = \det(\bb{V}_{\Gamma, \varphi}) \otimes \bb{L}^{\otimes 2}$. Then 
  \[\det(\bb{V}_{\Gamma, \varphi} \otimes \bb{L}) = e\]
  is torsion. It follows that $\bb{V}_{\Gamma, \varphi} \otimes \bb{L}$ has finite monodromy projective quotient, and finite monodromy determinant. Since the map 
  \[\Rm{GL}_2 \to \Rm{PGL}_2 \times \GG_m, \tab g \mapsto ([g], \det g)\]
  is finite of degree 2, we get that $\bb{V}_{\Gamma, \varphi} \otimes \bb{L}$ has finite monodromy. In particular, $\bb{V}_{\Gamma, \varphi} \otimes \bb{L}$ is of geometric origin.\\
  \\
  \tab We can now assume that $\bb{V}_{\Gamma, \varphi}$ is of geometric origin. Pick an irreducible component of $\theta_{\bb{V}_{\Gamma, \varphi}}^{-1}(Z)$. It must be an irreducible component of the twisted jump loci 
  \[\Sigma^i_k(X, \bb{V}_{\Gamma, \varphi}, 1) = \set*{\chi \in A\ \big\vert\ \dim H^i(X, \bb{V}_{\Gamma, \varphi} \otimes \chi) \geq k}\]
  which is then a torsion-translated subtorus $N \subset A$ by \cref{lem:rank1-jump}. It follows that $Z = \theta_{\bb{V}_{\Gamma, \varphi}}(N)$. Next, denote $H = \set*{\chi \in K(\rho)\ \vert\ \chi N = N}$, then since $H$ acts by translation, $N/H$ is a translated subtorus of $A/H$. Finally, $N/H \to \theta_{\bb{V}_{\Gamma, \varphi}}(N)$ is finite of degree 1 (since the extra identifications come from the fixed loci of $\chi \in K(\rho)$ with $\chi N \neq N$; these fixed loci are closed proper subset of $N$). As a result, $N/H \to Z$ is the normalization map. 
\end{proof}

\begin{prop}
  \label{thm:rigid-jump}
  Let $Z$ be an irreducible component of $\ol{\Sigma}^i_k(X, 2)$, the reduced good moduli image of $\Sigma^i_k(X, 2)$. If the general point of $Z$ is of type $\mathbf{(D)}$. Then there exists a translated subtorus $N \subseteq A$ such that $Z^{\nu} \simeq N/H$ with $\set*{\chi \in K(\rho)\ \vert\ \chi N = N}$.  
\end{prop}
\begin{proof}
  Here we must have $Z \subseteq \theta_{\bb{V}_{\Gamma, \varphi}}(A)$, with $\bb{V}_{\Gamma, \varphi}$ a $\Rm{GL}_2$-lift of an isolated point in $M_B(X, \Rm{PGL}_2)$. Hence it is the image under $\theta_{\bb{V}_{\Gamma, \varphi}}(A)$ of an irreducible component of 
  \[\Sigma^i_k(X, \bb{V}_{\Gamma, \varphi}, 1) = \set*{\chi \in A\ \big\vert\ \dim H^i(X, \bb{V}_{\Gamma, \varphi} \otimes \chi) \geq k}\]
  hence it is $\theta_{\bb{V}_{\Gamma, \varphi}}(N)$ with $N \subseteq A$ a translated subtorus by \cref{lem:rank1-jump} (we lose that the translation is torsion). The remaining statements are exactly as in the proof of \cref{thm:finite-jump}.  
\end{proof}

\section{Nonabelian jump component}
This section will handle the remaining case from the previous section, i.e., $Z$ is an irreducible component of $\ol{\Sigma}^i_k(X, 2)$ with general point of type $\Bf{(E)}$ (see \cref{lem:rep-type}). Let $A = M_B(X, 1)$. Following \cref{sec:jump-product}, $Z$ is contained in some $M_B(X, 2, \det \in \epsilon A^2)$ with torsion $\epsilon \in A$. Pick an irreducible component
\begin{equation}
  \begin{tikzcd}
	{} & {C \subseteq \Phi_{\epsilon}^{-1}(Z) \subseteq M_B(X, 2, \det = \epsilon) \times A} & \\
	& {M_B(X, 2, \det = \epsilon)} & A
	\arrow["{\Rm{pr}_1}"', from=1-2, to=2-2]
	\arrow["{\Rm{pr}_2}", from=1-2, to=2-3]
  \end{tikzcd}\label{eq:det-sheet}
\end{equation}
then by \cref{lem:normalized-det-sheet}, the normalized $C^{\nu} \to Z^{\nu}$ is a finite \'{e}tale torsor. 

First we take a detour to talk about local systems and cohomology on a Deligne-Mumford quotient stack. 
\subsection{Local systems on finite quotient curves}
We are following the exposition in \cite[\textsection 2.1 and 2.2]{corlette2008classification}. For more details, see \cite{noohi2004fundamental}. 
\label{sec:finite-quot}
\begin{definition}
  \label{def:DM-curve}
  A \textit{DM curve} is a connected smooth one-dimensional complex Deligne-Mumford stack, and an \textit{orbicurve} is a DM curve with trivial generic stabilizer. 
\end{definition}
In this paper, we work only with DM curves $\sh{Y}$ of the form $\sh{Y} = [Y/G]$ where $Y$ is a connected smooth projective curve and $G$ is a finite group acting on $Y$. 
\begin{definition}
  \label{def:fund-grp-DM}
  Let $\bar{y} \in \sh{Y}$ be a point, then we have the fundamental group $\pi_1(\sh{Y}, \bar{y})$. For $y \in Y$ lying over $\bar{y} \in \sh{Y}$ in the \'{e}tale map 
  \[\pi: Y \to \sh{Y} = [Y/G]\]
  an element of $\pi_1(\sh{Y}, \bar{y})$ is a pair $(g, \gamma)$ with $g \in G$ and $\gamma$ a path from $y$ to $g \cdot y$. Thus we have an exact sequence 
  \[1 \to \pi_1(Y, y) \to \pi_1(\sh{Y}, \bar{y}) \to G \to 1\]
\end{definition}
\begin{definition}
  \label{def:local-sys-DM}
  A complex local system on $\sh{Y}$ is a $G$-equivariant local system on $Y$, i.e. a pair $(\bb{V}, \alpha)$ where $\bb{V}$ is a local system on $Y$ and $\alpha$ is an action associating to each $g \in G$ an isomorphism $\alpha_g: g^*\bb{V} \xrightarrow{\simeq} \bb{V}$ such that 
  \[\alpha_g \circ g^*\alpha_h = \alpha_{gh}, \tab g, h \in G\]
  \tab Then we recover the familiar equivalence \cite[Theorem 18.19]{noohi2005foundationsi}
  \[\Rm{Loc}^{G}(Y) \simeq \Rm{Loc}_{\CC}(\sh{Y}) \simeq \Rm{Rep}_{\CC}(\pi_1(\sh{Y}, \bar{y}))\]
  between complex local systems on $\sh{Y}$ and complex representations of $\pi_1(\sh{Y}, \bar{y})$. \\
  \tab We have similar definitions for constructible complexes and perverse sheaves, i.e. $D^b_c(\sh{Y}, \CC) = D^{b, G}_c(Y, \CC)$ and $\Rm{Perv}(\sh{Y}) = \Rm{Perv}^{G}(Y)$. For $K \in D^{b}_c(\sh{Y}, \CC)$, $\Rm{R}\Gamma(\sh{Y}, K) \simeq \Rm{R}\Gamma(Y, \pi^*K)^{G}$. 
\end{definition}
\begin{definition}
  \label{def:char-var-DM}
  Since $\Rm{Loc}_{\CC}(\sh{Y}) \simeq \Rm{Rep}_{\CC}(\pi_1(\sh{Y}, \bar{y}))$, we have framed character varieties of $\sh{Y}$,
  \[R_B(\sh{Y}, \bar{y}, r) \coloneq \Hom(\pi_1(\sh{Y}, \bar{y}), \Rm{GL}_r)\]
  which are affine schemes. Its affine good moduli space is the GIT quotient 
  \[M_B(\sh{Y}, r) \coloneq R_B(\sh{Y}, \bar{y}, r)//\Rm{GL}_r = \Spec H^0(R_B(\sh{Y}, \bar{y}, r), \sh{O}_{R_B})^{\Rm{GL}_r}\]
  by the conjugation action.
\end{definition}
For $y \in Y$ let $G_y = \set*{g \in G\ \vert\ g\cdot y = y}$ be the stabilizer group. Notice that replacing $y$ by $g\cdot y$ identifies $G_y$ with $G_{gy}$ by conjugation, which then gives a canonical isomorphism 
\[M_B(BG_y, r) \simeq M_B(BG_{gy}, r)\]
where $BH = [\Spec \CC/H]$. Denote $I_{\eta} = \ker(G \to \Rm{Aut}(Y))$ the generic inertia, then for general $y \in Y$, $G_y = I_{\eta}$. There are finitely many exceptional points $y \in Y$ with $I_{\eta} \subsetneq G_y$; pick their $G-$orbit representatives $y_1, \dots, y_l$ with $G_j = G_{y_j}$. Pulling back along the residual gerbes $\iota_j: BG_{j} \to \sh{Y}$ gives 
\[\iota_j^*: M_B(\sh{Y}, r) \to M_B(BG_j, r), \tab 0 \leq j \leq l\]
where $\iota_0: BI_{\eta} \to \sh{Y}$ be the residual gerbe of a generic $y \in Y$. Since $I_{\eta}$ acts trivially on $Y$, $\iota_0^*$ does not depend on the choice of general point. 
\begin{definition}
  \label{def:fixed-inertia}
  For $\bm{\tau} = (\tau_j) \in \Prod{1 \leq j \leq l} M_B(BG_j, r)$, we define the \textit{fixed-inertia locus} to be 
  \[M_B(\sh{Y}, r, \bm{\tau}) = \bigcap_{0 \leq j \leq j} \iota_j^*(\tau_j)\]
  \tab Since $M_B(BG_j, r)$ is finite and reduced for $G_j$ finite, we get that every irreducible subvariety of $M_B(\sh{Y}, r)$ lies in one such fixed-inertia locus. 
\end{definition}

\begin{prop}
  \label{thm:DM-pi1-misc}
  For $\pi: Y \to \sh{Y} = [Y/G]$ we have the following 
  \begin{enumerate}
    \item If $j: \sh{U} \to \sh{Y}$ is a nonempty open substack, then $j_*:\pi_1(\sh{U}, \bar{u}) \to \pi_1(\sh{Y}, \bar{u})$ is surjective. Consequently restriction of local systems from $\sh{Y}$ to $\sh{U}$ is fully faithful and preserves irreducibility.  
    \item If $f: \sh{X} \to \sh{Y}$ induces surjection on fundamental groups, then 
    \[f^*: M_B(\sh{Y}, r) \to M_B(\sh{X}, r)\] is a closed immersion. 
  \end{enumerate}
\end{prop}
\begin{proof}
  Let $U = \sh{U} \times_{\sh{Y}} Y$ then $U$ is a $G-$invariant subset of $Y$, and we have a commutative diagram 
  \[
  \begin{tikzcd}
	1 & {\pi_1(U, u)} & {\pi_1(\sh{U}, \bar{u})} & G & 1 \\
	1 & {\pi_1(Y, u)} & {\pi_1(\sh{Y}, \bar{u})} & G & 1
	\arrow[from=1-1, to=1-2]
	\arrow[from=1-2, to=1-3]
	\arrow[two heads, from=1-2, to=2-2]
	\arrow[from=1-3, to=1-4]
	\arrow[from=1-3, to=2-3]
	\arrow[from=1-4, to=1-5]
	\arrow[no head, from=1-4, to=2-4]
	\arrow[shift left, no head, from=1-4, to=2-4]
	\arrow[from=2-1, to=2-2]
	\arrow[from=2-2, to=2-3]
	\arrow[from=2-3, to=2-4]
	\arrow[from=2-4, to=2-5]
\end{tikzcd}\]
  where the left vertical surjection comes from $Y \setminus U$ being finitely many points. Thus we get a surjection $\pi_1(\sh{U}, \bar{u}) \twoheadrightarrow \pi_1(\sh{Y}, \bar{u})$. The second statement (2) is exactly the same as for character varieties of varieties. 
\end{proof}

\subsection{A factorization theorem following Corlette-Simpson}
\label{sec:C-S-factorize}
With the notations in \eqref{eq:det-sheet}, let 
\[D = \ol{\Rm{pr}_1(C)}_{\red} \subset M_B(X, 2, \det = \epsilon), \tab T = \ol{\Rm{pr}_2(C)}_{\red} \subset A\]
with $T$ a torsion-translated subtorus, and $\epsilon \in A$ torsion. Notice that since a general point of $Z$ has dense $\Rm{PGL}_2$-image, the same is true for $D$. 
\begin{remark}
  \label{rmk:proj-nonrigidity} 
  We can look at the image of $D$ under $\PP: M_B(X, 2) \to M_B(X, \Rm{PGL}_2)$. If the image is zero-dimensional, then $D$ itself must be finite since two fixed-determinant lifts of a $\Rm{PGL}_2$-representation differs by an element of $A[2]$ which is finite. Since $C$ is irreducible, we must have $D = \set*{\bb{V}}$ and $C = \set*{\bb{V}} \times T$. From now on, we assume that $\PP(D)$ is positive-dimensional in $M_B(X, \Rm{PGL}_2)$. 
\end{remark}

The goal of this section is the following proposition:
\begin{prop}
  \label{thm:Shafarevich-map}
  There is a proper hyperbolic DM curve $\sh{Y} = [Y/G]$ and a morphism $s: X \to \sh{Y}$ and an irreducible component $Q$ of a fixed-determinant, fixed-inertia locus $M_B(\sh{Y}, 2, \det = \delta, \bm{\tau})$ such that
  \begin{enumerate}
    \item $s_*: \pi_1(X, x) \to \pi_1(\sh{Y}, \bar{y})$ is a surjection.
    \item with the closed embedding $s^*: M_B(\sh{Y}, 2) \hookrightarrow M_B(X, 2)$, $D \subseteq s^*(Q)$. 
  \end{enumerate}
\end{prop}
Let $p: X^{\epsilon} \to X$ be the connected finite \'{e}tale cover corresponding to $\ker \epsilon \triangleleft \pi_1(X, x)$. Then for every $\rho \in D \subset M_B(X, 2, \det = \epsilon)$, the restriction $\rho\vert_{X^{\epsilon}} \in M_B(X^{\epsilon}, \Rm{SL}_2)$. Call this map $\Rm{res}_{\epsilon}: M_B(X, 2, \det = \epsilon) \to M_B(X^{\epsilon}, \Rm{SL}_2)$. 

\begin{lemma}
  \label{lem:SL2-dense-nonrigid}
  The general point of $\Rm{res}_{\epsilon}(D)$ has Zariski-dense image in $\Rm{SL}_2$, and is non-rigid. 
\end{lemma}
\begin{proof}
  We start with some claims.
  \begin{claim}
    If $\bar{\rho}: \pi_1(X, x) \to \Rm{PGL}_2(\CC)$ has Zariski-dense image in $\Rm{PGL}_2$ and $\Gamma \triangleleft \pi_1(X, x)$ a normal finite-index subgroup, then $\bar{\rho}(\Gamma)$ is Zariski-dense in $\Rm{PGL}_2$. 
  \end{claim}
  \begin{claimproof} 
    Let $H = \ol{\bar{\rho}(\Gamma)}$ be the Zariski closure. Then $H$ is a normal finite-index subgroup of $\Rm{PGL}_2$. However, $\Rm{PGL}_2$ is connected, hence $H = \Rm{PGL}_2$.
  \end{claimproof}
   
  \begin{claim}
    On the locus of Zariski-dense representations, the restriction 
    \[\Rm{res}: M_B(X, \Rm{PGL}_2) \to M_B(X^{\epsilon}, \Rm{PGL}_2)\]
    is injective on closed points. Consequently, $\Rm{res}_{\epsilon}(D)$ has positive-dimensional image in $M_B(X^{\epsilon}, \Rm{SL}_2)$. 
  \end{claim}
  \begin{claimproof}
    Suppose otherwise, and $\bar{\rho}_1, \bar{\rho}_2$ are Zariski-dense projective representations such that their restrictions to $X^{\epsilon}$ are conjugate. Since $\bar{\rho}_2 \in M_B(X, \Rm{PGL}_2)$ only depends on conjugacy class, we might assume that $\bar{\rho}_1\vert_{X^{\epsilon}} = \bar{\rho}_2\vert_{X^{\epsilon}} = \bar{\rho}_{\epsilon}$. For $\gamma \in \pi_1(X, x)$ and $h \in \ker \epsilon$, denote 
  \[c_{\gamma} = \bar{\rho}_2(\gamma)^{-1} \circ \bar{\rho}_1(\gamma) \in \Rm{PGL}_2\]
  and we have, since $\gamma h \gamma^{-1} \in \ker(\epsilon)$, 
  \[c_{\gamma} \circ \bar{\rho}_{\epsilon}(h) \circ c_{\gamma}^{-1} = \bar{\rho}_{\epsilon}(h)\]
  \tab Now, by the previous claim, $\bar{\rho}_{\epsilon}$ has Zariski-dense image. Hence the centralizer of $\bar{\rho}_{\epsilon}(\ker \epsilon)$ agrees with the center $Z(\Rm{PGL}_2)$ which is trivial. Thus $c_{\gamma} = \Rm{Id}$ for all $\gamma \in \pi_1(X, x)$. It follows that
  \[\Rm{res}: M_B(X, \Rm{PGL}_2) \to M_B(X^{\epsilon}, \Rm{PGL}_2)\]
  is injective on closed points. We have a diagram 
  \[
  \begin{tikzcd}
	{M_B(X, 2, \det = \epsilon)} && {M_B(X^{\epsilon}, \Rm{SL}_2)} \\
	{M_B(X, \Rm{PGL}_2)} && {M_B(X^{\epsilon}, \Rm{PGL}_2)}
	\arrow["{\Rm{res}_{\epsilon}}", from=1-1, to=1-3]
	\arrow["\PP"', from=1-1, to=2-1]
	\arrow["\PP", from=1-3, to=2-3]
	\arrow["{\Rm{res}}", from=2-1, to=2-3]
\end{tikzcd}\] 
  where both vertical maps are finite, $\PP(D)$ is positive dimensional and its general point has Zariski-dense image, hence by the second claim $\Rm{res} \circ \PP(D)$ has positive dimensional image. It follows that $\Rm{res}_{\epsilon}(D)$ has positive dimensional image. 
  \end{claimproof}
  Thus the general point $\rho_{\epsilon} \in \Rm{res}_{\epsilon}(D)$ is non-rigid. By the diagram, such a general point has $\PP(\rho_{\epsilon})$ with Zariski-dense image in $\Rm{PGL}_2$. On the other hand, let $\ol{\rho_{\epsilon}(\ker \epsilon)} \subseteq \Rm{SL}_2$ be the Zariski closure of the image of $\rho_{\epsilon}$, then $\ol{\rho_{\epsilon}(\ker \epsilon)} \twoheadrightarrow \Rm{PGL}_2$. As a result, $\dim \ol{\rho_{\epsilon}(\ker \epsilon)} \geq 3$ hence $\ol{\rho_{\epsilon}(\ker \epsilon)}  = \Rm{SL}_2$. 
\end{proof}
For a general $\rho \in D$, by \cref{lem:SL2-dense-nonrigid}, the restriction $\rho\vert_{X^{\epsilon}}$ has Zariski-dense image and is non-rigid in $M_B(X^{\epsilon}, \Rm{SL}_2)$. Hence we have the following as immediate consequences of \cite{corlette2008classification}; see also \cite[Theorem 10]{simpson1991ubiquity}. 
\begin{lemma}[{\cite[Corollary 3.3 and Theorem 6.8]{corlette2008classification}}]
  \label{lem:proj-factorize1}
  For a general $\rho \in D$, the restricted projective representation $\bar{\rho}\vert_{X^{\epsilon}}$ factors through an algebraic map $X^{\epsilon} \to Y$ where $Y$ is a proper hyperbolic orbicurve.
\end{lemma}
To get a factorization for $\bar{\rho}$ instead of $\bar{\rho}\vert_{X^{\epsilon}}$ we have the following lemma.
\begin{lemma}[{\cite[Lemma 3.6]{corlette2008classification}$+\epsilon$}]
  \label{lem:proj-factorize2}
  Let $f: X' \to X$ be a finite \'{e}tale map of smooth projective varieties, and $\bar{\rho}: \pi_1(X, x) \to \Rm{PGL}_2(\CC)$ is a Zariski-dense representation such that $\bar{\rho}\vert_{X'}$ factors through a proper hyperbolic orbicurve. Then $\bar{\rho}$ itself factors through a proper hyperbolic orbicurve. 
\end{lemma}
\begin{proof}
  By taking the Galois closure, we might assume $f:X' \to X$ is a Galois cover and $X = [X'/G]$. Let $s: X' \to Y$ be the factorization map for $\bar{\rho}\vert_{X'}$ with $Y$ proper hyperbolic orbicurve. Let $Y' \to Y$ be the finite cover by a smooth curve. Replacing $X'$ with $X' \times_{Y} Y'$ we can assume that the factorization map $s': X' \to Y$ has smooth curve target. Take the Stein factorization, and remove from $X$ the closure of the image of all fibers of $s'$ which are reducible or not smooth (passing to an open subset of $X$ is fine by \cite[Lemma 3.5]{corlette2008classification}). Thus we may assume that all fibers of $s'$ are irreducible, and smooth.\\
  \tab The fibers of $s'$ are maximal subvarieties $W$ of $X'$ such that $\bar{\rho}\vert_{W}$ is trivial. For $g \in G = \Rm{Gal}(X'/X)$, $f \circ g = f$ hence $g$ permutes the fibers of $s': X' \to Y$. Hence the $G$-action descends to $Y$ making $s'$ a $G-$equivariant map. Thus it descends to $s: X \to [Y/G]$. \\
  \tab $[Y/G]$ is an orbicurve if and only if the $G$-action is faithful. Let $I_{\eta} = \ker(G \to \Rm{Aut}(Y))$, then $I_{\eta}$ is finite, and a normal subgroup of $\pi_1([Y/G], \bar{y})$. By the previous paragraphs, $\bar{\rho}$ descends to some $\bar{\theta}: \pi_1([Y/G], \bar{y}) \to \Rm{PGL}_2(\CC)$ which must also have Zariski-dense image. Then the normalizer $N_{\Rm{PGL}_2}(\bar{\theta}(I_{\eta}))$ contains $\bar{\theta}(\pi_1([Y/G], \bar{y}))$, thus must be $\Rm{PGL}_2$. Hence $\bar{\theta}(I_{\eta})$ is finite and normal in $\Rm{PGL}_2$, which gives that $\bar{\theta}(I_{\eta})$ is trivial. We can then replace $[Y/G]$ by $[Y/(G/K)]$ which has trivial generic inertia. 
\end{proof}
\begin{coroll}
  \label{thm:proj-factorize3}
  For a general $\rho \in D$, the projective representation $\bar{\rho}$ factors through a nonconstant map $s_{\bar{\rho}}: X \to Y_{\bar{\rho}}$ where $Y_{\bar{\rho}}$ is a proper hyperbolic orbicurve. 
\end{coroll}
\begin{proof}
  Apply \cref{lem:proj-factorize2} to the cover $X^{\epsilon} \to X$, and use \cref{lem:proj-factorize1} to get factorization on $X^{\epsilon}$. 
\end{proof}
Notice that if we have a finite cover $\pi: Y \to Y_{\bar{\rho}}$, and a lift $s: X \to Y$ of $s_{\bar{\rho}}$, then $s^*_{\bar{\rho}} = s^* \circ \pi^*$. Hence $s: X \to Y$ also satisfies the factorization property. By the following lemma, we can assume that $(s_{\bar{\rho}})_*:\pi_1(X, x) \to \pi_1(Y_{\bar{\rho}}, y)$ is a surjection. 

\begin{lemma}
  \label{lem:pi1-surjection}
  Let $s: X \to Y$ be a nonconstant algebraic morphism from projective connected variety $X$ to a proper hyperbolic orbicurve $Y$. Then there is a finite \'{e}tale cover $\pi: Y' \to Y$ and an algebraic lift $s': X \to Y'$ such that $s'$ is surjective and $s_{*}': \pi_1(X, x) \to \pi_1(Y', y')$ is surjective. 
\end{lemma}
\begin{proof}
  By \cite[Proposition 18.18 and Corollary 18.20]{noohi2005foundationsi} there is a cover $\pi: Y' \to Y$ corresponding to $\im(s_*) \subseteq \pi_1(Y, y)$, and an analytic lift $s': X \to Y'$ with $s_*'$ surjective. Since $X$ is a projective variety and $s'$ is nonconstant, we get that $\im(s') = Y'$, and $Y'$ is also proper. Thus $Y' \to Y$ is a finite analytic covering map. By \cite[Theorem 20.4]{noohi2005foundationsi} this is algebraizable. 
\end{proof}

\begin{coroll}
  \label{thm:proj-factorize4}
  For a general $\rho \in D$, the projective representation $\bar{\rho}$ factors through a map $s_{\bar{\rho}}: X \to Y_{\bar{\rho}}$ where $Y_{\bar{\rho}}$ is a proper hyperbolic orbicurve, $s_{\bar{\rho}}$ is surjective, and the induced map $(s_{\bar{\rho}})_*: \pi_1(X, x) \to \pi_1(Y_{\bar{\rho}}, y)$ is surjective. 
\end{coroll}

\begin{prop}
  \label{thm:proj-common}
  There is a proper hyperbolic orbicurve $Y$ and an algebraic $s: X \to Y$ such that $s_*: \pi_1(X, x) \to \pi_1(Y, y)$ is surjective, and for every $\rho \in D$, the projective representation $\bar{\rho}$ factors through $s$. 
\end{prop}
\begin{proof}
  By \cite[Proposition 2.8]{corlette2008classification}, there are only finitely many maps $s_j: X \to Y_j$ with $Y_j$ a hyperbolic orbicurve, and $s_j$ is surjective. By \cref{lem:pi1-surjection}, we can assume all $(s_j)_*$ are surjective, thus by \cref{thm:DM-pi1-misc}, we have closed embeddings 
  \[s_j^*: M_B(Y_j, \Rm{PGL}_2) \hookrightarrow M_B(X, \Rm{PGL}_2)\]
  \tab Let $\ol{\PP(D)}$ be the closure of the image of $D$ under 
  \[\PP: M_B(X, 2) \to M_B(X, \Rm{PGL}_2)\]
  then by \cref{thm:proj-factorize4}, a dense open subset of $\ol{\PP(D)}$ must be contained in the finite union $\bigcup_j \im(s_j^*)$ which is Zariski-closed. Hence $\ol{\PP(D)} \subseteq \bigcup_j \im(s_j^*)$. Since $D$ is irreducible, so is $\ol{\PP(D)}$ and thus must lie entirely in some $\im(s_j^*)$, and we are done. 
\end{proof}
Next we recover a factorization for every $\rho \in D$ through the projective factorization. See also \cite[Lemma 3.1]{corlette2008classification} for the strategy in the case $\epsilon$ is trivial, which we will adapt. 
\begin{lemma}
  \label{lem:shafarevich-map}
  Let $Y$ and $s$ be as in \cref{thm:proj-common}. There is a finite gerbe $\sh{Y} \to Y$ and a lift $s_{D}: X \to \sh{Y}$ such that $s_{D, *}: \pi_1(X, x) \to \pi_1(\sh{Y}, \bar{y})$ is surjective, and every $\rho \in D$ factors through $s_D: X \to \sh{Y}$.
\end{lemma}
\begin{proof}
  For $\rho \in D$, notice that $\bar{\rho}$ factors through $s: X \to Y$ hence $\ker(s_*) \subseteq \ker(\bar{\rho})$. It follows that for every $k \in \ker(s_*)$, $\rho(k)$ is scalar, and we can write 
  \[\rho(k) = \lambda_{\rho}(k)\cdot \Rm{Id}, \tab \lambda_{\rho}(k)^2 = \det \rho(k) = \epsilon(k)\] 
  \tab This $\lambda_{\rho}$ is independent of $\rho \in D$. 
  \begin{claim}
    There is a $\pi_1(X, x)$-invariant character $\lambda: \ker(s_*) \to \mu_{2m}$ such that $\rho(k) = \lambda(k)\cdot\Rm{Id}$ for every $\rho \in D$ and $k \in \ker(s_*)$. 
  \end{claim}
  \begin{claimproof}
    Indeed, for each $k \in \ker(s_*)$ consider the regular function $f_k \in \sh{O}_D(D)$ by 
  \[f_k(\rho) = \frac{1}{2} \Rm{tr}(\rho(k))\]
  then $f_k(\rho) = \lambda_{\rho}(k)$, $f_k^2(\rho) = \epsilon(k)$ for all $\rho \in D$. Since $\sh{O}_D(D)$ is integral, $f_k(\rho)$ must be constant and equal to one of the two square roots of $\epsilon(k)$; denote this square root $\lambda(k) = f_k$. For $k_1, k_2 \in \ker(s_*)$ and $\rho \in D$, 
  \[\lambda(k_1k_2) \cdot \Rm{Id} = \lambda_{\rho}(k) \cdot \Rm{Id} = \rho(k_1k_2) = \rho(k_1)\rho(k_2) = \lambda(k_1)\lambda(k_2)\cdot \Rm{Id}\]
  hence $\lambda$ is multiplicative. Similarly, one can check that $\lambda$ is invariant under conjugation by $\pi_1(X, x)$. Finally, since $\epsilon$ is torsion, there is $m \in \ZZ_+$ such that $\epsilon^m = 1$ hence $\lambda^{2m} = 1$. 
  \end{claimproof}
  We thus have a central (since $\sqbrac*{\pi_1(X, x), \ker(s_*)} \subseteq \ker(\lambda)$) extension 
  \[1 \to \ker(s_*)/\ker(\lambda) \to \pi_1(X, x)/\ker(\lambda) \xrightarrow{s_*} \pi_1(Y, y) \to 1\]
  with $\ker(s_*)/\ker(\lambda) \simeq \im(\lambda)$ finite. By \cite[Lemma 2.5]{corlette2008classification}, this determines a finite $\ker(s_*)/\ker(\lambda)$-gerbe $\sh{Y} \to Y$ with $\pi_1(\sh{Y}, \bar{y}) = \pi_1(X, x)/\ker(\lambda)$. Furthermore, the quotient $\pi_1(X, x) \to \pi_1(\sh{Y}, \bar{y})$ is a lift of $s_*: \pi_1(X, x) \to \pi_1(Y, y)$ which gives a lift $s_D: X \to \sh{Y}$ of $s: X \to Y$. By construction, $s_{D, *}$ is surjective. Notice that every $\rho \in D$ has $\rho(\ker(\lambda)) = 1$ hence $\rho$ factors through $\pi_1(\sh{Y}, \bar{y})$ as desired. 
\end{proof}
\begin{proof}[Proof of \cref{thm:Shafarevich-map}] 
  By \cref{lem:shafarevich-map} we get a proper hyperbolic DM curve $\sh{Y}$ and $s_D: X \to \sh{Y}$ such that every $\rho \in D$ factors through $s_D$, i.e. $D$ lies in the image of the closed embedding
  \[s_D^*: M_B(\sh{Y}, 2) \to M_B(X, 2)\]
  \tab The preimage $(s_D^*)^{-1}(D)$ is irreducible. Furthermore, for every $\psi \in (s_D^*)^{-1}(D)$, 
  \[s_D^*(\det \psi) = \det(s_D^* \psi) = \epsilon\]
  and $s_D^*: M_B(\sh{Y}, 1) \hookrightarrow M_B(X, 1)$ is also injective hence $\det(\psi)$ is constant for all $\psi \in (s_D^*)^{-1}(D)$. It follows that $\psi \in (s_D^*)^{-1}(D)$ has fixed-determinant, and thus must lie in a fixed-determinant, fixed-inertia locus by \cref{def:fixed-inertia}. Finally, the construction of \cref{lem:proj-factorize2} gives an proper hyperbolic orbicurve of the form $[Y/G]$. The finite gerbe $\sh{Y} \to [Y/G]$ from \cref{lem:shafarevich-map} remains a proper hyperbolic finite quotient curve. 
\end{proof}

\subsection{$C \to T$ is either a product or generically finite}
\label{sec:tbd1}
Recall \eqref{eq:det-sheet} that we have $C$ an irreducible component of $\Phi_{\epsilon}^{-1}(Z) \subseteq M_B(X, 2, \det = \epsilon) \times A$ with $Z$ an irreducible component of the reduced jump loci $\ol{\Sigma}^i_k(X, 2) \subseteq M_B(X, 2)$. And the two projections give
\[D = \ol{\Rm{pr}_1(C)}_{\red} \subseteq M_B(X, 2), \tab T = \ol{\Rm{pr}_2(C)}_{\red} \subseteq A\] 
and by \cref{thm:Shafarevich-map}, there is a proper hyperbolic DM curve $\sh{Y} = [Y/G]$ (with $Y$ smooth and $G$ finite), an irreducible component $Q$ of a fixed-determinant, fixed-inertia locus of $M_B(\sh{Y}, 2)$ such that, and a map $s: X \to \sh{Y}$ such that $D \subseteq s^*(Q)$. Since $C \subseteq D \times T$, and 
\[(s^*, \Rm{id}): M_B(\sh{Y}, 2) \times A \to M_B(X, 2) \times A\]
is a closed embedding, there is some $\widetilde{C} \subseteq Q \times T$ isomorphic to $C$. We might as well work with $C \subseteq Q \times T$. Denote 
\[\widetilde{\Sigma}^i_k = \set*{(\bb{W}, \bb{L}) \in Q \times T\ \big\vert\ \dim H^i(s^*\bb{W} \otimes \bb{L}) \geq k}\]
then $C \subseteq \widetilde{\Sigma}^i_k$. We start with a decomposition lemma for the proper map $s: X \to \sh{Y}$. We have the following diagram 
\begin{equation}
  \begin{tikzcd}
	{X'} && X \\
	\\
	Y && {\sh{Y}}
	\arrow["{\pi'}", from=1-1, to=1-3]
	\arrow["{s'}"', from=1-1, to=3-1]
	\arrow["s", from=1-3, to=3-3]
	\arrow["\pi", from=3-1, to=3-3]
  \end{tikzcd}\label{eq:pullback-DM}
\end{equation}
where $\pi$ and $\pi'$ are both finite \'{e}tale $G$-torsor. $X = [X'/G]$, and $s': X' \to Y$ is $G-$equivariant. We have base change for $K \in D^b_c(X, \CC)$,
\[\pi^* \Rm{R}s_* K \simeq \Rm{R}s_*' (\pi')^*K\]
\begin{lemma}
  \label{lem:DM-decomp}
  Let $\bb{L}$ be a rank one local system on $X$, then there is $G$-invariant dense open $j: U \hookrightarrow Y$ and a $G$-equivariant decomposition
  \[\Rm{R}s_*'(\pi')^*\bb{L} \simeq \bigoplus_{\alpha} j_{*} \bb{U}_{\alpha}[-a_{\alpha}] \oplus \bigoplus_{\beta} \iota_{\beta, *} \sh{V}_{\beta}[-b_{\beta}]\]
  where $a_{\alpha}, b_{\beta} \in \ZZ$, $\bb{U}_{\alpha}$ is an irreducible $G$-equivariant local system on $U$, and $\sh{V}_{\beta}$ is an irreducible local system on the orbit $G\cdot y_{\beta}$ of a point $y_{\beta} \in Y \setminus U$ with $G$-equivariant inclusion $\iota_{\beta}: G \cdot y_{\beta} \hookrightarrow Y$. \\
  \tab For $\bb{W}$ a local system on $\sh{Y}$, let $V_{\beta} = (\sh{V}_{\beta})_{y_{\beta}}$ and $G_{\beta}$ the stabilizer of $y_{\beta}$. There is a decomposition of cohomology 
  \begin{align*}
    H^i(X, s^*\bb{W} \otimes \bb{L}) \simeq &\bigoplus_{\alpha} H^{i-a_{\alpha}}(Y, (\pi^*\bb{W}) \otimes j_{*} \bb{U}_{\alpha})^{G}\\
    &\oplus \bigoplus_{\beta, b_{\beta} = i} ((\pi^*\bb{W})_{y_{\beta}} \otimes V_{\beta})^{G_{\beta}}
  \end{align*}
\end{lemma}
\begin{proof}
  We give a sketch here. $\bb{L}' = (\pi')^*\bb{L}$ has rank one, and thus is irreducible. By \cite[Main Theorem 1]{sabbah2005polarizable}, we have a decomposition
  \begin{equation}
    \Rm{R}s_*'(\bb{L}'[-\dim X]) \simeq \bigoplus_{i} {}^p\sh{H}^i(\Rm{R}s_*'(\bb{L}'[-\dim X]))[-i] \label{eq:sabbah-decomp}
  \end{equation} 
  where each ${}^p\sh{H}^i(\Rm{R}s_*'(\bb{L}'[-\dim X]))$ further decomposes into intersection complexes supported on closed irreducible analytic subset of $Y$. Since $Y$ is a smooth projective curve here, these are either of the form $j_*\bb{U}$ where $\bb{U}$ is a local system on a Zariski open subset of $Y$, or skycraper sheaves. \\
  \tab Notice that $\bb{L}' = (\pi')^*\bb{L}$ is $G$-equivariant, and $s'$ is $G-$equivariant hence the same is true for 
  \[K \coloneq \Rm{R}s_*'(\bb{L}'[-\dim X])\] 
  \tab For $g \in G$, $g: Y \to Y$ is an isomorphism hence $g^*$ preserves perverse t-structure. Thus the truncation triangle
  \begin{equation}
    {}^p\tau_{\leq i - 1} K \to {}^p\tau_{\leq i} K \to {}^p\sh{H}^i(K)[-i] \to {}^p\tau_{\leq i - 1} K[1] \label{eq:trunc-triangle}
  \end{equation}
  is also $G$-equivariant. Sabbah's decomposition \eqref{eq:sabbah-decomp} gives a splitting of this triangle. Since $G$ is finite, by averaging we get a $G$-equivariant splitting. Thus we might assume that \eqref{eq:sabbah-decomp} is $G$-equivariant. Each ${}^p\sh{H}^i(K)$ is $G-$equivariant, hence the $G$-action will permute the skyscraper sheaves. Thus the skyscraper sheaves come in $G$-orbits, as stated in the lemma. \\
  \tab Next, the full-support part of ${}^p\sh{H}^i(K)$ has to look like $j_{0, *}(\bb{U}_0[1])$ for some semisimple local system $\bb{U}_0$ on Zariski open $j_0: U_0 \hookrightarrow Y$. Let $U \subseteq U_0$ be a $G$-invariant Zariski open subset of $Y$ with $G$-equivariant inclusion $j: U \hookrightarrow Y$, and put 
  \[\bb{U} = \sh{H}^{-1}(j^{-1}j_{0, *}(\bb{U}[1])) \simeq \bb{U}_0\vert_{U}\]
  then $\bb{U}$ is a $G$-equivariant local system on $U$. It is semisimple since $\pi_1(U) \twoheadrightarrow \pi_1(U_0)$. By uniqueness of intermediate extension, we get a $G$-equivariant isomorphism 
  \[j_*(\bb{U}[1]) \simeq j_{0, *}(\bb{U}_0[1])\]
  \tab For the cohomology statement, by projection formula we have 
  \begin{align*}
    H^i(X, s^*\bb{W} \otimes \bb{L}) &\simeq H^i(X', (\pi')^*s^*\bb{W} \otimes (\pi')^*\bb{L})^{G}\\
    &\simeq H^i(X', (s')^*\pi^*\bb{W} \otimes (\pi')^*\bb{L})^G\\
    &\simeq \HH^i(Y, \pi^*\bb{W} \otimes \Rm{R}(s')_*(\pi')^*\bb{L})^G
  \end{align*}
  and we can apply the decomposition of $\Rm{R}(s')_*(\pi')^*\bb{L}$ to conclude. 
\end{proof}
\begin{prop}
  \label{thm:fin-or-prod}
  $C$ is an irreducible component of $\widetilde{\Sigma}^i_k$. For every $\bb{L} \in T$, the fiber $\widetilde{\Sigma}^i_k\vert_{\bb{L}}$ is either $Q$, a finite set, or empty. Consequently, $C = Q \times T$, or $C \to T$ is generically finite. 
\end{prop}
\begin{proof}
  If an irreducible component $C'$ of $\widetilde{\Sigma}^i_k$ contains $C$, then 
  \[Z = \Phi_{\epsilon}(C) \subseteq \Phi_{\epsilon}(C') \subseteq \ol{\Sigma}_k^i(X, 2)\]
  and $\Phi_{\epsilon}$ is finite hence $\Phi_{\epsilon}(C')$ is irreducible. But $Z$ is an irreducible component, hence $Z = \Phi_{\epsilon}(C')$. Finiteness then gives $\dim C = \dim C'$ thus $C = C'$.

  \medskip
  Next, recall that the fiber $\widetilde{\Sigma}_{k}^i \vert_{\bb{L}}$ is the locus 
  \[\set*{\bb{W} \in Q\ \big\vert\ \dim H^i(s^*\bb{W} \otimes \bb{L}) \geq k}\]
  and by \cref{lem:DM-decomp}, it suffices to look at the jumping behavior of each 
  \[H^{i - a}(Y, (\pi^*\bb{W}) \otimes j_* \bb{U})^{G}, \tab ((\pi^*\bb{W})_{y} \otimes V)^{G_y}\]
  with $j: U \hookrightarrow Y$ the inclusion of a $G$-invariant open, $\bb{U}$ an irreducible $G$-equivariant local system on $U$, $y \in Y \setminus U$ with $G_y$ the stabilizer, and $V$ an irreducible $G_y$-representation. Notice that $j: U \to Y$, $\bb{U}$, $G_y$ and $V$ do not depend on $\bb{W}$. Also, since $Y$ is a smooth curve, $H^{i - a} = 0$ unless $0 \leq i - a \leq 2$. 
  \begin{claim}
    The $G_y$-representation $(\pi^* \bb{W})_y$ does not depend on $\bb{W}$. Consequently, 
    \[\dim ((\pi^*\bb{W})_{y} \otimes V)^{G_y}\]
    is constant as $\bb{W}$ varies in $Q$. 
  \end{claim}  
  \begin{claimproof}
    We have injective $G_y \hookrightarrow \pi_1(\sh{Y}, \bar{y})$ hence each $g \in G_y$ determines an element of $\gamma_g \in \pi_1(\sh{Y}, \bar{y})$. Thus it induces a regular function on $Q$
    \[\bb{W} \mapsto \Rm{tr}(g\vert_{(\pi^* \bb{W})_y}) = \Rm{tr}(\rho_{\bb{W}}(\gamma_g))\]
    where $\rho_{\bb{W}}$ is the $\pi_1(\sh{Y}, \bar{y})$-representation corresponding to $\bb{W}$. Since $G_y$ is finite, the function defined above has discrete image. Since $Q$ is irreducible, it must be constant. It follows that all $G_y$-representations $(\pi^* \bb{W})_y$, as $\bb{W}$ varies in $Q$, have the same characters hence must be isomorphic. 
  \end{claimproof}
  Moving on, notice that $(\pi^*\bb{W}) \otimes j_*\bb{U} = j_*\brac*{(\pi^*\bb{W})\vert_{U} \otimes \bb{U}}$, and denote 
  \[h^l_{G}(\bb{W}, \bb{U}) = \dim H^l\brac*{Y, j_*\brac*{(\pi^*\bb{W})\vert_{U} \otimes \bb{U}}}^{G}\]
  and $\chi_{G}(\bb{W}, \bb{U}) = h^0_G(\bb{W}, \bb{U}) - h^1_G(\bb{W}, \bb{U}) + h^2_G(\bb{W}, \bb{U})$. 
  \begin{claim}
    $\chi_{G}(\bb{W}, \bb{U})$ is constant as $\bb{W}$ varies in $Q$. 
  \end{claim}
  \begin{claimproof}
    Let $\ol{Y} = Y/G$ be the coarse space for $[Y/G]$, with the quotient $q: Y \to \ol{Y}$ which is a finite map. Denote the constructible sheaf 
    \[\sh{L}(\bb{W}, \bb{U}) \coloneq j_*\brac*{(\pi^*\bb{W})\vert_{U} \otimes \bb{U}}, \tab \ol{\sh{L}}(\bb{W}, \bb{U}) \coloneq (q_*\sh{L}(\bb{W}, \bb{U}))^{G}\]
    then $H^l(Y, \sh{L}(\bb{W}, \bb{U}))^{G} \simeq H^l(\ol{Y}, \ol{\sh{L}}(\bb{W}, \bb{U}))$. Denote $I_{\eta} = \ker(G \to \Rm{Aut}(Y))$ the generic inertia. Denote $E \subset \ol{Y}$ a finite set of points, consisting of $q(Y \setminus U)$ and images of points on $Y$ with larger then generic stabilizer. Then $\ol{\sh{L}}(\bb{W}, \bb{U})\vert_{\ol{Y}\setminus E}$ is a local system of dimension 
    \[r(\bb{W}, \bb{U}) = \dim ((\pi^* \bb{W})_{y} \otimes \bb{U}_y)^{G_y} = \dim ((\pi^* \bb{W})_{y} \otimes \bb{U}_y)^{I_{\eta}}, \tab y \in q^{-1}(\ol{Y}\setminus E)\]
    and for $z \in E$, $\ol{\sh{L}}(\bb{W}, \bb{U})\vert_{z}$ is a vector space of dimension 
    \[\dim  ((\pi^* \bb{W})_{y} \otimes (j_*\bb{U})_y)^{G_y}, \tab y \in q^{-1}(z)\]
    \tab By the same argument as in previous claim, both of these values are constant as $\bb{W}$ varies in $Q$. On the other hand, since compactly supported Euler characteristic is additive for the stratification $\ol{Y} = E \sqcup (\ol{Y}\setminus E)$ we get that 
    \begin{align*}
      \chi_{G}(\bb{W}, \bb{U}) &= \chi(\ol{Y}, \ol{\sh{L}}(\bb{W}, \bb{U}))\\
      &= \chi_{c}(\ol{Y}\setminus E, \ol{\sh{L}}(\bb{W}, \bb{U})\vert_{\ol{Y}\setminus E}) + \Sum{z \in E} \dim \ol{\sh{L}}(\bb{W}, \bb{U})\vert_{z}\\
      &= \chi_{c}(\ol{Y}\setminus E) \cdot r(\bb{W}, \bb{U}) + \Sum{z \in E} \dim \ol{\sh{L}}(\bb{W}, \bb{U})\vert_{z} \\
      &= (2 - 2g(\ol{Y}) - \abs*{E}) \cdot r(\bb{W}, \bb{U}) + \Sum{z \in E} \dim \ol{\sh{L}}(\bb{W}, \bb{U})\vert_{z}
    \end{align*}
    which is constant as $\bb{W}$ varies in $Q$, as desired. 
  \end{claimproof}
  Next, we examine these cohomology groups more explicitly: 
  \begin{align*}
    H^0\brac*{Y, j_*\brac*{(\pi^*\bb{W})\vert_{U} \otimes \bb{U}}}^{G} &\simeq H^0(U, (\pi^*\bb{W})\vert_{U} \otimes \bb{U})^{G}\\
    &\simeq \Hom_{\Rm{Loc}^G(U)}(\bb{U}^{\vee}, (\pi^* \bb{W})\vert_{U})
  \end{align*}
  and by Poincar\'{e} duality \cite[Corollary 3.3.12]{Dimca2004Sheaves}
  \[H^2\brac*{Y, j_*\brac*{(\pi^*\bb{W})\vert_{U} \otimes \bb{U}}}^{G} \simeq \Hom_{\Rm{Loc}^G(U)}((\pi^* \bb{W})\vert_{U}, \bb{U}^{\vee})^{\vee}\]
  \tab Since $\bb{W}$ is semisimple on $\sh{Y}$, $\pi^*\bb{W}$ is semisimple in $\Rm{Loc}^G(Y)$. Thus $(\pi^* \bb{W})\vert_{U}$ is semisimple in $\Rm{Loc}^G(U)$. It follows that the two cohomology groups above both count the multiplicities of $\bb{U}^{\vee}$ in $(\pi^* \bb{W})\vert_{U}$ and $h^0_{G}(\bb{W}, \bb{U}) = h^2_{G}(\bb{W}, \bb{U})$. 
  \begin{claim}
    $h^0_G(\bb{W}, \bb{U}) \neq 0$ for only finitely many $\bb{W} \in Q$.  
  \end{claim}
  \begin{claimproof}
    The observation which will be used throughout is that the functor
    \[\Rm{Loc}(\sh{Y}) \to \Rm{Loc}_G(U), \tab \bb{W} \mapsto (\pi^*\bb{W})\vert_{U}\]
    is fully faithful. If $\bb{U}$ has rank $\geq 2$, then if $h^0_G(\bb{W}, \bb{U}) \neq 0$ we must have $(\pi^*\bb{W})\vert_{U} \simeq \bb{U}^{\vee}$. There is at most one such $\bb{W} \in Q$. \\
    \tab Now suppose $\bb{U}$ has rank 1. If $h^0_G(\bb{W}, \bb{U}) \neq 0$ then $\bb{U}^{\vee}$ must be a summand of $(\pi^*\bb{W})\vert_{U}$ and $(\pi^*\bb{W})\vert_{U}$ must be reducible in $\Rm{Loc}^{G}(U)$. It follows that $(\pi^*\bb{W})$ must be reducible in $\Rm{Loc}^G(Y)$. Since $\bb{W} \in Q \subseteq M_B(\sh{Y}, 2, \det = \delta, \bm{\tau})$, $\pi^*\bb{W}$ has fixed determinant $\pi^*\delta$. It follows that 
    \[\pi^* \bb{W} \simeq \bb{M} \oplus (\bb{M}^{\vee} \otimes \pi^*\delta)\]
    and $\bb{U}^{\vee} = \bb{M}$ or $\bb{M}^{\vee} \otimes \pi^*\delta$. Again, finitely many $\bb{W}$ giving this. 
  \end{claimproof}
  Each $\bb{U}$ in \cref{lem:DM-decomp} gives a finite set of $\bb{W} \in Q$ where $h^0_G(\bb{W}, \bb{U}) \neq 0$. Take $F(\bb{L})$ to be the union as $\bb{U}$ runs over all the finitely many summands in the decomposition of $\Rm{R}s_*'(\pi')^*\bb{L}$, then $F(\bb{L})$ is also finite. To summarize, for a fixed $\bb{L}$
  \begin{align*}
    \dim H^0\brac*{Y, (\pi^*\bb{W}) \otimes j_*\bb{U}}^{G} &= h^0_G(\bb{W}, \bb{U}) = 0 \tab \forall\ \bb{W} \in Q\setminus F(\bb{L})\\
    \dim H^2\brac*{Y, (\pi^*\bb{W}) \otimes j_*\bb{U}}^{G} &= h^2_G(\bb{W}, \bb{U}) = h^0_G(\bb{W}, \bb{U}) \tab \forall\ \bb{W} \in Q\setminus F(\bb{L})
  \end{align*}
  and thus for all $\bb{W} \in Q \setminus F(\bb{L})$,
  \begin{align*}
    \dim H^1\brac*{Y, (\pi^*\bb{W}) \otimes j_*\bb{U}}^{G} &= -\chi_G(\bb{W}, \bb{U}) + 2h^0_G(\bb{W}, \bb{U})\\
    &= -\chi_G(\bb{W}, \bb{U}) 
  \end{align*}
  which is constant for $\bb{W} \in Q$. It follows that there is a constant $k_{\bb{L}}$ such that 
  \begin{align}
    \dim H^i(X, s^*\bb{W} \otimes \bb{L}) &= k_{\bb{L}} \tab \forall\ \bb{W} \in Q\setminus F(\bb{L}) \label{eq:generic-bound}\\
    \dim H^i(X, s^*\bb{W} \otimes \bb{L}) &\geq k_{\bb{L}} \tab \forall\ \bb{W} \in Q 
  \end{align}
  and consequently we have proven:
  \begin{claim}
    For every $\bb{L} \in T$, the fiber $\widetilde{\Sigma}^i_k\vert_{\bb{L}}$ is either $Q$, finite, or empty. 
  \end{claim}
  Back to $C \subseteq \widetilde{\Sigma}^i_k \subseteq Q \times T$. The map $C \to T$ is dominant by definition. Here we can assume $\dim Q > 0$ by \cref{rmk:proj-nonrigidity}. Suppose that the general fiber of $C \to T$ has positive dimension, then on a Zariski open subset $T^{\circ} \subseteq T$, since $\widetilde{\Sigma}^i_k \supset C$, $\widetilde{\Sigma}^i_k\vert_{\bb{L}}$ has fiber $Q$. Hence,
  \[Q \times T^{\circ} \subset \widetilde{\Sigma}^i_k \subset Q \times T\]
  and since $\widetilde{\Sigma}^i_k$ is closed, we get $\widetilde{\Sigma}^i_k = Q \times T$. This is irreducible, so $C = Q \times T$. The other case is that $C \to T$ is generically finite. 
\end{proof}

\subsection{Finiteness of $C \to T$} 
\label{sec:finite-map}
In this section, we will prove that if $C \xrightarrow{\Rm{pr}_2} T$ is generically finite, then it must in fact be finite. To differentiate, let the $G-$invariant dense open in \cref{lem:DM-decomp} be denoted $j_{\bb{L}}: U_{\bb{L}} \hookrightarrow Y$ for the decomposition of $\Rm{R}s_*'(\pi')^* \bb{L}$. Also denote $\bb{L}' = (\pi')^*\bb{L}$. \\
\tab For $U \subseteq Y$ a $G-$equivariant Zariski open subset, we have a diagram by restricting \eqref{eq:pullback-DM}
\[\begin{tikzcd}
	{X_{U}'} && {X_U} \\
	\\
	U && {[U/G]}
	\arrow["{\pi_{U}'}", from=1-1, to=1-3]
	\arrow["{s_{U}'}"', from=1-1, to=3-1]
	\arrow["{s_U}", from=1-3, to=3-3]
	\arrow["{\pi_{U}}"', from=3-1, to=3-3]
\end{tikzcd}\]
and we can pick $U$ so that the map $s_{U}': X_{U}' \to U$ is smooth and proper. Then 
\begin{equation}
  \sh{H}^m_{\bb{L}} = \Rm{R}^m(s_{U}')_*\brac*{\bb{L}'\vert_{X_{U}'}} \label{eq:H^m_L-locsys}
\end{equation}
is a $G$-equivariant local system on $U$. 
\begin{lemma}
  \label{lem:equiv-quot}
  Let $\bb{L} \in T$, and suppose that the fiber $\widetilde{\Sigma}^i_k\vert_{\bb{L}}$ is finite. Let $(\bb{W}, \bb{L}) \in C \subseteq \widetilde{\Sigma}^i_k$ with $\bb{W}$ irreducible in $\Rm{Loc}(\sh{Y})$, then, for some $m$, there is a quotient 
  \[\sh{H}^m_{\bb{L}} \twoheadrightarrow ((\pi^*\bb{W})\vert_{U})^{\vee}\]
  in $\Rm{Loc}^{G}(U)$. In fact, $(\pi^*\bb{W})\vert_{U}^{\vee}$ is a direct summand of $\sh{H}^m_{\bb{L}}$. 
\end{lemma}
\begin{proof}
  Look at the $G$-equivariant decomposition of \cref{lem:DM-decomp}, 
  \begin{equation}
    \Rm{R}(s_*')\bb{L}' \simeq \bigoplus_{\alpha} j_{\bb{L}, *} \bb{U}_{\alpha}[-a_{\alpha}] \oplus \bigoplus_{\beta} \iota_{\beta, *}\sh{V}_{\beta}[-b_{\beta}] \label{eq:L-decomp}
  \end{equation}
  with $\bb{U}_{\alpha}$ irreducible in $\Rm{Loc}^{G}(U_{\bb{L}})$. By the proof of \cref{thm:fin-or-prod} (see in particular the paragraphs around \eqref{eq:generic-bound}), if $\widetilde{\Sigma}^i_k\vert_{\bb{L}}$ is finite, then it consists of points $(\bb{W}, \bb{L})$ such that 
  \[\dim H^0(Y, (\pi^*\bb{W}) \otimes j_{\bb{L}, *} \bb{U}_{\alpha})^{G} > 0\]
  which implies that $\Hom_{\Rm{Loc}^G(U_{\bb{L}})}(\bb{U}^{\vee}_{\alpha}, (\pi^* \bb{W})\vert_{U_{\bb{L}}}) > 0$. Since $\bb{W} \mapsto (\pi^* \bb{W})\vert_{U_{\bb{L}}}$ is a fully faithful functor which preserves irreducibility, we have that 
  \[\bb{U}_{\alpha} \simeq ((\pi^* \bb{W})\vert_{U_{\bb{L}}})^{\vee} \tab \text{ in }\Rm{Loc}^G(U_{\bb{L}})\]
  \tab Next, by restricting the decomposition \eqref{eq:L-decomp} to $U$ we get that $(j_{\bb{L}, *} \bb{U}_{\alpha})\vert_{U}$ is a direct summand of $\sh{H}^{a_{\alpha}}_{\bb{L}}$ in $\Rm{Loc}^G(U)$; in particular, $(j_{\bb{L}, *} \bb{U}_{\alpha})\vert_{U}$ is a local system on $U$. The previous argument show that 
  \[(j_{\bb{L}, *} \bb{U}_{\alpha})\vert_{U \cap U_{\bb{L}}} \simeq \bb{U}_{\alpha}\vert_{U \cap U_{\bb{L}}} \simeq ((\pi^*\bb{W})\vert_{U \cap U_{\bb{L}}})^{\vee}\]
  and since restriction $\Rm{Loc}^G(U) \to \Rm{Loc}^G(U \cap U_{\bb{L}})$ is fully faithful, we get that 
  \[(j_{\bb{L}, *} \bb{U}_{\alpha})\vert_{U} \simeq ((\pi^* \bb{W})\vert_{U})^{\vee}\tab \text{ in }\Rm{Loc}^G(U)\]
  and thus $((\pi^* \bb{W})\vert_{U})^{\vee}$ is a direct summand of $\sh{H}^{a_{\alpha}}_{\bb{L}}$. Since $\sh{H}^{a_{\alpha}}_{\bb{L}}$ is semisimple we get the quotient statement as well. 
\end{proof}
Let $\sh{O}(T) = \Gamma(T, \sh{O}_T)$. Pulling back the universal rank 1 local system along $T \subseteq M_B(X, 1)$ we get $\bb{L}_{\univ}$ a local system of free $\sh{O}(T)$-modules on $X$. For $\bb{L} \in T$, 
\[\bb{L}_{\univ} \otimes_{\sh{O}(T)} \kappa(\bb{L}) \simeq \bb{L}\]
where we are taking tensor products with, and over, constant local systems on $X$. Again, denote $\bb{L}_{\univ}' = (\pi')^*\bb{L}_{\univ}$, then $\bb{L}_{\univ}'$ is a $G$-equivariant local system of $\sh{O}(T)$-modules on $X'$ (the $G$-action on $T$ is trivial). We also have that 
\[\sh{H}^m_{\univ} = \Rm{R}^m(s_{U}')_*(\bb{L}_{\univ}'\vert_{X_{U}'})\]
is a $G$-equivariant local system of $\sh{O}(T)$-modules on $U$, since $s_{U}'$ is a smooth and proper $G$-equivariant map. Equivalently, fixing $u \in U$ with $\bar{u}$ the image in $[U/G]$, $(\sh{H}^m_{\univ})_{u}$ is a $\sh{O}(T)$-module, and we have a $\sh{O}(T)$-linear representation 
\[\pi_1([U/G], \bar{u}) \to \Rm{GL}_{\sh{O}(T)}((\sh{H}^m_{\univ})_{u})\]
\tab On the other hand, we can also think of $(\sh{H}^m_{\univ})_u$ as a $\sh{O}_T$-module on $T$. There is then a Zariski dense open $T^{\Rm{bc}} \subseteq T$ where cohomology and base change is available, i.e., for $\bb{L} \in T^{\Rm{bc}}$, there is an isomorphisms of $\pi_1([U/G], \bar{u})$-representations
\[(\sh{H}^m_{\univ})_{u} \otimes \kappa(\bb{L}) \simeq (\sh{H}^m_{\bb{L}})_u\]
with $\sh{H}^m_{\bb{L}}$ defined in \eqref{eq:H^m_L-locsys}. 
\begin{remark}
  \label{rmk:equiv-quot}
  By \cref{lem:equiv-quot}, a point $(\bb{W}, \bb{L})$ in a finite fiber of $C \to T$ with $\bb{L} \in T^{\Rm{bc}}$ and $\bb{W}$ irreducible correspond to a rank two quotient
  \[(\sh{H}^m_{\univ})_{u} \otimes \kappa(\bb{L}) \twoheadrightarrow ((\pi^*\bb{W})\vert_{U})^{\vee}_{u}\]
  as $\pi_1([U/G], \bar{u})$-representations. 
\end{remark}
\begin{prop}[{\cite[Corollary 2.4]{HaimanSturmfels2004}}]
  \label{thm:equiv-grass}
  Let $T$ be an affine noetherian variety, $\sh{H}$ a coherent $\sh{O}_T$-module, and $\Pi$ a finitely generated group. Suppose that there is a $\sh{O}_T$-linear action of $\Pi$ on $\sh{H}$. The functor assigning to an affine $T$-scheme $\phi: S \to T$ the set 
  \[\set*{\Pi-\text{invariant locally free rank 2 quotients of }\phi^*\sh{H}}\] 
  is representable by a projective $T$-scheme. 
\end{prop}
\begin{proof}
  Since we are on affine schemes, we will abuse notation between sheaves and modules. Denote the $\sh{O}(T)$-linear action of $\Pi$ on $\sh{H}$ by
  \[\rho: \Pi \to \Rm{GL}_{\sh{O}(T)}(\sh{H})\]
  then our functor assigns to a commutative $\sh{O}(T)$-algebra $R$ the set of $\rho(\Pi)-$invariant locally free $R$-module quotients of $\sh{H} \otimes_{\sh{O}(T)} R$. To apply the cited result, let $k = \sh{O}(T)$, $A$ be the singleton set, $T = \sh{H}$, and $F = \rho(\Pi)$.  
\end{proof}
Apply the previous proposition to $\Pi = \pi_1([U/G], \bar{u})$ and $\sh{H} = (\sh{H}^m_{\univ})_u$. Let the resulting projective $T$-scheme be denoted $\scr{G}_m$. Notice that a closed point of $\scr{G}_m$, over $\bb{L} \in T^{\Rm{bc}}$, is a $\pi_1([U/G], \bar{u})$-invariant quotient of $(\sh{H}^m_{\univ})_{u} \otimes \kappa(\bb{L}) \simeq (\sh{H}^m_{\bb{L}})_u$. If $\bb{L} \notin T^{\Rm{bc}}$ then the quotient is still a rank 2 $\pi_1([U/G], \bar{u})$-representation. We have a natural map over $T$
\[\Phi_m: \scr{G}_m \to M_B([U/G], 2) \times T, \tab [(\sh{H}^m_{\univ})_{u} \otimes \kappa(\bb{L}) \twoheadrightarrow \bb{V}] \mapsto ([\bb{V}^{\vee}], \bb{L})\]
which has image finite over $T$ since $M_B([U/G], 2)$ is affine. 
\begin{prop}
  \label{thm:finite-closure}
  There is a closed subscheme $\scr{F}_Q \subset Q \times T$, finite over $T$, with the following property: if $\bb{W} \in Q$ and $\bb{L} \in T^{\Rm{bc}}$, and $((\pi^*\bb{W})\vert_{U})^{\vee}$ is a rank two $G$-equivariant quotient of $\sh{H}_{\bb{L}}^m$ for some $m$, then $(\bb{W}, \bb{L}) \in \scr{F}_Q$. 
\end{prop}
\begin{proof}
  Denote the disjoint union
  \[\scr{G} = \bigsqcup_{m = 0}^{2\dim X - 1} \scr{G}_m\]
  and let $\scr{F}$ be the scheme-theoretic image of $\Phi: \scr{G} \to M_B([U/G], 2) \times T$. Then $\scr{F}$ is finite over $T$. We have a closed embedding 
  \[M_B(\sh{Y}, 2) \hookrightarrow M_B([U/G], 2), \tab \bb{W} \mapsto (\pi^*\bb{W})\vert_U\]
  hence we can view $Q \times T$ as a closed subscheme of $M_B([U/G], 2) \times T$, and define 
  \[\scr{F}_Q \coloneq \scr{F} \times_{M_B([U/G], 2) \times T} (Q \times T)\]
  which is finite over $T$. On the other hand, a point $(\bb{W}, \bb{L})$ with $\bb{L} \in T^{\Rm{bc}}$ and $((\pi^*\bb{W})\vert_{U})^{\vee}$ is a rank two $G$-equivariant quotient of $\sh{H}_{\bb{L}}^m$ for some $m$ must be in the image of $\Phi_m$ by \cref{rmk:equiv-quot}. As a result it must be in $\scr{F}_Q$. 
\end{proof}
\begin{prop}
  \label{thm:not-gen-fin}
  If $\Rm{pr}_2: C \to T$ is generically finite, then it must be finite and surjective.
\end{prop}
\begin{proof}
  There is an Zariski open subset $T^{\circ} \subseteq T$ such that 
  \begin{itemize}
    \item $T^{\circ} \subseteq T^{\Rm{bc}}$;
    \item $\widetilde{\Sigma}^i_k\vert_{\bb{L}}$ is finite or empty for every $\bb{L} \in T^{\circ}$;
    \item every $(\bb{W}, \bb{L}) \in C$ with $\bb{L} \in T^{\circ}$ has $\bb{W}$ irreducible. 
  \end{itemize}
  Let $(\bb{W}, \bb{L})$ be a closed point of $C \times_{T} T^{\circ}$. Then $\widetilde{\Sigma}^i_k\vert_{\bb{L}}$ is a finite fiber. By \cref{lem:equiv-quot} we get a quotient 
  \[\sh{H}_{\bb{L}}^m \twoheadrightarrow ((\pi^*\bb{W})\vert_{U})^{\vee}\]
  for some $m$. By \cref{thm:finite-closure}, $(\bb{W}, \bb{L}) \in \scr{F}_Q$. It follows that 
  \[C \times_{T} T^{\circ} \subseteq \scr{F}_Q\]
  \tab On the other hand, since $\Rm{pr}_2$ is dominant, $C \times_{T} T^{\circ}$ is dense in $C$. Taking closures in $Q \times T$ then gives $C$ is a closed subscheme of $\scr{F}_Q$. Since $\scr{F}_Q$ is finite over $T$, we get that $C \to T$ is finite. Since it is dominant, it must be surjective as well. 
\end{proof}

\subsection{\'{E}taleness of $C^{\nu} \to T$} 
\label{sec:etaleness}
We have that $\Rm{pr}_2: C \to T$ is finite. Let $\nu: C^{\nu} \to C$ be the normalization, and $\Rm{pr}_2^{\nu} \coloneq \Rm{pr}_2 \circ \nu: C^{\nu} \to T$. Then $\Rm{pr}_2^{\nu}$ is also finite and surjective. Let $C^{\text{\'{e}t}} \subseteq C$ be Zariski open locus where $\Rm{pr}_2$ is \'{e}tale, and denote
\[\Delta(\Rm{pr}_2) \coloneq \Rm{pr}_2(C \setminus C^{\text{\'{e}t}})_{\red} \subset T\]
which is closed in $T$. 
\begin{lemma}
  \label{lem:non-et=abs}
  The set $\Delta(\Rm{pr}_2)$ is absolute $\ol{\QQ}$-closed in $A$.
\end{lemma}
\begin{proof}
  See \eqref{eq:det-sheet} for the setup. By the first paragraph of the proof of \cref{lem:proj=torus}, but with $M_B$, we have that $C \subseteq M_B(X, 2, \det = \epsilon) \times A$ and $T \subseteq A$ are both absolute $\ol{\QQ}$-closed. In the notation of \cref{def:absolute-set}, $C^{\sigma}$ and $T^{\sigma}$ are both $\ol{\QQ}$-closed for any $\sigma \in \Rm{Aut}(\CC/\QQ)$. By \cite[Proposition 7.8]{simpson1994moduli2} we have an analytic isomorphism 
  \[\Sf{RH}: M_{\DR}(X, r)^{\Rm{an}} \xrightarrow{\simeq} M_B(X, r)^{\Rm{an}}\]
  and $C_{\DR} = \Sf{RH}^{-1}(C), T_{\DR} = \Sf{RH}^{-1}(T)$ are both algebraic by \cite[Proposition 7.4.5]{budur2020absolutesets}. We also have diagrams 
  \begin{equation}
    \label{eq:RH-square}
  \begin{tikzcd}
	{(C_{\DR})^{\Rm{an}}} && {(T_{\DR})^{\Rm{an}}} && {(C_{\DR}^{\sigma})^{\Rm{an}}} && {(T_{\DR}^{\sigma})^{\Rm{an}}} \\
	{} \\
	{C^{\Rm{an}}} && {T^{\Rm{an}}} && {(C^{\sigma})^{\Rm{an}}} && {(T^{\sigma})^{\Rm{an}}}
	\arrow["{(\Rm{pr}_{\DR, 2})^{\Rm{an}}}", from=1-1, to=1-3]
	\arrow["{\Sf{RH}}", from=1-1, to=3-1]
	\arrow["{\Sf{RH}}", from=1-3, to=3-3]
	\arrow["{(\Rm{pr}_{\DR, 2}^{\sigma})^{\Rm{an}}}", from=1-5, to=1-7]
	\arrow["{\Sf{RH}}"', from=1-5, to=3-5]
	\arrow["{\Sf{RH}}"', from=1-7, to=3-7]
	\arrow["{(\Rm{pr}_2)^{\Rm{an}}}", from=3-1, to=3-3]
	\arrow["{(\Rm{pr}_2^{\sigma})^{\Rm{an}}}", from=3-5, to=3-7]
  \end{tikzcd}
  \end{equation}
  where $(\Rm{pr}_{\DR, 2}^{\sigma})^{\Rm{an}}$ is induced by $(\Rm{pr}_{\DR, 2})^{\Rm{an}}$ through an isomorphism of schemes
  \[\sigma: M_{\DR}(X, r) \xrightarrow{\simeq} M_{\DR}(X^{\sigma}, r)\] 
  which is itself the pullback along the base change $X^{\sigma} \to X$. \\
  \tab We have that analytification detects \'{e}taleness \cite[Last paragraph of \href{https://stacks.math.columbia.edu/tag/025L}{Tag 025L}]{stacks-project}, i.e. an algebraic map is \'{e}tale if and only if its analytification is a local isomorphism. We also have that $(\Rm{pr}_{\DR, 2}^{\sigma})^{\Rm{an}}$ is a local isomorphism if and only if $(\Rm{pr}_{\DR, 2})^{\Rm{an}}$ is one. It follows that $\Delta(\Rm{pr}_2)^{\sigma} = \Delta(\Rm{pr}_2^{\sigma})$. Thus $\Delta(\Rm{pr}_2)$ is absolute $\ol{\QQ}$-closed. 
\end{proof}

\begin{lemma}
  \label{lem:ram-div}
  If $\Rm{pr}_2^{\nu}$ is not \'{e}tale, then there is a prime divisor $E \subset C^{\nu}$ such that $\Rm{pr}_2^{\nu}$ is ramified at the generic point $\eta_{E}$ of $E$. Furthermore, $\Rm{pr}_2^{\nu}(E)$ is an irreducible component of $\Delta(\Rm{pr}_2)$.
\end{lemma}
\begin{proof}
  The existence of $\eta_{E}$ follows from (5) of \cite[\href{https://stacks.math.columbia.edu/tag/0BMB}{Tag 0BMB}]{stacks-project} since $\Rm{pr}_2^{\nu}$ satisfies all (1-4). Since $\Rm{pr}_2^{\nu}$ is ramified at $\eta_{E}$, $\Rm{pr}_2$ must not be \'{e}tale at $\nu(\eta_{E})$ because otherwise $C$ is smooth at $\nu(\eta_{E})$ (since $T$ is smooth) and $\nu$ is an isomorphism at $\eta_{E}$. Thus $\Rm{pr}_2^{\nu}(E) \subseteq \Delta(\Rm{pr}_2)$. A finite map preserves dimensions, so $\Rm{pr}_2^{\nu}(E)$ is a prime divisor in $T$ hence must be an irreducible component of $\Delta(\Rm{pr}_2)$. 
\end{proof}

\begin{prop}
  \label{thm:etaleness}
  The finite map $\Rm{pr}_2^{\nu}: C^{\nu} \to T$ is \'{e}tale. 
\end{prop}
\begin{proof}
  Suppose otherwise, then by \cref{lem:ram-div} we get a prime divisor $E \subset C^{\nu}$ such that $\Rm{pr}_2^{\nu}$ is ramified at every point of $E$, and $\Rm{pr}_2^{\nu}(E) = T_0$ is an irreducible component of $\Delta(\Rm{pr}_2)$. Since $\Delta(\Rm{pr}_2)$ is absolute $\ol{\QQ}$-closed, the same is true for $T_0$ and by \cref{thm:rank1-absolute}, $T_0$ is a torsion-translated subtorus of $A$.\\
  \tab Since $C^{\nu}$ is normal, the singular locus is codimension 2, and $C^{\nu}$ is smooth on a Zariski open subset of $E$. $T_0$ is a torsion-translated subtorus, hence unitary points are Zariski dense in $T_0$. It follows that there exists a point $c \in E \subset C^{\nu}$ lying over $(\bb{V}, \bb{L})$ such that $\Rm{pr}_2^{\nu}(c) = \bb{L} \in T_0$ is unitary. We have a $\CC^*$-action on $M_B(X, 2) \times A \supset C$. Now, $C$ is the good moduli image of a Hodge substack of $M_B(X, 2) \times A$, hence by \cite[Corollary 8.24]{bakker2024linear}, $C$ is $\CC^*$-stable. Now look at the $\CC^*$-orbit of $(\bb{V}, \bb{L})$. Since $\bb{L}$ is unitary, it is $\CC^*$-fixed, hence the image under $\Rm{pr}_2$ of this orbit is just $\bb{L}$. But $\Rm{pr}_2: C \to T$ is finite, so the $\CC^*$-orbit is finite, thus $(\bb{V}, \bb{L})$ must be $\CC^*$-fixed to begin with. To summarize, we have $c \in C^{\nu}$ a smooth point lying over $(\bb{V}, \bb{L})$ with both $\bb{V}, \bb{L}$ underlying weight 0 variation of complex Hodge structures. By \cref{lem:normal-Hodge} and \cref{lem:GIT-proj=torus} we get a morphism of pro-$\Bf{T}(0)$-MHS-algebras
  \[\Rm{pr}_{2, c}^{\nu, \#}: \widehat{\sh{O}}_{T, \bb{L}} \to \widehat{\sh{O}}_{C, (\bb{V}, \bb{L})} \to \widehat{\sh{O}}_{C^{\nu}, c}\]
  which is a finite local map between regular complete local rings of the same dimension. By \cite[\href{https://stacks.math.columbia.edu/tag/00R4}{Tag 00R4}]{stacks-project} it must be flat, hence faithfully flat. In particular $\Rm{pr}_{2, c}^{\nu, \#}$ is injective. Next look at the map on cotangent spaces 
  \[\der(\Rm{pr}_2^{\nu})^{\vee}_c: \Ide{m}_{T, \bb{L}}/\Ide{m}_{T, \bb{L}}^2 \to \Ide{m}_{C^{\nu}, c}/\Ide{m}_{C^{\nu}, c}^2\]
  which must have nontrivial kernel since $\Rm{pr}_2^{\nu}$ is ramified at $c$. Take the functorial Deligne weight splitting (similar to the proof of \cref{lem:finite-weight}) to talk about homogeneous elements. Notice that $\Ide{m}_{T, \bb{L}}/\Ide{m}_{T, \bb{L}}^2$ is pure of weight $-1$. Pick a homogeneous weight $-1$ element in the kernel of $\der(\Rm{pr}_2^{\nu})^{\vee}_c$. Its image must either be 0 or also homogeneous of weight $-1$. The latter is impossible since $\Ide{m}^2_{C^{\nu}, c} \subseteq W_{-2}\widehat{\sh{O}}_{C^{\nu}, c}$, and the former is also impossible since $\widehat{\sh{O}}_{T, \bb{L}} \to \widehat{\sh{O}}_{C^{\nu}, c}$ is injective. Thus we reach a contradiction, and $\Rm{pr}_2^{\nu}: C^{\nu} \to T$ is \'{e}tale. 
\end{proof}

\begin{coroll}
  \label{lem:Cnu=torus}
  $C^{\nu}$ is an algebraic torus. 
\end{coroll}

Notice that $H_C$ (as in \cref{lem:normalized-det-sheet}) acts on $T$ by translations, and the map $C^{\nu} \to T$ is $H_C$-equivariant, hence $Z^{\nu}/H_C$ is finite \'{e}tale over $T/H_C$, which is an algebraic torus. Thus we end up with:

\begin{coroll}
  $Z^{\nu} \simeq C^{\nu}/H_C$ is an algebraic torus. 
\end{coroll}

\printbibliography

@book{hartshorne1977algebraic,
  author        = {Hartshorne, Robin},
  title         = {Algebraic geometry},
  series        = {Graduate Texts in Mathematics},
  number        = {52},
  publisher     = {Springer-Verlag},
  address       = {New York},
  year          = {1977}
}

@article{budur2015cohomology,
  author        = {Budur, Nero and Wang, Botong},
  title         = {Cohomology jump loci of differential graded {Lie} algebras},
  journal       = {Compositio Mathematica},
  volume        = {151},
  number        = {8},
  pages         = {1499--1528},
  year          = {2015}
}

@article{simpson1992higgs,
  author        = {Simpson, Carlos T.},
  title         = {{Higgs} bundles and local systems},
  journal       = {Publications Math{\'e}matiques de l'IH{\'E}S},
  volume        = {75},
  pages         = {5--95},
  year          = {1992}
}

@article{corlette1988flat,
  author        = {Corlette, Kevin},
  title         = {Flat {$G$}-bundles with canonical metrics},
  journal       = {Journal of Differential Geometry},
  volume        = {28},
  number        = {3},
  pages         = {361--382},
  year          = {1988}
}

@article{saito1990mixed,
  author        = {Saito, Morihiko},
  title         = {Mixed {Hodge} modules},
  journal       = {Publications of the Research Institute for Mathematical Sciences},
  volume        = {26},
  number        = {2},
  pages         = {221--333},
  year          = {1990}
}

@book{sabbah2005polarizable,
  author        = {Sabbah, Claude},
  title         = {Polarizable twistor {$\mathcal{D}$}-modules},
  series        = {Ast{\'e}risque},
  number        = {300},
  pagetotal     = {vi+208},
  publisher     = {Soci{\'e}t{\'e} Math{\'e}matique de France},
  address       = {Paris},
  year          = {2005}
}

@book{mochizuki2015mixed,
  author        = {Mochizuki, Takuro},
  title         = {Mixed twistor {$D$}-modules},
  series        = {Lecture Notes in Mathematics},
  number        = {2125},
  publisher     = {Springer},
  address       = {Cham},
  year          = {2015}
}

@article{simpson1994moduli2,
  author        = {Simpson, Carlos T.},
  title         = {Moduli of representations of the fundamental group of a smooth projective variety {II}},
  journal       = {Publications Math{\'e}matiques de l'IH{\'E}S},
  volume        = {80},
  pages         = {5--79},
  year          = {1994}
}

@article{simpson1993subspaces,
  author        = {Simpson, Carlos},
  title         = {Subspaces of moduli spaces of rank one local systems},
  journal       = {Annales scientifiques de l'{\'E}cole Normale Sup{\'e}rieure},
  series        = {4},
  volume        = {26},
  number        = {3},
  pages         = {361--401},
  year          = {1993}
}

@article{budur2020absolutesets,
  author        = {Budur, Nero and Wang, Botong},
  title         = {Absolute sets and the decomposition theorem},
  journal       = {Annales scientifiques de l'{\'E}cole Normale Sup{\'e}rieure},
  series        = {4},
  volume        = {53},
  number        = {2},
  pages         = {469--536},
  year          = {2020}
}

@misc{bakker2024linear,
  author        = {Bakker, Benjamin and Brunebarbe, Yohan and Tsimerman, Jacob},
  title         = {The linear {Shafarevich} conjecture for quasiprojective varieties and algebraicity of {Shafarevich} morphisms},
  year          = {2024},
  eprint        = {2408.16441},
  archivePrefix = {arXiv},
  primaryClass  = {math.AG}
}

@incollection{simpson1996hodgefilt,
  author        = {Simpson, Carlos},
  title         = {The {Hodge} filtration on nonabelian cohomology},
  editor        = {Koll{\'a}r, J{\'a}nos and Lazarsfeld, Robert and Morrison, David R.},
  booktitle     = {Algebraic geometry---{Santa Cruz} 1995, Part 2},
  series        = {Proceedings of Symposia in Pure Mathematics},
  number        = {62},
  pages         = {217--281},
  publisher     = {American Mathematical Society},
  address       = {Providence, RI},
  year          = {1997}
}

@misc{simpson1997mixedtwistor,
  author        = {Simpson, Carlos},
  title         = {Mixed twistor structures},
  year          = {1997},
  eprint        = {alg-geom/9705006},
  archivePrefix = {arXiv},
  primaryClass  = {math.AG}
}

@misc{stacks-project,
  author       = {The {Stacks project authors}},
  title        = {The Stacks project},
  howpublished = {\url{https://stacks.math.columbia.edu}},
  year         = {2026},
}

@article{nguyen2008algebraic,
  author        = {Nguyen, K. A. and van der Put, Marius and Top, Jaap},
  title         = {Algebraic subgroups of {$\mathrm{GL}_2(\mathbb{C})$}},
  journal       = {Indagationes Mathematicae},
  series        = {N.S.},
  volume        = {19},
  number        = {2},
  pages         = {287--297},
  year          = {2008}
}

@article{alper2020luna,
  author        = {Alper, Jarod and Hall, Jack and Rydh, David},
  title         = {A {Luna} {\'e}tale slice theorem for algebraic stacks},
  journal       = {Annals of Mathematics},
  series        = {2},
  volume        = {191},
  number        = {3},
  pages         = {675--738},
  year          = {2020}
}

@incollection{luna1973slices,
  author        = {Luna, Domingo},
  title         = {Slices {\'e}tales},
  booktitle     = {Sur les groupes alg{\'e}briques},
  series        = {Bulletin de la Soci{\'e}t{\'e} math{\'e}matique de France. M{\'e}moire},
  number        = {33},
  pages         = {81--105},
  publisher     = {Soci{\'e}t{\'e} Math{\'e}matique de France},
  address       = {Paris},
  year          = {1973}
}

@book{PetersSteenbrink2008,
  author        = {Peters, Chris A. M. and Steenbrink, Joseph H. M.},
  title         = {Mixed {Hodge} structures},
  series        = {Ergebnisse der Mathematik und ihrer Grenzgebiete. 3. Folge},
  number        = {52},
  publisher     = {Springer-Verlag},
  address       = {Berlin},
  year          = {2008}
}

@book{atiyah1969introduction,
  author        = {Atiyah, Michael F. and Macdonald, Ian G.},
  title         = {Introduction to commutative algebra},
  publisher     = {Addison-Wesley Publishing Co.},
  address       = {Reading, MA},
  year          = {1969}
}

@book{grauert1984coherent,
  author        = {Grauert, Hans and Remmert, Reinhold},
  title         = {Coherent analytic sheaves},
  series        = {Grundlehren der mathematischen Wissenschaften},
  number        = {265},
  publisher     = {Springer-Verlag},
  address       = {Berlin},
  year          = {1984}
}

@article{corlette2008classification,
  author        = {Corlette, Kevin and Simpson, Carlos},
  title         = {On the classification of rank-two representations of quasiprojective fundamental groups},
  journal       = {Compositio Mathematica},
  volume        = {144},
  number        = {5},
  pages         = {1271--1331},
  year          = {2008}
}

@article{noohi2004fundamental,
  author        = {Noohi, Behrang},
  title         = {Fundamental groups of algebraic stacks},
  journal       = {Journal of the Institute of Mathematics of Jussieu},
  volume        = {3},
  number        = {1},
  pages         = {69--103},
  year          = {2004}
}

@misc{noohi2005foundationsi,
  author        = {Noohi, Behrang},
  title         = {Foundations of topological stacks {I}},
  year          = {2005},
  eprint        = {math/0503247},
  archivePrefix = {arXiv},
  primaryClass  = {math.AG}
}

@incollection{simpson1991ubiquity,
  author        = {Simpson, Carlos T.},
  title         = {The ubiquity of variations of {Hodge} structure},
  editor        = {Carlson, James A. and Clemens, C. Herbert and Morrison, David R.},
  booktitle     = {Complex geometry and {Lie} theory},
  series        = {Proceedings of Symposia in Pure Mathematics},
  number        = {53},
  pages         = {329--348},
  publisher     = {American Mathematical Society},
  address       = {Providence, RI},
  year          = {1991}
}

@book{Dimca2004Sheaves,
  author        = {Dimca, Alexandru},
  title         = {Sheaves in topology},
  series        = {Universitext},
  publisher     = {Springer-Verlag},
  address       = {Berlin},
  year          = {2004}
}

@article{HaimanSturmfels2004,
  author        = {Haiman, Mark and Sturmfels, Bernd},
  title         = {Multigraded {Hilbert} schemes},
  journal       = {Journal of Algebraic Geometry},
  volume        = {13},
  number        = {4},
  pages         = {725--769},
  year          = {2004}
}

@article{Alper2012,
  author        = {Alper, Jarod},
  title         = {Local properties of good moduli spaces},
  journal       = {Tohoku Mathematical Journal},
  series        = {2},
  volume        = {64},
  number        = {1},
  pages         = {105--123},
  year          = {2012}
}

@article{GreenLazarsfeld1987,
  author        = {Green, Mark and Lazarsfeld, Robert},
  title         = {Deformation theory, generic vanishing theorems, and some conjectures of {E}nriques, {C}atanese and {B}eauville},
  journal       = {Inventiones Mathematicae},
  volume        = {90},
  number        = {2},
  pages         = {389--407},
  year          = {1987}
}

@article{GreenLazarsfeld1991,
  author        = {Green, Mark and Lazarsfeld, Robert},
  title         = {Higher obstructions to deforming cohomology groups of line bundles},
  journal       = {Journal of the American Mathematical Society},
  volume        = {4},
  number        = {1},
  pages         = {87--103},
  year          = {1991}
}

@article{Arapura1997,
  author        = {Arapura, Donu},
  title         = {Geometry of cohomology support loci for local systems. {I}},
  journal       = {Journal of Algebraic Geometry},
  volume        = {6},
  number        = {3},
  pages         = {563--597},
  year          = {1997}
}

@article{BudurWang2020,
  author        = {Budur, Nero and Wang, Botong},
  title         = {Cohomology jump loci of quasi-compact {K{\"a}hler} manifolds},
  journal       = {Pure and Applied Mathematics Quarterly},
  volume        = {16},
  number        = {4},
  pages         = {981--999},
  year          = {2020}
}

@article{EsnaultKerz2020,
  author        = {Esnault, H{\'e}l{\`e}ne and Kerz, Moritz},
  title         = {Arithmetic subspaces of moduli spaces of rank one local systems},
  journal       = {Cambridge Journal of Mathematics},
  volume        = {8},
  number        = {3},
  pages         = {453--478},
  year          = {2020}
}

@article{DimcaPapadima2014,
  author        = {Dimca, Alexandru and Papadima, {\c{S}}tefan},
  title         = {Non-abelian cohomology jump loci from an analytic viewpoint},
  journal       = {Communications in Contemporary Mathematics},
  volume        = {16},
  number        = {4},
  eid           = {1350025},
  year          = {2014},
  note          = {47 pp.}
}

@article{PapadimaSuciu2020,
  author        = {Papadima, {\c{S}}tefan and Suciu, Alexander I.},
  title         = {Rank two topological and infinitesimal embedded jump loci of quasi-projective manifolds},
  journal       = {Journal of the Institute of Mathematics of Jussieu},
  volume        = {19},
  number        = {2},
  pages         = {451--485},
  year          = {2020}
}

@article{PapadimaSuciu2019,
  author        = {Papadima, {\c{S}}tefan and Suciu, Alexander I.},
  title         = {Naturality properties and comparison results for topological and infinitesimal embedded jump loci},
  journal       = {Advances in Mathematics},
  volume        = {350},
  pages         = {256--303},
  year          = {2019}
}

@misc{lefevre2021mixed,
  author        = {Lef{\`e}vre, Louis-Cl{\'e}ment},
  title         = {Mixed {Hodge} structures on cohomology jump ideals},
  year          = {2021},
  eprint        = {2106.12339},
  archivePrefix = {arXiv},
  primaryClass  = {math.AG}
}

@article{EyssidieuxSimpson2011,
  author        = {Eyssidieux, Philippe and Simpson, Carlos},
  title         = {Variations of mixed {Hodge} structure attached to the deformation theory of a complex variation of {Hodge} structures},
  journal       = {Journal of the European Mathematical Society},
  volume        = {13},
  number        = {6},
  pages         = {1769--1798},
  year          = {2011}
}

@article{Deligne1971,
  author        = {Deligne, Pierre},
  title         = {Th{\'e}orie de {Hodge}: {II}},
  journal       = {Publications Math{\'e}matiques de l'IH{\'E}S},
  volume        = {40},
  pages         = {5--57},
  year          = {1971}
}

@book{Sabbah2018,
  author        = {Sabbah, Claude},
  title         = {Irregular {Hodge} theory},
  series        = {M{\'e}moires de la Soci{\'e}t{\'e} Math{\'e}matique de France (N.S.)},
  number        = {156},
  pagetotal     = {vi+126},
  publisher     = {Soci{\'e}t{\'e} Math{\'e}matique de France},
  address       = {Paris},
  year          = {2018},
  note          = {With the collaboration of Jeng-Daw Yu}
}

@article{BudurWang2015QuasiProjective,
  author        = {Budur, Nero and Wang, Botong},
  title         = {Cohomology jump loci of quasi-projective varieties},
  journal       = {Annales scientifiques de l'{\'E}cole Normale Sup{\'e}rieure},
  series        = {4},
  volume        = {48},
  number        = {1},
  pages         = {227--236},
  year          = {2015}
}
\end{document}